\documentclass{article}
\usepackage{latexsym}
\usepackage{graphics}
\usepackage{amsmath}
\usepackage{amsfonts}
\begin{document}
\bibliographystyle{plain}
\title{The decomposition of global conformal invariants IV:
A proposition on local Riemannian invariants.}
\author{Spyros Alexakis\thanks{University of Toronto, alexakis@math.toronto.edu.
\newline
This work has absorbed
the best part of the author's energy over many years. 
This research was partially conducted during 
the period the author served as a Clay Research Fellow, 
an MSRI postdoctoral fellow,
a Clay Liftoff fellow and a Procter Fellow.  
\newline
The author is immensely indebted to Charles
Fefferman for devoting twelve long months to the meticulous
proof-reading of the present paper. He also wishes to express his
gratitude to the Mathematics Department of Princeton University
for its support during his work on this project.}} \date{}
\maketitle
\newtheorem{proposition}{Proposition}
\newtheorem{theorem}{Theorem}
\newcommand{\Sum}{\sum}
\newtheorem{lemma}{Lemma}
\newtheorem{observation}{Observation}
\newtheorem{formulation}{Formulation}
\newtheorem{definition}{Definition}
\newtheorem{conjecture}{Conjecture}
\newtheorem{corollary}{Corollary}
\numberwithin{equation}{section}
\numberwithin{lemma}{section}
\numberwithin{theorem}{section}
\numberwithin{definition}{section}
\numberwithin{proposition}{section}

\begin{abstract}
This is the fourth in a series of papers where 
we  prove a conjecture of Deser and Schwimmer
regarding the algebraic structure of ``global conformal
invariants''; these  are defined to
be conformally invariant integrals of geometric scalars.
 The conjecture asserts that the integrand of
any such integral can be expressed as a linear
combination of a local conformal invariant, a divergence and of
the Chern-Gauss-Bonnet integrand. 

 The present paper lays out the second half of this entire work: The 
second half proves certain purely algebraic statements regarding 
local Riemannian invariants; these were used 
extensively in \cite{alexakis1,alexakis2}. 
These results may be of independent 
interest, applicable to related problems. 
\end{abstract}
\tableofcontents

\section{Introduction}
 This is the fourth in a series of papers 
\cite{alexakis1,alexakis3,alexakis4,alexakis5,alexakis6} 
where we prove a conjecture of Deser-Schwmimmer \cite{ds:gccaad} regarding the algebraic 
structure of global conformal invariants. We recall that a global 
conformal invariant is an integral of a natural scalar-valued function of 
Riemannian metrics, $\int_{M^n}P(g)dV_g$, with the property that this integral 
remains invariant under conformal re-scalings of the underlying 
metric.\footnote{See the introduction of \cite{alexakis1}
for a detailed discussion of the Deser-Schwimmer 
conjecture, and for background on scalar Riemannian invariants.} 
More precisely, $P(g)$ is assumed to be a linear combination, $P(g)=\sum_{l\in L} a_l C^l(g)$, 
where each $C^l(g)$ is a complete contraction in the form:

\begin{equation}
\label{contraction} 
contr^l(\nabla^{(m_1)}R\otimes\dots\otimes\nabla^{(m_s)}R);
\end{equation}
 here each factor $\nabla^{(m)}R$ stands for the $m^{th}$ iterated 
covariant derivative of the curvature tensor $R$. $\nabla$ is the Levi-Civita 
connection of the metric $g$ and $R$ is the curvature associated to this connection. 
The contractions are taken with respect to the quadratic form $g^{ij}$.
In this series of papers we prove:

\begin{theorem}
\label{thetheorem} 
Assume that $P(g)=\sum_{l\in L} a_l C^l(g)$, where each $C^l(g)$ is a 
complete contraction in the form (\ref{contraction}), with weight $-n$. 
Assume that for every closed Riemannian manifold $(M^n,g)$ and every $\phi\in C^\infty (M^n)$:
$$\int_{M^n}P(e^{2\phi}g)dV_{e^{2\phi}g}=\int_{M^n}P(g)dV_g.$$

We claim that $P(g)$ can then be expressed in the form:
$$P(g)=W(g)+div_iT^i(g)+\operatorname{Pfaff}(R_{ijkl}).$$
Here $W(g)$ stands for a local conformal invariant of weight $-n$  (meaning 
that $W(e^{2\phi}g)=e^{-n\phi}W(g)$ for every $\phi\in C^\infty (M^n)$), 
$div_iT^i(g)$ is the divergence of a Riemannian vector field of 
weight $-n+1$, and   $\operatorname{Pfaff}(R_{ijkl})$ is the Pfaffian of the curvature tensor. 
\end{theorem}

We now discuss the position of the present paper in this series. 

 We recall from the introduction of \cite{alexakis1} that this series of papers can be natually 
 subdivided into two parts: Part I (consisting of 
\cite{alexakis1,alexakis2,alexakis3}) proves the Deser-Schwimmer conjecture, 
 subject to establishing certain ``Main algebraic Propositions'', 
namely Proposition 5.2 in \cite{alexakis1} and Propositions 3.1, 3.2 in \cite{alexakis2}. 
Part II, consisting of the present paper and 
\cite{alexakis5,alexakis6} prove 
these main algebraic propositions. 
\newline

 \par The first task that we undertake in the present paper 
is to reduce the ``main algebraic Propositions'' in \cite{alexakis1,alexakis2}
 to a single Proposition, \ref{giade} below, which we 
call the ``fundamental Proposition \ref{giade}'' and which will 
be proven by an elaborate induction on four parameters.  
 This fundamental Proposition is actually a 
 {\it generalization} of the ``Main algebraic Propositions'' 
in \cite{alexakis1,alexakis2},
in the sense that the ``Main algebraic Propositions'' 
are special cases of Proposition \ref{giade}; in fact 
they are the ultimate or penultimate steps of the aforementioned  
induction with respect to certain of the parameters.  
\newline

{\bf An outline of the goals of the papers 
\cite{alexakis4,alexakis5,alexakis6}:} 
The main goal in the present paper is to introduce Proposition \ref{giade} below, 
which will imply the ``main algebraic Propositions'' 
in \cite{alexakis1,alexakis2}, and then to reduce {\it the inductive 
step} in the proof of Proposition \ref{giade} 
to three Lemmas, \ref{zetajones}, \ref{pool2}, \ref{pskovb} below (along with two preparatory claims
needed for Lemma \ref{pskovb}, namely Lemmas \ref{oui}, \ref{oui2}):
We prove in the present paper that the three Lemmas \ref{zetajones}, \ref{pool2}, \ref{pskovb} 
imply the inductive step of Proposition \ref{giade}, apart from certain special cases. 
In this derivation we employ certain technical
 Lemmas.\footnote{(More on this in the outline of the present paper below).} 
In the next paper in the series, \cite{alexakis5} we derive the 
inductive step of Proposition \ref{giade} in the special cases, and we also 
provide a proof of the  aforementioned technical Lemmas. Thus, 
the present paper and \cite{alexakis5} reduce the task of  proving 
the Deser-Schwimmer conjecture to 
proving the Lemmas \ref{zetajones}--\ref{pskovb} below.

Then, Lemmas \ref{zetajones}--\ref{pskovb}  are proven in 
the final paper in this series, \cite{alexakis6}. 
\newline

{\bf Outline of the ``Fundamental Proposition \ref{giade}'':} 
In section \ref{fundprop} we set up the considerable notational and language
conventions needed to state our fundamental
 Proposition \ref{giade}; we then state the  fundamental 
Proposition and explain how the ``main algebraic Propositions''
 5.2 and 3.1, 3.2  in \cite{alexakis1,alexakis2} are special cases of it. 
We also explain how the fundamental Proposition will be
proven via an induction on four parameters. 
 In section \ref{startB} we 
distinguish three cases I,II,III on the 
hypothesis of Proposition \ref{giade}  and 
 claim  three Lemmas, \ref{zetajones}, \ref{pool2}, 
\ref{pskovb} which correspond to these three cases. 
 Finally in section \ref{laderivation} we prove that these three Lemmas {\it imply} Proposition 
 \ref{giade}.\footnote{We make use of the inductive assumption 
of Proposition \ref{giade} in this derivation.}
 In section \ref{laderivation} we also assert certain important technical 
Lemmas which will also be used in the subequent papers in this series;
some of these technical papers are proven in the 
present paper and some in \cite{alexakis5}. 
 \newline

Now, since the fundamental proposition is very complicated to even write out, 
we reproduce here the claim of the first ``main algebraic Proposition''
from \cite{alexakis1}, for the reader's convenience. As explained above, 
this first ``main algebraic Proposition'' 
is a special case of Proposition \ref{giade} below.
\newline

{\it A simplified description of the main algebraic Proposition 5.2 in \cite{alexakis1}:}
Given a Riemannian metric $g$ over an $n$-dimensional 
 manifold $M^n$ and auxilliary $C^\infty$ scalar-valued functions 
$\Omega_1,\dots,\Omega_p$ defined over $M^n$, the 
objects of study are linear combinations of tensor fields 
$\sum_{l\in L} a_l C^{l,i_1\dots i_\alpha}_g$, where each 
$C^{l,i_1\dots i_\alpha}_g$ is a {\it partial contraction} 
with $\alpha$ free indices, in the form: 
\begin{equation}
\label{gen.form1}
 pcontr(\nabla^{(m)}R\otimes\dots\otimes\nabla^{(m_s)}R\otimes
\nabla^{(b_1)}\Omega_1\otimes\dots\otimes\nabla^{(b_m)}\Omega_p);
\end{equation}
here $\nabla^{(m)}R$ stands for the $m^{th}$ covariant derivative of 
the curvature tensor $R$,\footnote{In other words it is an $(m+4)$-tensor; 
 if we write out its free 
indices it would be in the form 
$\nabla^{(m)}_{r_1\dots r_m}R_{ijkl}$.} and 
$\nabla^{(b)}\Omega_h$ stands for the $b^{th}$ covariant 
derivative of the function $\Omega_h$. A {\it partial contraction} 
means that we have list of pairs of indices $({}_a,{}_b),\dots, ({}_c,{}_d)$ 
in (\ref{gen.form1}), 
which are contracted against each other using 
the  metric $g^{ij}$. The remaining 
 indices (which are not contracted against 
another index in (\ref{gen.form1})) are the 
{\it free indices} ${}_{i_1},\dots, {}_{i_\alpha}$. 

The ``main algebraic Proposition'' of \cite{alexakis1} (roughly) asserts the following: Let 
$\sum_{l\in L_\mu} a_l C^{l,i_1\dots i_\mu}_g$ stand for a 
linear combination of partial contractions in the form (\ref{gen.form1}), 
where each $C^{l,i_1\dots i_\mu}_g$ has a given number $\sigma_1$ of factors and a given number 
$p$ of factor $\nabla^{(b)}\Omega_h$. Assume also that $\sigma_1+p\ge 3$, 
each $b_i\ge 2$,\footnote{This means that 
each function $\Omega_h$ is differentiated at least twice.} and that 
for each pair of contracting 
indices $({}_a,{}_b)$ in any given $C^{l,i_1\dots i_\mu}_g$, 
the indices ${}_a,{}_b$ do not belong to the 
same factor. Assume also the rank 
$\mu>0$ is fixed and each partial contraction $C^{l,i_1\dots i_\mu}_g, l\in L_\mu$ 
has a given {\it weight} $-n+\mu$.\footnote{See \cite{alexakis1} 
for a precise definition of weight.} 
Let also $\sum_{l\in L_{>\mu}} a_l C^{l,i_1\dots i_{y_l}}_g$ stand for 
a (formal) linear combination of partial contractions of 
weight $-n+y_l$, with all the properties of the 
terms indexed in $L_\mu$, {\it except} that now all the partial 
contractions have a different rank $y_l$, and each $y_l>\mu$. 

The assumption of the ``main algebraic Proposition'' is 
a local equation:

\begin{equation}
\label{assumption.1} 
\sum_{l\in L_\mu} a_l Xdiv_{i_1}\dots Xdiv_{i_\mu}C^{l,i_1\dots i_\mu}_g+
\sum_{l\in L_{>\mu}} a_l Xdiv_{i_1}\dots Xdiv_{i_{y_l}}C^{l,i_1\dots i_{y_l}}_g=0,
\end{equation}
which
is assumed to hold {\it modulo} complete contractions with $\sigma+1$ factors.
Here given a partial contraction $C^{l,i_1\dots i_\alpha}_g$ 
in the form (\ref{gen.form1}) $Xdiv_{i_s}[C^{l,i_1\dots i_\alpha}_g]$
stands for sum of $\sigma-1$ terms in $div_{i_s}[C^{l,i_1\dots i_\alpha}_g]$
where the derivative $\nabla^{i_s}$ is {\it not} 
allowed to hit the factor to which the 
free index ${}_{i_s}$ belongs.\footnote{Recall that given a partial contraction
$C^{l,i_1\dots i_\alpha}_g$ in the form (\ref{gen.form1}) 
with $\sigma$ factors, $div_{i_s}C^{l,i_1\dots i_\alpha}_g$ 
is a sum of $\sigma$ partial contractions of rank $\alpha-1$. 
The first summand arises by adding a derivative $\nabla^{i_s}$ 
onto the first factor $T_1$ and then contracting the upper index ${}^{i_s}$
against the free index ${}_{i_s}$; the second summand 
arises by adding a derivative $\nabla^{i_s}$ 
onto the second factor $T_2$ and then contracting the upper index ${}^{i_s}$
against the free index ${}_{i_s}$ etc.} 

The main algebraic Proposition in \cite{alexakis1} 
then claims that there will exist a linear combination of partial 
contactions in the form (\ref{gen.form1}), 
$\sum_{h\in H} a_h C^{h,i_1\dots i_{\mu+1}}_g$ with all the properties 
of the terms indexed in $L_{>\mu}$, and all with rank $(\mu+1)$, so that:

\begin{equation}
\label{conclusion.1} 
\sum_{l\in L_\mu} a_l C^{l,(i_1\dots i_\mu)}_g+
\sum_{h\in H} a_h  Xdiv_{i_{\mu+1}}
C^{l,(i_1\dots i_\mu)i_{\mu+1}}_g=0;
\end{equation}
the above holds modulo terms of length $\sigma+1$. Also the symbol $(\dots)$ 
means that we are {\it symmetrizing} over the indices between parentheses. 
\newline

{\bf  The proof of the ``fundamental Proposition \ref{giade}'' 
via an induction, and a brief description of Lemmas \ref{zetajones},\ref{pool2},\ref{pskovb}:} 

The fundamental Proposition \ref{giade} is a generalization 
of Proposition 5.2 from \cite{alexakis1}, 
in the sense that it deals with partial contractions 
in the form (\ref{gen.form1}), which {\it in addition} 
contain factors $\nabla\phi_h$;\footnote{See the forms (\ref{form1}), (\ref{form2}) below.}
 these are assumed to contract 
against the different factors  $\nabla^{(m)}R, \nabla^{(p)}\Omega_x$ 
according to a given {\it pattern}.\footnote{This encoding is described by 
the notions of {\it weak} and {\it simple} character--see the informal 
discussion after Definition \ref{raellenth}.} The Proposition \ref{giade} 
also groups up the different partial contractions of minimum rank indexxed in  $L_\mu$ 
according to the distribution of the free indices among its different factors 
$\nabla^{(m)}R, \nabla^{(p)}\Omega_x$.\footnote{This encoding is described by 
the notions of {\it double} and {\it refined double} character--see the informal 
discussion after Definition \ref{raellenth}.} The claim of Proposition \ref{giade}
is then an adaptation of (\ref{conclusion.1}), restricted to a particular 
subset of the partial contractions indexed in $L_\mu$. A discussion of 
how the Proposition \ref{giade} is proven via an induction on four
parameters, as well as how the inductive step is reduced to the 
Lemmas \ref{zetajones}, \ref{pool2}, \ref{pskovb} is provided in subsections 
\ref{sectioninduction} and 
 \ref{reduction.lemmas}. The reader is also refered to those 
subsections for a more conceptual outline of the ideas 
in the present paper. 
\newline

Before proceeding to give the strict formulation of the 
fundamental Proposition, we 
digress to discuss the relationship of the whole 
work \cite{alexakis1}--\cite{alexakis6} and of the 
papers \cite{alexakis4}--\cite{alexakis6} in particular with the 
study of local scalar Riemannian and conformal invariants.

{\bf Broad Discussion:} The theory of {\it local} invariants of Riemannian structures 
(and indeed, of more general geometries,
e.g.~conformal, projective, or CR)  has a long history. 
As discussed in \cite{alexakis1}, the original foundations of this 
field were laid in the work of Hermann Weyl and \'Elie Cartan, see \cite{w:cg, cartan}. 
The task of writing out local invariants of a given geometry is intimately connected
with understanding polynomials in a space of tensors with  given symmetries;  
these polynomials are required to remain invariant under the action of a Lie group
on the components of the tensors. 
In particular, the problem of writing down all 
 local Riemannian invariants reduces to understanding 
the invariants of the orthogonal group. 

 In more recent times, a major program was laid out by C.~Fefferman in \cite{f:ma}
aimed at finding all scalar local invariants in CR geometry. This was motivated 
by the problem of understanding the  
local invariants which appear in the asymptotic expansion of the 
Bergman and Szeg\"o kernels of strictly pseudo-convex CR manifolds,
 in a similar way to which Riemannian invariants appear in the asymptotic expansion  
of the heat kernel; the study of the local invariants
in the singularities of these kernels led to important breakthroughs 
in \cite{beg:itccg} and more recently by Hirachi in \cite{hirachi1}.
 This program was later extended  to conformal geometry in \cite{fg:ci}. 
Both these geometries belong to a 
broader class of structures, the
{\it parabolic geometries}; these admit a principal bundle with 
structure group a parabolic subgroup $P$ of a semi-simple 
Lie group $G$, and a Cartan connection on that principle bundle 
(see the introduction in \cite{cg1}). 
An important question in the study of these structures 
is the problem of constructing all their local invariants, which 
can be thought of as the {\it natural, intrinsic} scalars of these structures.

  In the context of conformal geometry, the first (modern) landmark 
in understanding {\it local conformal invariants} was the work of Fefferman 
and Graham in 1985 \cite{fg:ci},
where they introduced the {\it ambient metric}. This allows one to 
construct local conformal invariants of any order in odd 
dimensions, and up to order $\frac{n}{2}$ in even dimensions. 
The question is then whether {\it all} invariants arise via this construction. 

The subsequent work of Bailey-Eastwood-Graham \cite{beg:itccg} proved that 
this is indeed true in odd dimensions; in even dimensions, 
they proved that the result holds  
when the weight (in absolute value) is bounded by the dimension. The ambient metric construction 
in even dimensions was recently extended by Graham-Hirachi, \cite{grhir}; this enables them to 
identify in a satisfactory way {\it all} local conformal invariants, 
even when the weight (in absolute value) exceeds the dimension.  

 An alternative 
construction of local conformal invariants can be obtained via the {\it tractor calculus} 
introduced by Bailey-Eastwood-Gover in \cite{bego}. This construction bears a strong 
resemblance to the Cartan conformal connection, and to 
the work of T.Y.~Thomas in 1934, \cite{thomas}. The tractor 
calculus has proven to be very universal; 
tractor buncles have been constructed \cite{cg1} for an entire class of parabolic geometries. 
The relation betweeen the conformal tractor calculus and the Fefferman-Graham 
ambient metric  has been elucidated in \cite{cg2}.

The present work \cite{alexakis1}--\cite{alexakis6}, while pertaining to the question above
(given that it ultimately deals with the algebraic form of local 
{\it Riemannian} and {\it conformal} invariants), nonetheless addresses a different 
{\it type} of problem:  We here consider Riemannian invariants $P(g)$ for 
which the {\it integral} $\int_{M^n}P(g)dV_g$ remains invariant 
under conformal changes of the underlying metric; we then seek to understand 
the possible algebraic form of the {\it integrand} $P(g)$, 
ultimately proving that it can be de-composed 
in the way that Deser and Schwimmer asserted. 
It is thus not surprising that the prior work on 
 the construction and understanding of local {\it conformal} 
invariants plays a central role in this endeavor, in 
the papers \cite{alexakis2,alexakis3}.

On the other hand, as explained above, 
a central element of our proof is the (roughly outlined above)
``fundamental Proposition \ref{giade}'',\footnote{This 
Proposition is a generalization of the main algebraic 
Propositions 5.1, 3.1, 3.2 in \cite{alexakis1,alexakis2}.} which deals {\it exclusively} 
with algebraic properties of the {\it classical} scalar Riemannian invariants.\footnote{We 
refer the reader to the introuction of \cite{alexakis1} for a detailed discussion of these.}
The ``fundamental Proposition \ref{giade}'' makes no reference to 
integration; it is purely a statement concerning 
algebraic properties of 
{\it local Riemannian invariants}. 
While the author was led to led to the main algebraic 
Propositions in \cite{alexakis1,alexakis2} 
out of the strategy that he felt was necessary to 
solve the Deser-Schwimmer conjecture, they can 
be thought of as results with an independent interest. 
The {\it proof} of these Propositions, presented
 in the present paper and in \cite{alexakis5,alexakis6} is in fact 
not particularly intuitive. It is the author's 
sincere hope that deeper insight (and hopefuly a more intuitive proof) 
will be obtained in the future as to  {\it why} these algebraic 
 Propositions hold.

\section{The fundamental Proposition.}
\label{fundprop}

\par In order to state and prove the fundamental proposition we will
 need to introduce a lot of terminology. 
 
\subsection{Definitions and Terminology.}
 \par We will be considering (complete or partial) contractions
$C^{i_1\dots i_\alpha}_{g}(\Omega_1,\dots , \Omega_p,\phi_1,\dots
,\phi_u)$ of length $\sigma +u$ (with no internal contractions) in
the form:

\begin{equation}
\label{form1}
\begin{split}
&pcontr(\nabla^{(m_1)}R_{ijkl}\otimes\dots\otimes\nabla^{(m_s)}R_{ijkl}
\otimes
\\& \nabla^{(b_1)}\Omega_1\otimes\dots\otimes \nabla^{(b_p)}\Omega_p
\otimes\nabla\phi_1\otimes\dots \otimes\nabla\phi_u);
\end{split}
\end{equation}
here $\sigma=s+p$ and ${}_{i_1},\dots ,{}_{i_\alpha}$ are the free indices in 
$C^{i_1\dots i_\alpha}_{g}(\Omega_1,\dots , \Omega_p,\phi_1,\dots
,\phi_u)$.

\begin{definition}
\label{accept} Any (complete or partial) contraction in the form (\ref{form1})
will be called acceptable if: 

\begin{enumerate}

 \item  Each of the free indices must 
belong to a factor $\nabla^{(m)}R_{ijkl}$ or $\nabla^{(b)}\Omega_h$.

\item Each of the factors
 $\nabla\phi_h$ is contracting against a factor $\nabla^{(m)}R_{ijkl}$ or
$\nabla^{(a)}\Omega_f$. 

\item Each of the factors
$\nabla^{(a)}\Omega_f$ should satisfy $a\ge 2$. 
\end{enumerate}
\end{definition}

\par More generally, we will also be  considering tensor fields
\\ $C^{i_1\dots i_\alpha}_{g} (\Omega_1,\dots ,
\Omega_p,\phi_{z_1},\dots ,\phi_{z_u},\phi'_{z_{u+1}},\dots ,
\phi'_{z_{u+d}},\tilde{\phi}_{z_{u+d+1}}, \dots
,\tilde{\phi}_{z_{u+d+y}})$ of length $\sigma+u$ (with no internal
contractions)  in the form:\footnote{We recall that 
 $S_{*}\nabla^{(\nu)}_{r_1\dots r_\nu}R_{ijkl}$ stands for the symmetrization of the tensor 
$\nabla^{(\nu)}_{r_1\dots r_\nu}R_{ijkl}$ over the indices ${}_{r_1},\dots,{}_{r_\nu},{}_j$.}

\begin{equation}
\label{form2}
\begin{split}
&pcontr(\nabla^{(m_1)}R_{ijkl}\otimes\dots\otimes\nabla^{(m_{\sigma_1})}
R_{ijkl}\otimes \\&S_{*}\nabla^{(\nu_1)}R_{ijkl}\otimes\dots\otimes
S_{*}\nabla^{(\nu_t)} R_{ijkl}\otimes
\\& \nabla^{(b_1)}\Omega_1\otimes\dots\otimes
\nabla^{(b_p)}\Omega_p\otimes
\\& \nabla\phi_{z_1}\dots \otimes\nabla\phi_{z_w}\otimes\nabla
\phi'_{z_{w+1}}\otimes
\dots\otimes\nabla\phi'_{z_{w+d}}\otimes\dots \otimes
\nabla\tilde{\phi}_{z_{w+d+1}}\otimes\dots\otimes\nabla\tilde{\phi}_{z_{w+d+y}}).
\end{split}
\end{equation}
Here,   
the functions
$\tilde{\phi}_a$, $\phi'_b$ are the same as the functions
$\phi_a,\phi_b$. The symbols $\tilde{}$ and $'$ are used only to
illustrate the {\it kind} of indices that these factors are contracting
against (we will explain this in the next definition).

\par The notion of {\it acceptability} for contractions in the form (\ref{form2}) is a generalization of Definition
\ref{accept}:

\begin{definition}
\label{accept2} We will call a (complete or partial) contraction in the form
(\ref{form2}) acceptable if the following conditions hold:

\begin{enumerate}
\item{$\{z_1,\dots , z_{w+d+y}\} =\{1,\dots ,u\}$. Also, each of the
free indices must belong to a factor $\nabla^{(m)}R_{ijkl}$,
 $\nabla^{(b)}\Omega_h$ or $S_{*}\nabla^{(\nu)}R_{ijkl}$. In addition, the factors
 $\nabla^{(b)}\Omega_f$ must have $b\ge 2$.} 
 
\item{The factors $\nabla\phi_h,\nabla\phi'_h,\nabla\tilde{\phi}_h$ contract according to the 
following pattern:  Each of the factors
 $\nabla\phi_h$ is contracting against a
 derivative index in a factor $\nabla^{(m)}R_{ijkl}$ or $\nabla^{(b)}
\Omega_f$. Each of the factors $\tilde{\phi}_h$ must be 
contracting against the index ${}_i$ of some factor
$S_{*}\nabla^{(\nu)}R_{ijkl}$. Conversely, each index ${}_i$ in
any factor $S_{*}\nabla^{(\nu)}R_{ijkl}$ must contract against
some factor $\nabla\tilde{\phi}_h$. Lastly, each factor
$\nabla\phi'_h$ is contracting against some factor
$S_{*}\nabla^{(\nu)}_{r_1\dots r_\nu} R_{ijkl}$, but necessarily
against one of the indices ${}_{r_1},\dots ,{}_{r_\nu},{}_j$.}
\end{enumerate}
\end{definition}

\begin{definition}
\label{raellenth} For any (complete or partial) contraction in the form (\ref{form1})
or (\ref{form2}), we define its {\it real length} to be the number
of its factors {\it if we exclude} the factors $\nabla\phi_h$,
$\nabla\tilde{\phi}_h$, $\nabla\phi'_h$. (So for contractions in the form (\ref{form2}) the real length is $\sigma_1+t+p$).
\end{definition}

\par We now introduce the notions of {\it weak}, {\it simple}, {\it
double} and {\it refined-double}  characters for acceptable 
contractions $C^{i_1\dots i_a}_g$ in the form (\ref{form2}). Since these definitions are
rather technical, we give the gist of these notions here: The weak
character encodes the pattern of {\it which} factors in $C^{i_1\dots i_a}_g$ the various
terms $\nabla\phi_h,h=1,\dots ,u$ are contracting against. The
simple character encodes the above, but also encodes whether each
given factor $\nabla\phi_h$ that contracts against a factor $T=
 S_{*}\nabla^{(\nu)}R_{ijkl}$ is contracting against  the
index ${}_i$, or one of the indices ${}_{r_1},\dots ,{}_{r_\nu},{}_j$. 
The double character encodes the simple character, but
also encodes how the free indices are distributed among the different 
factors in $C^{i_1\dots i_a}_g$ (i.e.~how many free indices belong to each factor).
Finally, the refined double character encodes the same information
as the double character, but also encodes whether the free indices
are {\it special indices}.\footnote{Meaning that they are internal indices
in some $\nabla^{(m)}R_{ijkl}$ or indices ${}_k,{}_l$ in some
$S_{*}\nabla^{(\nu)}R_{ijkl}$.}

\par Now we present the proper definitions of the different notions of ``character''. 

\begin{definition}
\label{charf1} Consider any acceptable complete or partial
contraction in
 either of the forms (\ref{form1}) or (\ref{form2}).
The weak character $\vec{\kappa}_{weak}$ is defined
to be a pair of two lists of sets: $(L_1,L_2)$. 
$L_1$ stands for the list $(S_1,\dots , S_p)$ where each
$S_t$ stands for the set of numbers $r$ for which $\nabla\phi_r$
contracts against the factor $\nabla^{(b_t)}\Omega_t$. 
$L_2$ stands for the list  of sets $(S_1,\dots ,S_{\sigma-p})$
where each $S_t$ stands for the set of numbers $r$ for which
$\nabla\phi_r$ contracts against the $t^{th}$ curvature factor in
(\ref{form1}) or (\ref{form2}) (in the latter case the curvature factor may be in the form
$\nabla^{(m)}R_{ijkl}$ or $S_{*}\nabla^{(\nu)}R_{ijkl}$).
\end{definition}

\begin{definition}
\label{kenya} Consider complete or partial contractions in the form (\ref{form2}); we
define the simple character $\vec{\kappa}_{simp}$ to be a triplet of
 lists: $(L_1,L_2,L_3)$.

\par $L_1$ is as above. $L_2$
stands for the list of of sets $(S_1,\dots ,S_{\sigma_1})$ where
each $S_t$ stands for the set of numbers $r$ for which
$\nabla\phi_r$ contracts against the  $t^{th}$ factor
$\nabla^{(m_t)}R_{ijkl}$ in the first line of (\ref{form2}). $L_3$
is a sequence of pairs of sets: $L_3=[(\{\alpha_1\}, S_1), \dots
,(\{\alpha_t\},S_t)]$, where $\alpha_w$ stands for the one number
for which the index ${}_i$ in the $w^{th}$ factor in the second line
of (\ref{form2}) is contracting
 against the factor $\nabla\tilde{\phi}_{a_w}$. $S_w$ stands for the set
 of numbers $r$ for which the $w^{th}$ factor is contracting against the factors
$\nabla\phi'_r$ in (\ref{form2}).
\end{definition}

\par Now, we define the double character. We note that this notion
is defined for tensor fields in  the form (\ref{form2}) that {\it
do not} have both indices ${}_i,{}_j$ or ${}_k,{}_l$ in a factor $\nabla^{(m)}R_{ijkl}$ 
or $S_{*}\nabla^{(\nu)}R_{ijkl}$ being free.

\begin{definition}
\label{kenya2}
Consider complete or partial contractions in the form (\ref{form2}); we
define the  double character to be the  union of two triplets of lists:
$\vec{\kappa}_{doub}=(L_1,L_2,L_3)|(H_1,H_2,H_3)$. Here
$L_1,L_2,L_3$ are as above. $H_1$ stands for the list
$(h_1,\dots ,h_p)$, where $h_t$ stands for the number of free
indices that belong to the factor $\nabla^{(b_t)}\Omega_t$. $H_2$
stands for the list of numbers $(h_1,\dots ,h_{\sigma_1})$, where
each $h_i$ stands for the number of free indices that belong to
the $i^{th}$ factor in the first line of (\ref{form2}). $H_3$
stands for the set of
 numbers $(h_1,\dots ,h_t)$ where $h_u$ stands for the number of free
 indices that belong to the $u^{th}$ factor on the second line of
(\ref{form2}).
\end{definition}

\par Now, one more definition before stating our main  Proposition for the present paper.
 We will be defining the {\it refined} double
 character of tensor fields in the form (\ref{form2}).
For that purpose, we will
 be paying special attention to the free indices ${}_{i_f}$
that are internal indices in some factor $\nabla^{(m)}R_{ijkl}$ or are
one of the indices ${}_k,{}_l$ in one of the factors
$S_{*}\nabla^{(\nu)}R_{ijkl}$. We will be calling those free indices
{\it special free indices}.

\begin{definition}
\label{refdoub}
Consider complete or partial contractions in the form (\ref{form2}); we
define  its {\it refined double
character} to be the union of two triplets of 
 sets: $\vec{\kappa}_{ref-doub}=(L_1,L_2,L_3)|
(H_1,\tilde{H}_2,\tilde{H}_3)$ where the sets $L_1,L_2, L_3,H_1$ are as before,
whereas:

\par $\tilde{H}_2$ stands for the list of sets $(\tilde{h}_1,\dots
\tilde{h}_{\sigma_1})$ where $\tilde{h}_k$ stands for the
following: If the $k^{th}$ factor $\nabla^{(m)}R_{ijkl}$ has no
special free indices then $\tilde{h}_k=h_k$ (same as for the double character).
 If it contains one special free index then
$\tilde{h}_k=\{h_k\}\bigcup \{*\}$. Finally, if it contains
 two special free indices then
$\tilde{h}_k=\{h_k\}\bigcup \{**\}$.

\par $\tilde{H}_3$ stands for the list of sets $(\tilde{h}_1,\dots
\tilde{h}_{\sigma_2})$ where $\tilde{h}_k$ stands for the
following: If the $k^{th}$ factor $S_{*}\nabla^{(\nu)}R_{ijkl}$ has
no special free indices then $\tilde{h}_k=h_k$ (same as for the double character).
 If it contains one special free index then
$\tilde{h}_k=\{h_k\}\bigcup \{*\}$.

\par The elements $\{**\}, \{*\}$ above will be called marks.
\end{definition}

\par Now, we will introduce a notion of {\it equivalence}
 for the characters (weak, simple, double or refined double)
of tensor fields.

\begin{definition}
\label{equivchar}
 We say that two (complete or partial) contractions in the form (\ref{form2})
have equivalent (simple, double or refined double)
characters if
their (simple, double or refined double) characters can be made equal by permuting factors
among
 the first two lines of (\ref{form2}).

\par More generally, we will say that two (complete or partial) contractions in the more general
form (\ref{form1}) (possibly of different rank)
have equivalent weak characters if we
can permute their curvature factors and make their weak characters equal.
\end{definition}

\par We thus see that the various ``characters'' we have defined can be thought of as abstract lists,
which are equipped with a natural notion of equivalence. We note that we will
occasionally be speaking of a
(weak, simple, double or refined double) character abstractly, without specifying
a (complete or partial) contraction or tensor field that it represents. Furthermore,
we have seen that the notions of weak, simple, double and refined double
characters are {\it graded}, in the sense that 
each of these four notions contains all the
information of the previous ones.  Now, given a 
simple character $\vec{\kappa}_{simp}$ we define
$Weak(\vec{\kappa}_{simp})$ to stand for the weak 
character that corresponds to that simple character.
Analogously, given any refined double character 
$\vec{\kappa}_{ref-doub}$, we let $Simp(\vec{\kappa}_{ref-doub})$
stand for the simple character that 
corresponds to that refined double character and also
$Weak(\vec{\kappa}_{ref-doub})$ to
stand for the weak character that corresponds to that refined double character.
\newline

\par Now, we will introduce a weak notion of {\it ordering} for simple and refined double characters.
This notion is ``weak'' in the  sense that we will be
specifying a particular simple or refined double character
$\vec{\kappa}_{simp}$ and $\vec{\kappa}_{ref-doub}$ respectively,
and we will define what it means for complete or partial contractions 
 (in the form (\ref{form1}) or (\ref{form2})) to be
{\it subsequent} to $\vec{\kappa}_{simp}$ and
$\vec{\kappa}_{ref-doub}$ respectively. This relation is {\it not}
transitive.

\begin{definition}
\label{kenya3} Given any contraction in the form (\ref{form2}), we
consider any simple character $\vec{\kappa}_{simp}$ or refined
double character $\vec{\kappa}_{ref-doub}$ and we let the
defining set $Def(\vec{\kappa}_{simp})$,
$Def(\vec{\kappa}_{ref-doub})$ to be the set of numbers $r$ for
which $\nabla\tilde{\phi}_r$ is contracting against an internal
index ${}_i$ in some factor $S_{*}\nabla^{(\nu)}R_{ijkl}$.

We now consider any general complete contraction
$C_{g}(\Omega_1,\dots ,\Omega_p,\phi_1,\dots
,\phi_u,\\\phi'_{u+1},\dots ,\phi'_{u+d},
\tilde{\phi}_{u+d+1},\dots ,\tilde{\phi}_{u+d+y})$ or partial contraction
$C_{g}^{i_1\dots i_a}$  in the form (\ref{form1}) or
(\ref{form2}), with a weak character $Weak(\vec{\kappa}_{simp})$
or $Weak(\vec{\kappa}_{ref-doub})$ respectively. We say that
$C_{g}$ or $C_{g}^{i_1\dots i_a}$ is simply subsequent to
$\vec{\kappa}_{simp}$ or $\vec{\kappa}_{ref-doub}$ if for at least
one number $\nu\in Def(\vec{\kappa}_{simp})$ (or $\nu\in
Def(\vec{\kappa}_{ref-doub})$), the factor
$\nabla\tilde{\phi}_\nu$ in $C_{g}$ or $C_{g}^{i_1\dots i_a}$ is
contracting against a derivative index. This terminology extends
to linear combinations.
\end{definition}

\par Now, we will introduce 
 a partial ordering among refined double
characters $\vec{\kappa}_h$ with
$Simp(\vec{\kappa}_h)=\vec{\kappa}_{simp}$, where $\vec{\kappa}_{simp}$ is a {\it fixed}
simple character.
 To do this, some more notation is needed:

\par For a given tensor field in the form (\ref{form2}), with a
refined double character $\vec{\kappa}_{ref-doub}$, we define
$Def^{*}(\vec{\kappa}_{ref-doub})$ to stand for the subset of
$Def(\vec{\kappa}_{ref-doub})$ which consists of those numbers
$a_w$ for which $\nabla\tilde{\phi}_{a_w}$ contracts
against a factor $S^{*}\nabla^{(\nu)}R_{ijkl}$ where one of the indices
${}_k$ or ${}_l$ is a free index.

\begin{definition}
\label{doubsub}
 We  compare two refined double characters 
$\vec{\kappa}_1,\vec{\kappa}_2$ with \\$Simp(\vec{\kappa}_1)=Simp(\vec{\kappa}_2)$  according
to their $*$-decreasing rearrangements, which means the following:

\par By a $*$-decreasing rearrangement of any refined double character, we mean the
 rearrangement of the lists $\tilde{H}_2,\tilde{H}_3$ (see Definition \ref{refdoub} above) so
that the elements in $\tilde{H}_2$ with a mark $\{**\}$ must come
first (and the elements with such factors are arranged in
decreasing rearrangement). Then among the rest of the elements,
 the ones with a mark $\{*\}$ must come first
(and the elements with such a mark are arranged in decreasing
rearrangement). Then, the elements without a mark will come in the
end, arranged in decreasing rearrangement. Furthermore, for the
lists $\tilde{H}_3$ $*$-decreasing rearrangement means that the
elements in $\tilde{H}_3$
 corresponding to a factor $S_{*}\nabla^{(\nu)}R_{ijkl}$ {\it which
 is contracting against some factor $\nabla\phi'_h$ in $\vec{\kappa}_{simp}$}
 will come first in the list (and they are arranged
 in decreasing rearrangement), and then come the elements corresponding
to a factor $S_{*}\nabla^{(\nu)}R_{ijkl}$ {\it which
 are not  contracting against any factor $\nabla\phi'_h$
in $\vec{\kappa}_{simp}$}, and those are also arranged in decreasing rearrangement.

\par Now, to compare the refined double characters
$\vec{\kappa}_1,\vec{\kappa}_2$ according to their $*$-decreasing rearrangements
means that we
 take their lists $\tilde{H}_3$ in $*$-decreasing rearrangement (see above) and
 order them lexicographically according to these rearranged lists.
 If they are still equivalent, we take the lists
$\tilde{H}_2$ in
 $*$-decreasing rearrangement and order them lexicographically. If
 they are still equivalent, we take the decreasing rearrangement
 of the lists $H_1$ and compare them lexicographically. ``Lexicographically''
  here means that we compare the first elements in the $*$-rearranged lists, then the second etc.
The list with the first larger element is precedent (the converse
of ``subsequent'') to the other list. 

\par If $\vec{\kappa}_1,\vec{\kappa}_2$ 
 are still equivalent after all these comparisons, we say that the refined 
 double characters $\vec{\kappa}_1,\vec{\kappa}_2$ are ``equipolent''. We remark that two 
 refined double characters $\vec{\kappa}_1,\vec{\kappa}_2$ can be equipolent 
 withount being the same.
 \end{definition}

\par A final note before stating our proposition. We remark that in the
definition of simple or of weak character, the {\it rank} of the tensor field
 plays no role. On the other hand, in the definitions of double character
 and of refined double character, the rank of the tensor fields does
 play a role: Two double characters or two refined double characters
 cannot be equivalent if the tensor fields do not have the same rank.
We will then extend the notion of double character and refined
double character as follows: We consider any tensor field
$C^{i_1\dots i_\beta}$ in the form (\ref{form2}), and also any
number $\alpha\le \beta$. We then define
 the $\alpha$-double character or the $\alpha$-refined double character
of $C^{i_1\dots i_\beta}$ in the same way as for definitions
\ref{kenya} and \ref{kenya2}, with the extra restriction that
whenever we refer to a free index ${}_{i_d}$, we will mean that $d\le
\alpha$. We notice that with this new definition,
 we can have two tensor fields $C^{i_1\dots i_\beta}$,
$C^{i_1\dots i_\alpha}$ with $\alpha<\beta$, so that the
double character of $C^{i_1\dots i_\alpha}$
and the $\alpha$-double character of $C^{i_1\dots i_\beta}$ are
 equivalent. The same is true of refined double characters. We
 note that this notion depends on the order of the free indices in
 $C^{i_1\dots i_\beta}$.

\par Furthermore, we note that we will sometimes be referring to a $u$-simple
character $\vec{\kappa}_{simp}$ to stress that the information encoded
will refer to the $u$ factors $\nabla\phi_1,\dots,\nabla\phi_u$. Analogously,
 we will sometimes refer to a $(u,\mu)$-refined double character
 to stress that the information encoded refers to the $u$ factors $\nabla\phi_1,\dots,\nabla\phi_u$
and the $\mu$ free indices ${}_{i_1},\dots ,{}_{i_\mu}$.
\newline

{\it Forbidden Cases:} Now, we introduce a last definition of certain ``forbidden cases'' in which 
the Proposition \ref{giade} will not apply. Firstly we introduce  a definition.

\begin{definition}
\label{simplefactor} Given a simple character
$\vec{\kappa}_{simp}$ and any factor
$T=S_{*}\nabla^{(\nu)}R_{ijkl}$ in $\vec{\kappa}_{simp}$, we will
say that $T$ is simple if it is not contracting against any
factors $\nabla\phi'_h$ in $\vec{\kappa}_{simp}$.

\par Also, given a factor $T=\nabla^{(B)}\Omega_k$, we will say that
$T$ is simple if it is not contracting against any factor
$\nabla\phi_h$ in $\vec{\kappa}_{simp}$.
\end{definition}

We recall that $\sigma_2$ stands for the number of factors 
$S_*\nabla^{(\nu)}R_{ijkl}$ in $\vec{\kappa}_{simp}$.

\begin{definition}
\label{forbidden} 
A tensor field in the form (\ref{form2}) will be called ``forbidden'' 
only when $\sigma_2>0$, under the following additional restrictions: 

If $\sigma_2=1$, it will be forbidden if: 

\begin{enumerate}
 \item {Any factor $\nabla^{(m)}R_{ijkl}$ must  
have all its derivative indices 
contracting against factors $\nabla\phi_x$ and contain no free indices.} 

\item{Any factor $\nabla^{(p)}\Omega_h$ must 
have $p=2$, be simple, and contain no free indices.}

\item{The factor $S_{*}\nabla^{(\nu)}R_{ijkl}$ 
must have $\nu=0$, be simple, 
and contain exactly one (special) free index.}
\end{enumerate} 

If $\sigma_2>1$, it will  be forbidden if: 

\begin{enumerate}
 \item {Any factor $\nabla^{(m)}R_{ijkl}$ must  
have all its derivative indices contracting against factors 
$\nabla\phi_x$ ans contain at most one (necesarily special) free index.} 

\item{Any factor $\nabla^{(p)}\Omega_h$ must 
have $p=2$. If it is simple, it can 
contain at most one free index; if it is non-simple, then 
it must contract against exactly one factor $\nabla\phi_h$ and contain no free indices.}

\item{The factor $S_{*}\nabla^{(\nu)}R_{ijkl}$ 
must have $\nu=0$, be simple, 
and contain at most one free index. Moreover at least one of 
the factors $S_*R_{ijkl}$ must contain a special free index.}
\end{enumerate}
\end{definition}

Finally, we note  that in stating Proposition \ref{giade} we will be
formally considering linear combinations of tensor fields of
different ranks.

\begin{proposition}
\label{giade}
Consider two linear combinations of acceptable tensor fields in
the form (\ref{form2}):

$$\Sum_{l\in L_\mu} a_l
C^{l,i_1\dots i_{\mu}}_{g} (\Omega_1,\dots
,\Omega_p,\phi_1,\dots ,\phi_u),$$

$$\Sum_{l\in L_{>\mu}} a_l
C^{l,i_1\dots i_{\beta_l}}_{g} (\Omega_1,\dots
,\Omega_p,\phi_1,\dots ,\phi_u),$$
 where each tensor field above has real length $\sigma\ge 3$ and a given
simple character $\vec{\kappa}_{simp}$. We assume  that for each
$l\in L_{>\mu}$,  $\beta_l\ge \mu+1$. We also assume 
 that none of the tensor fields of maximal refined double
character in $L_\mu$ are ``forbidden'' (see Definition (\ref{forbidden})).

 We denote by
$$\Sum_{j\in J} a_j C^j_{g}(\Omega_1,\dots ,\Omega_p,
\phi_1,\dots ,\phi_u)$$ a generic linear combination of complete
contractions (not necessarily acceptable) in the form
(\ref{form1}) that are simply subsequent to
$\vec{\kappa}_{simp}$.\footnote{Of course if 
$Def(\vec{\kappa}_{simp})=\emptyset$ then  by
definition $\Sum_{j\in J} \dots=0$.} We assume that:

\begin{equation}
\label{hypothese2}
\begin{split}
&\Sum_{l\in L_\mu} a_l Xdiv_{i_1}\dots Xdiv_{i_{\mu}}
C^{l,i_1\dots i_{\alpha}}_{g} (\Omega_1,\dots
,\Omega_p,\phi_1,\dots ,\phi_u)+
\\&\Sum_{l\in L_{>\mu}} a_l Xdiv_{i_1}\dots Xdiv_{i_{\beta_l}}
C^{l,i_1\dots i_{\beta_l}}_{g} (\Omega_1,\dots
,\Omega_p,\phi_1,\dots ,\phi_u)+
\\& \Sum_{j\in J} a_j
C^j_{g}(\Omega_1,\dots ,\Omega_p,\phi_1,\dots ,\phi_u)=0.
\end{split}
\end{equation}

\par We draw our conclusion with a little more notation: We break the index set
$L_\mu$ into subsets $L^z, z\in Z$, ($Z$ is finite)
with the rule that each $L^z$  indexes tensor fields
with the same refined double character, and conversely two tensor
fields with the same refined double character must be indexed in
the same $L^z$. For each index set $L^z$, we denote the
refined double character in question by $\vec{L}^z$. Consider
the subsets $L^z$ that index the tensor fields of {\it maximal} refined double
 character.\footnote{Note that in any set $S$ of $\mu$-refined double characters
 with the same simple character there is going to be a subset $S'$
 consisting of the maximal refined double characters.}
 We assume that the index set of those $z$'s
is $Z_{Max}\subset Z$.

We claim that for each $z\in Z_{Max}$ there is some linear
combination of acceptable $(\mu +1)$-tensor fields,

$$\Sum_{r\in R^z} a_r C^{r,i_1\dots i_{\alpha +1}}_{g}(\Omega_1,
\dots ,\Omega_p,\phi_1,\dots ,\phi_u),$$ where 
 each $C^{r,i_1\dots i_{\mu +1}}_{g}(\Omega_1, \dots
,\Omega_p,\phi_1,\dots ,\phi_u)$ has a $\mu$-double
character $\vec{L^z_1}$ and also the same set of factors 
$S_{*}\nabla^{(\nu)}R_{ijkl}$ as in $\vec{L}^z$ contain 
special free indices, so that:

\begin{equation}
\label{bengreen}
\begin{split}
& \Sum_{l\in L^z} a_l C^{l,i_1\dots i_\mu}_{g}
(\Omega_1,\dots ,\Omega_p,\phi_1,\dots
,\phi_u)\nabla_{i_1}\upsilon\dots\nabla_{i_\mu}\upsilon -
\\&\Sum_{r\in R^z} a_r X div_{i_{\mu +1}}
C^{r,i_1\dots i_{\mu +1}}_{g}(\Omega_1,\dots
,\Omega_p,\phi_1,\dots ,\phi_u)\nabla_{i_1}\upsilon\dots
\nabla_{i_\mu}\upsilon=
\\& \Sum_{t\in T_1} a_t
C^{t,i_1\dots i_\mu}_{g}(\Omega_1,\dots ,\Omega_p,,\phi_1,\dots
,\phi_u)\nabla_{i_1}\upsilon\dots \nabla_{i_\mu}\upsilon,
\end{split}
\end{equation}
modulo complete contractions of length $\ge\sigma +u+\mu +1$.
Here each
$$C^{t,i_1\dots i_\mu}_{g}(\Omega_1,\dots
,\Omega_p,\phi_1,\dots ,\phi_u)$$ is acceptable and is either
simply or doubly subsequent to $\vec{L}^z$.\footnote{Recall that
``simply subsequent'' means that the simple character of
$C^{t,i_1\dots i_\mu}_{g}$ is subsequent to $Simp(\vec{L}^z)$.}
\end{proposition}

{\it Trivial observation:} We recall that when a tensor field is acceptable,
then by definition it does not have two free indices (say
${}_{i_q},{}_{i_w}$) that are indices ${}_i,{}_j$ or ${}_k,{}_l$
in the same curvature factor. Thus, such tensor fields are not allowed
in our Proposition hypothesis, (\ref{hypothese2}). Nonetheless,
the conclusion of the above Proposition would still be true if we
did allow such tensor fields: It suffices to observe that this
sublinear combination would vanish both in the hypothesis of our
Proposition and in its conclusion. This is straightforward, 
by virtue of the antisymmetry of those indices.

\par Now, Proposition \ref{giade} has a Corollary which will be used more often
 than the Proposition itself:

\begin{corollary}
\label{corgiade} Assume equation (\ref{hypothese2}) (and 
again assume that the maximal refined double characters appearing
there are not ``forbidden''). We then claim that there is a linear
combination of acceptable $(\mu+1)$-tensor fields
$$\Sum_{h\in H} a_h C^{h,i_1,\dots ,i_{\mu +1}}_{g}
(\Omega_1,\dots ,\Omega_p,\phi_1,\dots ,\phi_u)$$ with simple
character $\vec{\kappa}_{simp}$, so that:

\begin{equation}
\label{conrconclu}
\begin{split}
&\Sum_{l\in L_\mu} a_l C^{l,i_1\dots i_\mu}_{g} (\Omega_1,\dots
,\Omega_p,\phi_1,\dots ,\phi_u)\nabla_{i_1}
\upsilon\dots\nabla_{i_\mu}\upsilon+
\\&\Sum_{h\in H} a_h Xdiv_{i_{\mu +1}}
C^{h,i_1\dots i_{\mu +1}}_{g} (\Omega_1,\dots
,\Omega_p,\phi_1,\dots ,\phi_u)\nabla_{i_1}
\upsilon\dots\nabla_{i_\mu}\upsilon=
\\&\Sum_{t\in T} a_t C^{t,i_1,\dots ,i_\mu}_{g}
(\Omega_1,\dots ,\Omega_p,\phi_1,\dots ,\phi_u)\nabla_{i_1}
\upsilon\dots\nabla_{i_\mu}\upsilon,
\end{split}
\end{equation}
modulo complete contractions of length $\ge\sigma +u+\mu +1$.
Here the right hand side stands for a generic linear combination
of acceptable tensor fields that are simply subsequent to
$\vec{\kappa}_{simp}$.
\end{corollary}

{\it Proof that Corollary \ref{corgiade} follows from Proposition
\ref{giade}.}
\newline

We will prove our claim by an induction.

We consider all the
 $(u,\mu)$-double
 characters $\vec{\kappa}$ with the property that \\$Simp(\vec{\kappa})=\vec{\kappa}_{simp}$.
  It follows by definition that there is a finite
 number of such refined double characters, so we denote their set by
$\{ Doub_1(\vec{L_\mu}), \dots ,Doub_U(\vec{L_\mu})\}$. We view the
above as an ordered set, with the additional restriction that for
each $a,b, 1\le a<b\le U$ $Doub_a(\vec{L_\mu})$ is not
doubly subsequent to $Doub_b(\vec{L_\mu})$. Accordingly, we break
the
 index set $L_\mu$ into subsets $L^1,\dots ,L^U$ (where if $l\in L^t$ then
 $C^{l,i_1\dots i_\mu}_{g}$ has a refined double character $Doub_t(\vec{L_\mu})$).

\par We then claim the following inductive statement: We inductively
 assume that for some $f, 1\le f\le U$, we have shown that there is a
 linear combination of acceptable $(\mu +1)$-tensor fields with simple characters $\vec{\kappa}_{simp}$,
 say
$$\Sum_{h\in H''} a_h C^{l,i_1\dots i_{\mu+1}}_{g}
(\Omega_1,\dots,\Omega_p,\phi_1,\dots , \phi_u),$$ so that:

\begin{equation}
\label{notharm}
\begin{split}
&\Sum_{w\le f}\Sum_{l\in L^w} a_l C^{l,i_1\dots
i_\mu}_{g}(\Omega_1,\dots ,\Omega_p,\phi_1,\dots ,
\phi_u)\nabla_{i_1}\upsilon\dots\nabla_{i_\mu}\upsilon-
\\&\Sum_{h\in H''} a_h Xdiv_{i_{\mu +1}}
C^{l,i_1\dots i_{\mu+1}}_{g}
(\Omega_1,\dots,\Omega_p,\phi_1,\dots ,
\phi_u)\nabla_{i_1}\upsilon\dots\nabla_{i_\mu}\upsilon=
\\&\Sum_{w=f+1}^U\Sum_{d\in D^w} a_d
C^{d,i_1\dots i_\mu}_{g}(\Omega_1,\dots ,\Omega_p,\phi_1,\dots
, \phi_u)\nabla_{i_1}\upsilon\dots\nabla_{i_\mu}\upsilon +
\\&\Sum_{t\in T} a_t
C^{t,i_1\dots i_\mu}_{g}(\Omega_1,\dots ,\Omega_p,\phi_1,\dots
, \phi_u)\nabla_{i_1}\upsilon\dots\nabla_{i_\mu}\upsilon,
\end{split}
\end{equation}
where each $C^{d,i_1\dots i_\mu}$, $d\in D^w$ has a refined double
 character $Doub_w(\vec{L})$. Write:

\begin{equation}
\label{hellenicus}
\begin{split}
&\Sum_{w= f+1}^U\Sum_{l\in L^w} a_l C^{l,i_1\dots
i_\mu}_{g}(\Omega_1,\dots ,\Omega_p,\phi_1,\dots , \phi_u)
\\&+\Sum_{w=f+1}^U\Sum_{d\in
D^w} a_d C^{d,i_1\dots i_\mu}_{g}(\Omega_1,\dots
,\Omega_p,\phi_1,\dots , \phi_u)
\\&=\Sum_{w=f+1}^U \Sum_{y\in Y^w} a_y C^{y,i_1\dots
i_\mu}_{g}(\Omega_1,\dots ,\Omega_p,\phi_1,\dots ,\phi_u),
\end{split}
\end{equation}
where the index sets $Y^w$ stand for the index sets that arise
 when we group up all the acceptable $\mu$-tensor fields  of
the same refined double character. We then claim that for
$w=f+1$ we have that there is a linear combination of acceptable
$(\mu +1)$-tensor fields, say
$$\Sum_{h\in H'''} a_h
C^{l,i_1\dots i_{\mu+1}}_{g}
(\Omega_1,\dots,\Omega_p,\phi_1,\dots , \phi_u),$$
 so that:

\begin{equation}
\label{teleiwne}
\begin{split}
&\Sum_{y\in Y^{f+1}} a_y C^{y,i_1\dots
i_\mu}_{g}(\Omega_1,\dots ,\Omega_p,\phi_1,\dots
,\phi_u)\nabla_{i_1}\upsilon\dots \nabla_{i_\mu}\upsilon-
\\&\Sum_{h\in H'''} a_h Xdiv_{i_{\mu +1}}
C^{l,i_1\dots i_{\mu+1}}_{g}
(\Omega_1,\dots,\Omega_p,\phi_1,\dots ,
\phi_u)\nabla_{i_1}\upsilon\dots\nabla_{i_\mu}\upsilon=
\\&\Sum_{t\in T} a_t
C^{t,i_1\dots i_\mu}_{g}(\Omega_1,\dots ,\Omega_p,\phi_1,\dots
,\phi_u)\nabla_{i_1}\upsilon\dots\nabla_{i_\mu} \upsilon+
\\&\Sum_{k>f+1} \Sum_{y\in Y^k} a_y C^{y,i_1\dots
i_\mu}_{g}(\Omega_1,\dots ,\Omega_p,\phi_1,\dots
,\phi_u)\nabla_{i_1}\upsilon\dots \nabla_{i_\mu}\upsilon.
\end{split}
\end{equation}
It is clear that if we can show the above, then since the set
$\{ Doub_1(\vec{L_\mu}),\dots $,\\$Doub_U(\vec{L_\mu})\}$ is finite,
 we will have shown our corollary.

\par But (\ref{teleiwne}) is not difficult to show: Since
(\ref{notharm}) holds formally we can replace the
$\nabla\upsilon$s by $Xdiv$s (see the last Lemma in the 
Appendix of \cite{alexakis1})
 and substitute into (\ref{hypothese2}) to obtain:

\begin{equation}
\label{dontharm}
\begin{split}
&\Sum_{w=f+1}^U \Sum_{y\in Y^w} a_y Xdiv_{i_1}\dots
Xdiv_{i_\mu} C^{y,i_1\dots i_\mu}_{g}(\Omega_1,\dots
,\Omega_p,\phi_1,\dots ,\phi_u)+
\\& +\Sum_{h\in H} a_h Xdiv_{i_1}\dots Xdiv_{i_\mu}
Xdiv_{i_{\mu +1}} C^{l,i_1\dots i_{\mu+1}}_{g}
(\Omega_1,\dots,\Omega_p,\phi_1,\dots , \phi_u)+
\\&\Sum_{j\in J} a_j
C^j_{g}(\Omega_1,\dots ,\Omega_p,\phi_1,\dots ,\phi_u)=0.
\end{split}
\end{equation}

\par Because the sum in the first line of (\ref{dontharm}) starts
at $w=f+1$, it follows that one of the maximal
sublinear combinations in the first line of (\ref{dontharm})
is the sublinear combination
$$\Sum_{y\in Y^{f+1}} a_y
C^{y,i_1\dots i_\mu}_{g}(\Omega_1,\dots ,\Omega_p,\phi_1,\dots
,\phi_u).$$
\par Therefore, (\ref{teleiwne})  follows immediately from the
 conclusion of Proposition \ref{giade}. $\Box$

\subsection{The main algebraic Propositions in \cite{alexakis1,alexakis2} follow from Corollary \ref{corgiade}.}

We will now show how Proposition 5.1 in \cite{alexakis1} and
 Propositions 3.1, 3.2 in \cite{alexakis2} follow  from Corollary \ref{giade}. 
Proposition 5.1 in \cite{alexakis1} and Proposition 3.1 in \cite{alexakis2}
 follow immediately from Corollary \ref{corgiade}:
Observe that in  case of Proposition 5.2 in \cite{alexakis1}, the simple character of 
the tensor fields in the equation (\ref{hypothese2}) just 
encodes the fact that there are $\sigma_1$ factors 
$\nabla^{(m)}R_{ijkl}$ and $p$ factors $\nabla^{(y)}\Omega_h, h=1,\dots,p$; 
in the setting of Proposition 3.1 in \cite{alexakis2} it additionaly encodes the fact that 
the tensor fields also contain a factor $\nabla\phi$ ($=\nabla\phi_1$) which either 
contracts against a factor $\nabla^{(y)}\Omega$ or against a derivative 
index of a factor $\nabla^{(m)}R_{ijkl}$, for each of the tensor fields in the hypothesis of that Proposition. 
There are no factors $S_{*}\nabla^{(\nu)}R_{ijkl}$ in this setting, thus $\sum_{j\in J} a_j \dots=0$,
both in the hypothesis and in the conclusion of Corollary \ref{corgiade}. 
Thus, the claims of these two  Propositions follow 
directly from the conclusion of Corollary \ref{corgiade}  
by just replacing the expression $\nabla_{i_1}\upsilon\dots\nabla_{i_\mu}\upsilon$
by a symmetrization over the indices ${}^{i_1},\dots,{}^{i_\mu}$.\footnote{See the remark after the statement of Proposition 5.1 in \cite{alexakis1}.} 

On the other hand, in order to derive the  Proposition 
3.2 in \cite{alexakis2} from Corollary \ref{corgiade}, 
we have to massage the hypothesis of that Proposition 
in order make it fit with the hypothesis of Corollary \ref{corgiade}. 
For each tensor field $C^{l,i_1\dots i_\mu}_g(\Omega_1,\dots,\Omega_p,\phi)$ and
$C^{l,i_1\dots i_{\beta_l}}_g(\Omega_1,\dots,\Omega_p,\phi)$ in Proposition 3.2 in \cite{alexakis2},
 we denote by $\tilde{C}^{l,i_1\dots i_\mu}_g(\Omega_1,\dots,\Omega_p,\phi)$ and
$\tilde{C}^{l,i_1\dots i_{\beta_l}}_g(\Omega_1,\dots,\Omega_p,\phi)$
the tensor fields that arise from them by
formally replacing the expression $\nabla^{(m)}_{r_1\dots r_m}R_{ijkl}\nabla^i\phi$ 
by $S_{*}\nabla^{(m)}_{r_1\dots r_m}R_{ijkl}\nabla^i\phi$. We then observe 
(by virtue of the second Bianchi identity) that:

\begin{equation}
\label{prongo1}
 C^{l,i_1\dots i_\mu}_g(\Omega_1,\dots,\Omega_p,\phi)=
\tilde{C}^{l,i_1\dots i_\mu}_g(\Omega_1,\dots,\Omega_p,\phi)+\sum_{j\in J} a_j 
C^{j,i_1\dots i_\mu}_g(\Omega_1,\dots,\Omega_p,\phi), 
\end{equation}

\begin{equation}
\label{prongo2}
C^{l,i_1\dots i_{\beta_l}}_g(\Omega_1,\dots,\Omega_p,\phi)=
\tilde{C}^{l,i_1\dots i_{\beta_l}}_g(\Omega_1,\dots,\Omega_p,\phi)+\sum_{j\in J} a_j 
C^{j,i_1\dots i_{\beta_l}}_g(\Omega_1,\dots,\Omega_p,\phi).
\end{equation}
Notice that the tensor fields $\tilde{C}^{l,i_1\dots i_\mu}_g(\Omega_1,\dots,\Omega_p,\phi)$,
$\tilde{C}^{l,i_1\dots i_{\beta_l}}_g(\Omega_1,\dots,\Omega_p,\phi)$ are all in the form (\ref{form2}) 
and they all all have the same {\it simple character}, 
which we denote by $\vec{\kappa}_{simp}$.
 The complete contractions in $\sum_{j\in J} a_j\dots$ in
the hypothesis of Proposition 3.2 in \cite{alexakis2} 
{\it and} the complete contractions in
\\ $Xdiv_{i_1}\dots Xdiv_{i_{\mu}}C^{j,i_1\dots i_\mu}_g(\Omega_1,\dots,\Omega_p,\phi)$,
 $Xdiv_{i_1}\dots Xdiv_{i_{\beta_l}}C^{j,i_1\dots i_{\beta_l}}_g(\Omega_1,\dots,\Omega_p,\phi)$
are all {\it simply subsequent}   to 
$\vec{\kappa}_{simp}$. 

Thus, replacing the above into the hypothesis of Proposition 3.2
in \cite{alexakis2} we obtain an equation to which 
Corollary \ref{corgiade} can be applied.\footnote{Notice that the 
extra requirement in Proposition 3.2 ensures that we do not fall under 
``a forbidden case'' of Corollary \ref{corgiade}.} 
We derive that there is a linear combination of acceptable 
tensor fields, $\sum_{h\in H} a_h C^{h,i_1\dots i_{\mu+1}}_g
(\Omega_1,\dots,\Omega_p,\phi)$
 in the form (\ref{form2}), each with a 
simple character $\vec{\kappa}_{simp}$ so that: 

\begin{equation}
 \label{ngo}
\begin{split}
&\sum_{l\in L_1} a_l \tilde{C}^{l,(i_1\dots i_\mu)}_g(\Omega_1,\dots,\Omega_p,\phi)-
Xdiv_{i_{\mu+1}}\sum_{h\in H} a_h C^{h,(i_1\dots i_\mu)i_{\mu+1}}_g
(\Omega_1,\dots,\Omega_p,\phi)=
\\&\sum_{j\in J} a_j C^{J,(i_1\dots i_\mu)}_g (\Omega_1,\dots,\Omega_p,\phi).
\end{split}
\end{equation}
Combined with equations (\ref{prongo1}), (\ref{prongo2}) above, 
(\ref{ngo}) is precisely our desired conclusion.

\section{Proof of Proposition \ref{giade}: Set up an induction and reduce to Lemmas
\ref{zetajones}, \ref{pool2}, \ref{pskovb} below.}
\label{startB}

\subsection{The proof of Proposition \ref{giade} via an induction:}
\label{sectioninduction}

\par We will prove Proposition \ref{giade} by a multiple induction on different parameters
(see the enumeration below).
 We re-write the hypothesis of the
Proposition \ref{giade}, for reference purposes. 

 We are given two
linear combinations of acceptable tensor fields in the form
(\ref{form2}) all with a given simple character
$\vec{\kappa}_{simp}$. We have the linear combination:

$$\Sum_{l\in L_\mu} a_l
C^{l,i_1\dots i_\mu}_{g}(\Omega_1,\dots ,\Omega_p, \phi_1,\dots
,\phi_u),$$ for which all the tensor fields have rank $\mu$, and also
the linear combination:
$$\Sum_{l\in L\setminus L_\mu} a_l
C^{l,i_1\dots i_a}_{g}(\Omega_1,\dots ,\Omega_p,\phi_1,\dots
,\phi_u),$$ for which all the tensor fields have rank strictly greater
than $\mu$ (we should denote the rank by $a_l>\mu$ instead of $a$,
to stress that
 the tensor fields in $L\setminus L_\mu$ have different
ranks--however we will write $C^{l,i_1\dots i_a}_{g}$, thus
abusing notation). We are assuming an equation:

\begin{equation}
\label{assumpcion}
\begin{split}
&\Sum_{l\in L_\mu} a_l Xdiv_{i_1}\dots Xdiv_{i_\mu} C^{l,i_1\dots
i_\mu}_{g}(\Omega_1,\dots ,\Omega_p, \phi_1,\dots ,\phi_u)+
\\&\Sum_{l\in L\setminus L_\mu} a_l Xdiv_{i_1}\dots Xdiv_{i_a}
C^{l,i_1\dots i_a}_{g}(\Omega_1,\dots ,\Omega_p,\phi_1,\dots
,\phi_u)+
\\&\Sum_{j\in J} a_j C^j_{g}(\Omega_1,\dots ,\Omega_p,\phi_1,\dots
,\phi_u)=0,
\end{split}
\end{equation}
which holds modulo complete contractions of length $\ge\sigma+u+1$.\footnote{We have 
now set $L_\mu\bigcup L_{>\mu}=L$.}
(Recall that $\sigma$ stands for the number of factors $\nabla^{(m)}R_{ijkl},
 S_{*}\nabla^{(\nu)}R_{ijkl}, \nabla^{(B)}\Omega_x$ in $\vec{\kappa}_{simp}$).
   We recall that each $C^j$ is simply subsequent to
$\vec{\kappa}_{simp}$.

{\it The inductive assumptions:}
We explain the inductive assumptions on Proposition \ref{giade}
in detail:

\par Denote
the left hand side of equation (\ref{assumpcion}) by
$L_{g}(\Omega_1,\dots ,\Omega_p,\phi_1,\dots ,\phi_u)$ or just
$L_{g}$ for short. For the complete contractions in $L_{g}$,
$\sigma_1$ will stand for the number of factors
$\nabla^{(m)}R_{ijkl}$ and $\sigma_2$ will stand for the number of
factors $S_{*}\nabla^{(\nu)} R_{ijkl}$. Also $\Phi$ will stand for
the total number of factors
$\nabla\phi,\nabla\tilde{\phi},\nabla\phi'$ and $-n$ will stand
for the weight of the complete contractions involved.
\newline

\begin{enumerate}

\item{We assume that Proposition \ref{giade} is true for all
linear combinations $L_{g^{n'}}$ with weight $-n'$,
$n'<n$, $n'$ even, that satisfy the
hypotheses of our Proposition.}

\item{We assume that Proposition \ref{giade} is true for all
linear combinations $L_{g}$ of weight $-n$ and real length
$\sigma'<\sigma$, that satisfy the
hypotheses of our Proposition.}

\item{We assume that Proposition \ref{giade} is true for all
linear combinations $L_{g}$ of weight $-n$ and real length
$\sigma$, with $\Phi'>\Phi$ factors $\nabla\phi,\nabla\tilde{\phi},\nabla\phi'$, 
that satisfy the hypotheses of our Proposition.}

\item{We assume that Proposition \ref{giade} is true for all
linear combinations $L_{g}$ of weight $-n$ and real length
$\sigma$, $\Phi$ factors $\nabla\phi,\nabla\tilde{\phi},\nabla\phi'$
 and with {\it fewer than $\sigma_1+\sigma_2$ curvature
factors} $\nabla^{(m)}R_{ijkl},S_{*}\nabla^{(\nu)}R_{ijkl}$,
provided $L_g$  satisfies the 
hypotheses of our Proposition.}
\end{enumerate}

\par We will then show Proposition \ref{giade} for the linear
 combinations $L_{g}$ with weight $-n$, real length $\sigma$, $\Phi$  factors
 $\nabla\phi,\nabla\phi',\nabla\tilde{\phi}$ and with
$\sigma_1+\sigma_2$ curvature factors $\nabla^{(m)}R_{ijkl}$,
$S_{*}\nabla^{(\nu)}R_{ijkl}$. So we are proving our Proposition
by
 a multiple induction on the parameters $n,\sigma, \Phi,\sigma_1+\sigma_2$
of the linear combination $L_{g}$. A trivial observation: For each
weight $-n$, there are obvious (or assumed)  bounds on the numbers $\sigma$ ($\ge 3$),
$\sigma-\sigma_1-\sigma_2$ ($\ge 0, <n$) and on the number $\Phi$ ($\le \frac{n}{2}$).
In view of this, we see that if we can show this inductive
statement, then Proposition \ref{giade} will follow by induction.
\newline

The rest of this series of papers is devoted to proving 
this inductive step of Proposition \ref{giade}. However, for simplicity we will still
refer to proving Proposition \ref{giade} rather 
than proving the inductive step of Proposition \ref{giade}. 

\subsection{Reduction of Proposition \ref{giade} to three Lemmas:} 
\label{reduction.lemmas}
We will
claim three Lemmas below: Lemmas \ref{zetajones}, \ref{pool2} and
\ref{pskovb}.\footnote{Lemma \ref{pskovb} also relies on certain preparatory 
Lemmas, \ref{oui}, \ref{oui2} which will be proven in \cite{alexakis6}.} 
We will then prove in the next section that
Proposition \ref{giade}  follows from these three Lemmas (apart
from some exceptional cases, where we will derive
Proposition \ref{giade} directly--they will be presented in the
paper \cite{alexakis5} in this series). As these Lemmas are rather technical, we give here
the gist of their claims, and also indicate, very roughly, how
they will imply Proposition \ref{giade}, by virtue of our
inductive assumptions above. A rigorous proof 
of how Lemmas \ref{zetajones}, \ref{pool2} and
\ref{pskovb} imply Proposition
 \ref{giade} will be given in section 
\ref{laderivation} of the present paper.
\newline

{\bf General Discussion of Ideas:}  We distinguish three cases
 regarding the tensor fields of rank $\mu$ appearing in
(\ref{assumpcion}). In the first case, some of the $\mu$-tensor
 fields (indexed in $L_\mu$) have special free indices belonging to
factors $S_{*}\nabla^{(\nu)}R_{ijkl}$.\footnote{Recall that a special
 free index in a factor $S_{*}\nabla^{(\nu)}R_{ijkl}$ is
 one of the indices ${}_k,{}_l$.} In the second case,
none of the $\mu$-tensor fields (indexed in $L_\mu$) have special free indices belonging to factors
$S_{*}\nabla^{(\nu)}R_{ijkl}$, but some have
special free indices in factors $\nabla^{(m)}R_{ijkl}$.\footnote{Recall that
a special free index that belongs to $\nabla^{(m)}R_{ijkl}$
is one of the indices ${}_i,{}_j,{}_k,{}_l$} In the
third case, there are no special free indices in any
$\mu$-tensor field in (\ref{assumpcion}). The three Lemmas
\ref{zetajones}, \ref{pool2} and \ref{pskovb} correspond to these three cases.

\par A note: It follows that in the first
 case above, the $\mu$-tensor fields in (\ref{assumpcion})
 of {\it maximal} refined double character will have a special free index
 in some factor $S_{*}\nabla^{(\nu)}R_{ijkl}$ (this
follows from the definition of {\it maximal refined double character}, Definition \ref{refdoub}).
 It also follows that in the second case the
  $\mu$-tensor fields of {\it maximal} refined double character will have a special free index
 in some factor $\nabla^{(m)}R_{ijkl}$, while in the third case, the maximal $\mu$-tensor
 fields will have no special free indices.
We now outline the statements of the Lemmas
\ref{zetajones}, \ref{pool2}, \ref{pskovb}:
\newline

\par In the roughest terms, in each of the three cases above, the corresponding
Lemma states the following: We ``canonically'' pick out some sub-linear
combination of the maximal $\mu$-tensor fields (for this discussion
 we denote the index set of this sublinear combination by $\overline{L}^{Max}_\mu\subset L_\mu$). In the
 first two cases, we consider each $C^{l,i_1\dots i_\mu}_{g}$, $l\in \overline{L}^{Max}_\mu$ and
 canonically pick out one (or a set of) special free indices.
 For the purposes of this discussion, we will assume that we have
 canonically picked out one free index, and we will assume it is the index ${}^{i_1}$ in each
$C^{l,i_1\dots i_\mu}_{g}$, $l\in \overline{L}^{Max}_\mu$.

{\bf A rough description of the claim of Lemmas \ref{zetajones}
and \ref{pool2}:} For Lemmas \ref{zetajones} and
\ref{pool2}, our claim is roughly an equation of the form:

\begin{equation}
\label{gendisc}
\begin{split}
&\Sum_{l\in \overline{L}^{Max}_\mu} a_l Xdiv_{i_2}\dots
Xdiv_{i_\mu}C^{l,i_1\dots i_\mu}_{g}\nabla_{i_1}\phi_{u+1}
\\&+\Sum_{\nu\in N} a_\nu Xdiv_{i_2}\dots Xdiv_{i_\mu}C^{\nu,i_1\dots
i_\mu}_{g}\nabla_{i_1}\phi_{u+1}
\\&+\Sum_{p\in P} a_p
Xdiv_{i_2}\dots Xdiv_{i_{\mu+1}}C^{p,i_1\dots
i_{\mu+1}}_{g}\nabla_{i_1}\phi_{u+1}+ \Sum_{j\in J'} a_j
C^{j,i_1}_{g}\nabla_{i_1}\phi_{u+1}=0,
\end{split}
\end{equation}
which holds modulo longer complete contractions. Here the tensor fields
indexed in $\overline{L}^{Max}_\mu\bigcup N\bigcup P$ are ``
acceptable'' in a  suitable sense. Also, the $(\mu-1)$-tensor
fields indexed in  $\overline{L}^{Max}_\mu$ have a specified
``simple character'' $\vec{\kappa}'_{simp}$ (in a suitable sense),
where this ```simple character'' encodes the pattern of which
factors the different terms $\nabla\phi_1,\dots ,
\nabla\phi_u,\nabla\phi_{u+1}$ are contracting against. Also,
 all the $(\mu-1)$-tensor fields indexed in $N$
 have this ``simple character'' $\vec{\kappa}'_{simp}$,
but they are ``doubly subsequent'' (in a suitable sense)
 to the tensor fields indexed in $\overline{L}^{Max}_\mu$.
Finally, the complete contractions indexed in $J'$ are ``simply
subsequent'' (in a suitable sense) to $\vec{\kappa}'_{simp}$. {\it
Note:} The tensor fields in (\ref{gendisc}) will not always be in
the form (\ref{form2}), thus our usual definitions of
``character'', ``subsequent'' etc.
 do not immediately apply.
\newline

{\bf A rough description of the claim of Lemma \ref{pskovb}:} Now,
we can roughly describe the claim of Lemma \ref{pskovb} (which is
the hardest of the three): We ``canonically'' pick out a
sub-linear combination of the maximal $\mu$-tensor fields in
(\ref{assumpcion}),
 (we again denote the index set of this sublinear combination
by $\overline{L}^{Max}_\mu$). For each $\mu$-tensor field
$C^{l,i_1\dots i\mu}_{g}$, $l\in \overline{L}^{Max}_\mu$, there is
a ``canonical way'' of picking out two factors:
 The ``critical factor'' and the ``second critical factor''.\footnote{In fact, we may have a
 set of critical factors, and a set of second-critical factors, but for
 this discussion we will assume that they are unique, 
for each tensor field indexed in $\overline{L}^{Max}_\mu$.}
We distinguish cases based on the number of free indices belonging to the second critical factor. 
Case A corresponds to the case where there are at least two such;
in that case, we assume that the indices ${}_{i_1},{}_{i_2}$ belong to the
second critical factor. We then introduce a formal
 operation that erases the index ${}_{i_1}$, and adds a new derivative free
 index (denote it by $\nabla_{i_{*}}$) onto the critical factor.
 For each $l\in \overline{L}^{Max}_\mu$, we will denote the resulting tensor field by
 $\dot{C}^{l,i_2\dots i_\mu,i_{*}}_{g}$. We also consider the tensor field
$\dot{C}^{l,i_2\dots i_\mu,i_{*}}_{g}\nabla_{i_2}\phi_{u+1}$ that
is obtained from it by contracting the free index ${}^{i_2}$
against a new factor $\nabla_{i_2}\phi_{u+1}$. This new,
$(\mu-1)$-tensor field has a $(u+1)$-simple character which we will
again denote by $\vec{\kappa}'_{simp}$. The claim of Lemma
\ref{pskovb} is that an equation of the following form holds:

\begin{equation}
\label{gendisc2}
\begin{split}
&\Sum_{l\in \overline{L}^{Max}_\mu} a_l Xdiv_{i_3}\dots
Xdiv_{i_\mu}Xdiv_{i_{*}} \dot{C}^{l,i_2\dots
i_\mu,i_{*}}_{g}\nabla_{i_2}\phi_{u+1}+
\\&\Sum_{\nu\in N} a_\nu Xdiv_{i_2}\dots Xdiv_{i_\mu}C^{\nu,i_1\dots i_\mu}_{g}\nabla_{i_1}\phi_{u+1}
\\&+\Sum_{p\in P} a_p Xdiv_{i_2}\dots Xdiv_{i_{\mu+1}}C^{p,i_1\dots
i_{\mu+1}}_{g}\nabla_{i_1}\phi_{u+1}
\\& +\Sum_{j\in J'} a_j
C^{j,i_1}_{g}\nabla_{i_1}\phi_{u+1}=0;
\end{split}
\end{equation}
here the tensor fields indexed in $N$ are acceptable and have a
simple character $\vec{\kappa}'_{simp}$, but they are doubly
 subsequent to the tensor fields in the first line.
The tensor fields in $P$ have rank $>\mu-1$ (but they may fail to
be acceptable), while the complete contractions in $J'$   are
simply subsequent to $\vec{\kappa}'_{simp}$.

\subsection{The rigorous formulation of Lemmas \ref{zetajones},
\ref{pool2}, \ref{pskovb}:}
\label{rigorousform}

\par Consider equation (\ref{assumpcion}). Denote by
$L_\mu^{Max}\subset L_\mu$ the index set of the tensor fields
 of maximal refined double character (recall there may
 be many maximal refined double characters). A note is in order here:
As explained (roughly) in the above discussion,
 in order to state our Lemmas we will be ``canonically''
 picking out some factor, and contracting one of its free indices against
a new factor $\nabla\phi_{u+1}$.  In particular,
for Lemmas \ref{zetajones} and \ref{pool2}, we will
 be defining a ``critical factor'' for the tensor
 fields in the equation (\ref{assumpcion}),
while for Lemma \ref{pskovb} we will be defining
 both a ``critical factor'' and a ``second critical factor''
for the tensor fields in (\ref{assumpcion}). We will make a
preliminary note here regarding these notions:

{\it Note on ``critical factors'':} The ``critical factor''
 (or factors) will be defined for
 each of the tensor fields or complete contractions in
(\ref{assumpcion}). In other words, once we specify the
critical factor(s), we will be able to examine any tensor
 field $C^{l,i_1\dots i_\beta}_{g}$ or complete contraction
$C^j_{g}$  in (\ref{assumpcion})
 and unambiguously pick out the (set of) critical factor(s) in $C^{l,i_1\dots i_\beta}_{g}$ or $C^j_{g}$.
In particular, the critical factor (or set of critical factors)
will be defined to be one of
the following: Either it will be a factor $\nabla^{(y)}\Omega_a$,
for some particular value of $a, 1\le a\le p$, or it will be
defined to be the curvature factor that is contracting against a
given factor $\nabla\phi_b$, for some chosen value of $b$, $1\le
b\le u$, or we will define the set of critical factors (in each of
the contractions in (\ref{assumpcion})) to stand for the set of
curvature factors $\nabla^{(m)}R_{ijkl}$ that are not contracting
against {\it any} factors $\nabla\phi_b$.

\par In order to avoid confusion further down, we will also remark that the way the critical factor is
specified is by examining the $\mu$-tensor fields in
(\ref{assumpcion}) that have a
 {\it maximal} refined double character. Nonetheless,
 once the critical factor(s) has (have) been specified, we will be able to
look at {\it any} tensor field or complete contraction in
(\ref{assumpcion}) and pick it (them) out.  All this discussion is
also true for the ``second critical factor'' (which will only be
defined in the setting of Lemma \ref{pskovb}).
\newline

{\bf Rigorous formulation of Lemma \ref{zetajones}:}
\newline

\par Our first Lemma applies to the case where there are $\mu$-tensor
 fields in (\ref{assumpcion})
for which some factors $S_{*}\nabla^{(\nu)}R_{ijkl}$ have special
free indices (this will be called case I).

\par In each $C^{l,i_1\dots i_\mu}_{g}$,
$l\in L^{Max}_\mu$,\footnote{$L^{Max}_\mu:=\bigcup_{z\in Z_{Max}}L^z$ is the index set of $\mu$-tensor fields
of {\it maximal} refined double character in (\ref{assumpcion}).}
we pick out the factors $T_1=S_{*}\nabla^{(\nu)}R_{ijkl}$ with a
special free index\footnote{Recall that a special free index in a
factor $S_{*}\nabla^{(\nu)}R_{ijkl}$ is one of the indices
${}_k,{}_l$.}
 (observe that by the definition of ``maximal'' refined double character
and by the assumption in the previous paragraph there will be
such factors in each tensor field indexed in $L^{Max}_\mu$).

\par  Among those factors,
 we pick out the ones with the maximum
number of free indices (in total). We denote this maximum number
of free indices by $M$. It follows from the definition of {\it
ordering} among refined double characters that this number $M$ is
universal among all $C^{l,i_1\dots i_\mu}_{g}$, $l\in
L^{Max}_\mu$.

\begin{definition}
\label{critfactor} 
 For
each $l\in L^{Max}_\mu$ we list all the factors
$S_{*}\nabla^{(\nu)}R_{ijkl}$ which contain
a special free index and also have $M$ free indices in total:
 $\{F_1,\dots F_\alpha\}$. If at least one of these factors $F_h$
is contracting against at least one factor $\nabla\phi'_b$, then we define
the {\it critical factor} to be the factor from the list above which is contracting against
the factor $\nabla\phi'_b$ with the smallest value for $b$ (say $Min$).
If no factors $F_h$ in the above list are contracting
against $\nabla\phi'_b$'s, then we define the {\it critical factor} to be
the factor from the list above which is contracting against
the factor $\nabla\tilde{\phi}_r$ for
the smallest value for $r$ (say $Min$).
\end{definition}

 We then denote by
$L^{Max}_{\mu,Min}\subset L^{Max}_\mu$ the index set of the tensor
fields of maximal refined double character $C^{l,i_1\dots i_\mu}_{g}$ for which the
 factor $S_{*}\nabla^{(\nu)}R_{ijkl}$ that is  contracting against $\nabla\phi'_{Min}$ or
  $\nabla\tilde{\phi}_{Min}$, respectively, contains $M$ free indices
   in total and one of them is special.

Let us observe that there exists some subset $Z_{Max}'\subset  Z_{Max}$ so that:

\begin{equation}
\label{wei}  L_{\mu,Min}^{Max}=\bigcup_{z\in Z'_{Max}}L^z.
\end{equation}

\par Now, with no loss of generality (only for notational convenience),
 we assume that for each
$l\in L_{\mu,Min}^{Max}$, the critical factor $S_{*}\nabla^{(\nu)}R_{ijkl}$ against which
$\nabla\tilde{\phi}_{Min}$ contracts has the index ${}_k$ being the free
index ${}_{i_1}$. We also recall that $\vec{\kappa}_{simp}$
is the simple
character of the tensor fields indexed in $L_\mu$. We will then
denote by $\vec{\kappa}^z$ the refined double character of
 each $C^{l,i_1\dots ,i_\mu}_{g}$, $l\in L^z$, $z\in Z'_{Max}$.

 Under the assumptions above, our claim is the following:

\begin{lemma}
\label{zetajones}
\par Assume (\ref{assumpcion}), with weight $-n$, real length $\sigma$,
 $u=\Phi$ and $\sigma_1+\sigma_2$ factors
$\nabla^{(m)}R_{ijkl},S_{*}\nabla^{(\nu)}R_{ijkl}$--assume also
that the tensor fields of maximal refined double character are not
 ``forbidden'' (see Definition \ref{forbidden}).
Suppose that
there are $\mu$-tensor fields in (\ref{hypothese2}) with at least one 
 special free index in  a factor $S_{*}\nabla^{(\nu)}R_{ijkl}$.
We then claim that there is a linear combination of
acceptable tensor fields,
$$\Sum_{p\in P} a_p C^{p,i_1\dots i_{b}}_{g}
(\Omega_1,\dots ,\Omega_p,\phi_1,\dots ,\phi_u)$$ each with
$b\ge\mu+1$, with a simple character $\vec{\kappa}_{simp}$ and
where each \\$C^{p,i_1\dots i_b}_{g}(\Omega_1,\dots
,\Omega_p,\phi_1,\dots ,\phi_u)$ has the property that the free
index ${}_{i_1}$ is the index ${}_k$ in the critical factor
$S_{*}\nabla^{(\nu)}R_{ijkl}$ against which
$\nabla\tilde{\phi}_{Min}$ is contracting, so that modulo complete
contractions of length $\ge\sigma +u+2$:

\begin{equation}
\label{narod}
\begin{split}
&\Sum_{z\in Z'_{Max}}\Sum_{l\in L^z} a_l Xdiv_{i_2}\dots
Xdiv_{i_\mu}C^{l,i_1\dots i_\mu}_{g}(\Omega_1,\dots
,\Omega_p,\phi_1,\dots ,\phi_u)\nabla_{i_1}\phi_{u+1}+
\\&\Sum_{\nu\in N} a_\nu Xdiv_{i_2}\dots
Xdiv_{i_\mu}C^{\nu,i_1\dots i_\mu}_{g}(\Omega_1,\dots
,\Omega_p,\phi_1,\dots ,\phi_u)\nabla_{i_1}\phi_{u+1} -
\\&\Sum_{p\in P} a_p Xdiv_{i_2}\dots
Xdiv_{i_{b}}C^{p,i_1\dots i_{b}}_{g} (\Omega_1,\dots
,\Omega_p,\phi_1,\dots ,\phi_u)\nabla_{i_1} \phi_{u+1}=
\\& \Sum_{t\in T} a_t C^{t,i_{*}}_{g}
(\Omega_1,\dots , \Omega_p,\phi_1,\dots
,\phi_u)\nabla_{i_{*}}\phi_{u+1}.
\end{split}
\end{equation}
 Here each $C^{\nu,i_1\dots i_\mu}_{g}(\Omega_1,\dots
,\Omega_p,\phi_1,\dots ,\phi_u)\nabla_{i_1}\phi_{u+1}$ is
acceptable and has a simple character $\vec{\kappa}_{simp}$
 (and ${}_{i_1}$ is again the index ${}_k$ in the critical factor
  $S_{*}\nabla^{(\nu)}R_{ijkl}$), but
also has either strictly fewer than $M$ free indices in the critical factor
or is doubly subsequent to each $\vec{\kappa}^z, z\in Z'_{Max}$.
 Also, each  \\$C^{t,i_{*}}_{g} (\Omega_1,\dots , \Omega_p,\phi_1,\dots
,\phi_u)\nabla_{i_{*}}\phi_{u+1}$ is in the form (\ref{form1}) and is either simply subsequent to
$\vec{\kappa}_{simp}$ {\it or}
\\ $C^{t,i_{*}}_{g} (\Omega_1,\dots ,
\Omega_p,\phi_1,\dots ,\phi_u)$ has a  $u$-simple character
$\vec{\kappa}_{simp}$ but the index ${}_{i_{*}}$ is not a special index.
 All complete contractions have the same weak $(u+1)$-simple character.
\end{lemma}

\par If we can prove the above then (as we will show in the next subsection) by iterative
repetition we can reduce ourselves to proving Proposition
\ref{giade} with the extra assumption that for each $l\in L_\mu$
there are no special free indices in any factor
$S_{*}\nabla^{(\nu)}R_{ijkl}$. Lemma \ref{pool2} will apply to this subcase.
\newline

{\bf Rigorous formulation of Lemma \ref{pool2}:}
\newline

\par We now assume that all $\mu$-tensor fields in (\ref{assumpcion}) have no special free indices in factors
$S_{*}\nabla^{(\nu)}R_{ijkl}$, but certain $\mu$-tensor fields do
have special free indices in factors $\nabla^{(m)}R_{ijkl}$--this
will be called case II. We will then pick out those $\mu$-tensor
fields in (\ref{assumpcion}) that have special free indices in
factors $\nabla^{(m)}R_{ijkl}$. (If there are no such tensor
fields, we may proceed to Lemma \ref{pskovb}).

\par In order to state our Lemma we will need to define the
critical factor (or set of critical factors) for the tensor fields
appearing in (\ref{assumpcion}), in this setting:

\begin{definition}
\label{crtifact1} Firstly we consider the case where there are
factors $\nabla^{(m)}R_{ijkl}$ with two special free indices among
the $\mu$-tensor fields in (\ref{assumpcion}). We will define the
 critical factor(s) in that setting:

\par Among all the $\mu$-tensor fields with maximal refined double characters
in (\ref{assumpcion}),
we pick out all the factors $\nabla^{(m)}R_{ijkl}$ with
two special free indices. Among those factors, we pick out the ones  with
the maximal total number of free indices, say $M\ge 2$.
Denote that list by $\{T_1,\dots T_\pi\}$. (All the $T_i$'s are in
the form $\nabla^{(m)}R_{ijkl}$). We inquire whether any of the
factors $T_1,\dots T_\pi$ are contracting against a factor
$\nabla\phi_h$. If so, we define the critical factor to be the
factor $T_i$ that is contracting against the factor $\nabla\phi_o$
for the smallest $o$. If none of the factors $T_1,\dots ,T_\pi$
are contracting against a factor $\nabla\phi_h$ we define the set of 
critical factors to be the set of factors $\nabla^{(m)}R_{ijkl}$ 
which are not contracting against any factor $\nabla\phi_h$.

\par The same definition can be applied to define a critical
factor in the case where there are no factors $\nabla^{(m)}R_{ijkl}$
with two special free indices but there are factors
$\nabla^{(m)}R_{ijkl}$ with one free index (list them out as
$\{T_1,\dots ,T_{\pi'}\}$ and proceed as above).
\end{definition}

\par Now, we index the maximal $\mu$-tensor fields with $M$ free indices in
the critical factor in the index set $\bigcup_{z\in Z'_{Max}}L^z\subset L_\mu$.

\begin{definition}
\label{sspecialset}
 For each $z\in Z'_{Max}$ define
$I_{*,l} \subset I_l$ (recall that $I_l$ stands for the set of
free indices in the tensor field $C^{l,i_1\dots i_\mu}_{g}$)
to be the set of special free indices that belong to the critical
factor (if it is unique), or to one of the critical factors.
\end{definition}

{\it Note:} By virtue of the definition of the maximal refined
double characters and of the critical factors, we observe that for
any two $l_1,l_2\in \bigcup_{z\in Z'_{Max}}L^z$, we will have
 $|I_{*,l_1}|=|I_{*,l_2}|$.

\par Now, for each $z\in Z'_{Max}$ we consider the $(\mu-1)$-tensor fields in
the linear combination
$$\Sum_{l\in L^z} a_l \Sum_{i_h\in I_{*,l}}C^{l,i_1\dots i_\mu}_{g}
(\Omega_1,\dots ,\Omega_p,\phi_1,\dots
,\phi_u)\nabla_{i_h}\phi_{u+1}$$
 and we will write them as a linear combination of
 $(\mu-1)$-tensor fields in the form (\ref{form2}) (plus error terms--see below):

 \par  For each $l\in L^z, z\in Z'_{Max}$ and each  ${}_{i_h}\in I_{*,l}$ (we may assume
 with no loss of generality that ${}_{i_h}$ is the index ${}_i$
in some factor $\nabla^{(m)}R_{ijkl}$), we denote by
$\tilde{C}^{l,i_1\dots i_\mu}_{g}(\Omega_1,\dots
,\Omega_p,\phi_1,\dots ,\phi_u)\nabla_{i_h}\phi_{u+1}$ the tensor
field that arises from $C^{l,i_1\dots i_\mu}_{g}(\Omega_1,\dots
,\Omega_p,\phi_1,\dots ,\phi_u)\nabla_{i_h}\phi_{u+1}$ by
replacing the expression \\$\nabla^{(m)}_{r_1\dots r_m}R_{i_hjkl}\nabla^{i_h}\phi_{u+1}$
by an expression $S_{*}\nabla^{(m)}_{r_1\dots r_m}R_{i_hjkl}\nabla^{i_h}\phi_{u+1}$. By
the first and second Bianchi identity, it then follows that:

\begin{equation}
\label{joeharris}
\begin{split}
&C^{l,i_1\dots i_\mu}_{g}(\Omega_1,\dots ,\Omega_p,\phi_1,\dots
,\phi_u)\nabla_{i_h}\phi_{u+1}=
\\&\tilde{C}^{l,i_1\dots i_\mu}_{g}(\Omega_1,\dots
,\Omega_p,\phi_1,\dots ,\phi_u)\nabla_{i_h}\phi_{u+1}+
\\&\Sum_{t\in T} a_t C^{t,i_1\dots i_\mu}_{g}
(\Omega_1,\dots ,\Omega_p,\phi_1,\dots
,\phi_u)\nabla_{i_h}\phi_{u+1},
\end{split}
\end{equation}
where each $C^{t,i_1\dots i_\mu}_{g} (\Omega_1,\dots
,\Omega_p,\phi_1,\dots ,\phi_u)\nabla_{i_h}\phi_{u+1}$ has the
factor $\nabla\phi_{u+1}$ contracting against a derivative index
in a factor $\nabla^{(m)}R_{ijkl}$-see the statement of Lemma \ref{pool2}.

\par We denote by $\vec{\kappa}'_{simp}$ and $\vec{\kappa}^z$
 the $(u+1)$-simple character (respectively, the $(u+1,\mu-1)$-refined double character)
 of the tensor fields
$\tilde{C}^{l,i_1\dots i_\mu}_{g}\nabla_{i_h}\phi_{u+1}$, $l\in
L^z, z\in Z'_{Max}$. (We observe that for each $l\in L^z$, $z\in
Z'_{Max}$ the simple characters of the tensor fields
$\tilde{C}^{l,i_1\dots i_\mu}_{g} \nabla_{i_h}\phi_{u+1}$ will be
equal).

\begin{lemma}
\label{pool2} Assume (\ref{assumpcion}) with weight $-n$, real length $\sigma$,
 $u=\Phi$ and $\sigma_1+\sigma_2$ factors
$\nabla^{(m)}R_{ijkl},S_{*}\nabla^{(\nu)}R_{ijkl}$.
Suppose that no $\mu$-tensor fields
 have special free indices in factors $S_{*}\nabla^{(\nu)}R_{ijkl}$,
but some have special free indices in factors $\nabla^{(m)}R_{ijkl}$.
In the notation above we claim that
 there exists a linear combination $\Sum_{d\in D} a_d
C^{d,i_1,\dots ,i_b}_{g}(\Omega_1,\dots ,\Omega_p,\phi_1,\dots
,\phi_u,\phi_{u+1})$ of acceptable tensor fields with a
$(u+1)$-simple character $\vec{\kappa}'_{simp}$ and rank $\ge \mu$, so that:

\begin{equation}
\label{vecFskillb}
\begin{split}
&\Sum_{z\in Z'_{Max}}\Sum_{l\in L^z} a_l \Sum_{i_h\in
I_{*,l}}Xdiv_{i_1}\dots\hat{Xdiv}_{i_h}\dots
 Xdiv_{i_\mu}\tilde{C}^{l,i_1\dots
i_\mu}_{g}(\Omega_1,\dots ,\Omega_p,\phi_1,\dots
,\phi_u)\\&\nabla_{i_h}\phi_{u+1}
+\Sum_{\nu\in N} a_\nu Xdiv_{i_2}\dots Xdiv_{i_{\mu}} C^{\nu,i_1\dots i_{\mu}}_{g}(\Omega_1,\dots
,\Omega_p,\phi_1,\dots ,\phi_u)\nabla_{i_1}\phi_{u+1} 
\\&- \Sum_{d\in D} a_d Xdiv_{i_1}\dots Xdiv_{i_b} C^{d,i_1,\dots
,i_b}_{g}(\Omega_1,\dots ,\Omega_p,\phi_1,\dots
,\phi_u,\phi_{u+1})=
\\& \Sum_{t\in T} a_t C^{t,i_{*}}_{g}
(\Omega_1,\dots ,\Omega_p,\phi_1,\dots
,\phi_u,\phi_u)\nabla_{i_{*}}\phi_{u+1},
\end{split}
\end{equation}
where the $(\mu-1)$-tensor fields $C^{\nu,i_1\dots
i_{\mu}}_{g}(\Omega_1,\dots ,\Omega_p,\phi_1,\dots
,\phi_u)\nabla_{i_1}\phi_{u+1}$ are acceptable, have
$(u+1)$-simple character $\vec{\kappa}'_{simp}$ but also either
 have fewer than $M$ free indices in the factor against which
  $\nabla_{i_h}\phi_{u+1}$ contracts,\footnote{``Fewer than $M$ free indices''
  where we also count the free index ${}_{i_h}$.} or are
doubly subsequent
 to all the refined double characters $\vec{\kappa}^z$,
$z\in Z'_{Max}$. Moreover we require that each $C^{\nu,i_1\dots
,i_\mu}_{g}$ has the property that at least one of the indices
${}_{r_1},\dots ,{}_{r_\nu},{}_j$ in the factor $S_{*}\nabla^{(\nu)}_{r_1\dots
r_\nu}R_{ijkl}$ is neither free nor contracting against  a factor
$\nabla\phi_h'$, $h\le u$. The complete contractions
$C^{t,i_{*}}_{g} (\Omega_1,\dots ,\Omega_p,\phi_1,\dots
,\phi_u,\phi_u)\nabla_{i_{*}} \phi_{u+1}$ are in the form (\ref{form1}) and are simply subsequent to
$\vec{\kappa}'_{simp}$.
\end{lemma}

\par We will show in  section \ref{laderivation}  that if we can prove
 the above, then we can  reduce ourselves to
proving Proposition \ref{giade} when none of the $\mu$-tensor fields in
 the Lemma hypothesis have special free indices. Lemma \ref{pskovb} will apply to 
 that setting:
\newline

{\bf Rigorous formulation of Lemma \ref{pskovb}:}
\newline

\par Recall that we have grouped up the $\mu$-tensor fields
$C^{l,i_1\dots i_\mu}_{g}$ according to their refined double
characters: $\Sum_{l\in L_\mu} a_l C^{l,i_1\dots i_\mu}_{g}=
 \Sum_{z\in Z} \Sum_{l\in L^z} a_l\dots$.  We have then
picked out the sublinear combinations in $\Sum_{l\in L_\mu} a_l
C^{l,i_1\dots i_\mu}_{g}$ which consist of tensor fields
 with the same {\it maximal} refined double character.
Thus we obtain a sublinear combination $\Sum_{z\in Z_{Max}} \Sum_{l\in
L^z} a_l C^{l,i_1\dots i_\mu}_{g}$.
\newline

\par Now, in order to state our Lemma we will distinguish two further 
 subcases; first we must introduce some more terminology.

\par We will again define the critical factor for the tensor fields in (\ref{assumpcion}):

\begin{definition}
\label{critfact2} Consider all the $\mu$-tensor fields of 
maximal refined double character in
(\ref{assumpcion}), and let $M$ stand for the maximum number of
free indices that can belong to a given factor in such a
$\mu$-tensor field (we call these
``maximal'' $\mu$-tensor fields). We then list all the factors that appear with
$M$ free indices in some maximal $\mu$-tensor field in
(\ref{assumpcion}):  $\{ T_1,\dots ,T_\pi\}$. If at least one of
those factors
 $T_l$ is of the form $\nabla^{(p)}\Omega_h$, we define the critical factor to be the factor
$\nabla^{(p)}\Omega_h$ in the list $\{ T_1,\dots ,T_\pi\}$ with
the
 smallest value $h$. If none of the factors in that list are
in the form $\nabla^{(p)}\Omega_h$, we inquire whether any factors in
the list are contracting against factors $\nabla\phi_h$ (or
$\nabla\tilde{\phi}_h$). If so, we define the critical factor in
(\ref{assumpcion}) to be the factor in the list $\{ T_1,\dots
,T_\pi\}$ that is contracting against
 the factor $\nabla\phi_h$ (or $\nabla\tilde{\phi}_h$) with the
smallest value of $h$. Finally, if none of the factors in the list
$\{ T_1,\dots ,T_\pi\}$ are contracting against a factor
$\nabla\phi_h$ (or $\nabla\tilde{\phi}_h$) (so all of them
 must be in the form $\nabla^{(m)}R_{ijkl}$), then we declare the set of critical factors to be the set  of
 factors $\nabla^{(m)}R_{ijkl}$ that are not contracting
 against any $\nabla\phi_h$ to be critical factors.
\end{definition}

 In addition to the critical
factor, we now define the {\it second critical factor}
 in (\ref{assumpcion}). The definition goes as follows:

\begin{definition}
\label{seccrit} Consider any of the maximal $\mu$-tensor fields in
(\ref{assumpcion}), $C^{l,i_1\dots ,i_\mu}_g$, $l\in \bigcup_{z\in
Z_{Max}}L^z$. If the critical factor is unique, we
 construct a list of all the non-critical
factors $T$ that belong to one of the tensor fields
$C^{l,i_1\dots ,i_\mu}_g, l\in \bigcup_{z\in Z'_{Max}}L^z$
Suppose that list is
 $\{T_1,\dots, T_\pi\}$.

\par We then pick out the second critical factor from that list
in the same way that we pick out the critical factor in definition \ref{critfact2}.

\par If there are multiple critical factors, we just define the set
 of second critical factors to be the set of critical factors.

 \par In either case, we denote by $M'$ the total number of
 free indices that belong to the (a) second critical fator.
\end{definition}

\par Now, an important note: The 
 ``critical factor'' (or factors) in (\ref{assumpcion}) has been defined based on the
 maximal refined double characters $\vec{L^z}$,
$z\in Z_{Max}$. Nonetheless, once we have chosen a critical factor
(or a set of critical factors) for the set $C^{l,i_1\dots i_\mu}_{g}$,
$l\in \bigcup_{z\in Z'_{Max}}L^z$, we may then unambiguously speak
of the critical factor(s) for {\it all} the tensor fields and
complete contractions appearing in (\ref{assumpcion}).

{\it The two cases for Lemma \ref{pskovb}:} We now distinguish two
cases on (\ref{hypothese2}): We say that (\ref{hypothese2}) (where no 
tensor fields contain special free indices) fall under 
 case A if $M'\ge 2$.
It falls under case B if $M'\le 1$.
\newline

\par Now, we will state Lemma \ref{pskovb} after we first state a
few extra claims. These claims will be proven in the paper
 \cite{alexakis6} in this series.\footnote{These claims 
involve much notation and are rather technical. The reader 
may choose to disregard them in the first reading, as 
they are not central to the argument.}
\newline

{\it The extra claims needed to state Lemma \ref{pskovb}:}
\newline

\par In order to state Lemma \ref{pskovb},
we must first show some preliminary results. We introduce some definitions:

\par  We denote by
$L_\mu^{*}\subset L_\mu$ the index set of those tensor fields
$C^{l,i_1\dots i_\mu}_{g}$ in (\ref{assumpcion}) for which some
chosen factor $\nabla^{(A)}_{r_1\dots r_A}\Omega_x$ 
(the value $x$ which determizes this factor will be 
chosen at a later stage; we may also {\it not} 
choose any such factor $\nabla^{(A)}_{r_1\dots r_A}\Omega_x$, 
in which case we set $L_\mu^*=\emptyset$)
 has $A=2$ and both
 indices ${}_{r_1},{}_{r_2}$ are free indices.

\par Also,  
we define $L_\mu^{+}\subset L_\mu$ to
 stand for the index
set of those $\mu$-tensor fields that have a 
free index (${}_{i_\mu}$ say) belonging to a
factor $S_{*}R_{ijkl}\nabla^i \tilde{\phi}_h$ (without
derivatives) and in fact ${}_j={}_{i_\mu}$.

\par  We also denote by
$$\Sum_{l\in \tilde{L}} a_l C^{l,i_1\dots i_\mu}_{g}
(\Omega_1,\dots ,\Omega_p,\phi_1,\dots ,\phi_u)$$ a  linear
combination of acceptable $\mu$-tensor fields with simple
character $\vec{\kappa}_{simp}$ which do not have special free indices
and does not contain tensor fields in any of the above two forms.

\begin{lemma}
\label{oui} Assume (\ref{assumpcion}), where the terms in the LHS of that equation
have weigh $-n$, real length $\sigma$, $\Phi$ factors 
$\nabla\phi, \nabla\phi',\nabla\tilde{\phi}$ and $\sigma_1+\sigma_2$ curvature 
factors $\nabla^{(m)}R_{ijkl},S_{*}\nabla^{(\nu)}R_{ijkl}$;\footnote{See the discussion
on the {\it induction} in subsection \ref{sectioninduction}.}  assume also  that no
$\mu$-tensor field there has any special free indices.
 We claim that there is a linear combination of
acceptable $(\mu+1)$-tensor fields,  $\Sum_{p\in P} a_p
C^{p,i_1\dots i_{\mu +1}}_{g}(\Omega_1,\dots ,\Omega_p,
\phi_1,\dots ,\phi_u)$ with a simple character
$\vec{\kappa}_{simp}$ so that:

\begin{equation}
\label{dog}
\begin{split}
&\Sum_{l\in L_\mu^{*}\bigcup L_\mu^{+}} a_l C^{l,i_1\dots
i_\mu}_{g}(\Omega_1,\dots ,\Omega_p, \phi_1,\dots
,\phi_u)\nabla_{i_1}\upsilon\dots \nabla_{i_\mu}\upsilon+
\\& \Sum_{p\in P} a_p Xdiv_{i_{\mu +1}}
C^{p,i_1\dots i_{\mu +1}}_{g}(\Omega_1,\dots ,\Omega_p,
\phi_1,\dots ,\phi_u)\nabla_{i_1}\upsilon\dots
\nabla_{i_\mu}\upsilon=
\\& \Sum_{j\in J} a_j C^{j,i_1\dots i_\mu}_{g}(\Omega_1,\dots ,\Omega_p,
\phi_1,\dots ,\phi_u)\nabla_{i_1}\upsilon\dots
\nabla_{i_\mu}\upsilon +
\\& \Sum_{l\in \tilde{L}} a_l C^{l,i_1\dots
i_\mu}_{g} (\Omega_1,\dots ,\Omega_p,\phi_1,\dots ,\phi_u)
\nabla_{i_1}\upsilon\dots\nabla_{i_\mu}\upsilon,
\end{split}
\end{equation}
modulo complete contractions of length $\ge\sigma +u+\mu +1$. The
tensor fields indexed in $J$ on the right hand side are simply
subsequent to $\vec{\kappa}_{simp}$.
\end{lemma}

\par Assuming the above Lemma, by making the $\nabla\upsilon$'s into $Xdiv$s (see the last
 Lemma in the Appendix of \cite{alexakis1}) and replacing
into (\ref{assumpcion}) we are reduced to showing our Proposition
\ref{giade} under the additional assumption that $L_\mu^{*}\bigcup
L^{+}_\mu=\emptyset$. So for the rest of this subsection we will 
be assuming that $L_\mu^{*}\bigcup L_\mu^{+}=\emptyset$.
\newline

\par We now consider the sublinear combination indexed in
$L\setminus L_\mu(=L_{>\mu})$ in (\ref{assumpcion}). We define 
$L''_{+}\subset L_{\>\mu}$ to stand for the index set of tensor
fields with a
factor $R_{ijkl}\nabla^i\tilde{\phi}_h$ for which {\it both} indices
${}_j,{}_k$ are free.

\par Now, we denote by
$$\Sum_{l\in \tilde{L}'} a_l C^{l,i_1\dots i_{\mu+1}}_{g}
(\Omega_1,\dots ,\Omega_p,\phi_1,\dots ,\phi_u)$$ a generic linear
combination of acceptable $(\mu +1)$-tensor fields that do not
contain tensor fields in
  the form described above. We then claim:

\begin{lemma}
\label{oui2} Assume (\ref{assumpcion})  with weight $-n$, real length $\sigma$,
 $u=\Phi$ and $\sigma_1+\sigma_2$ factors
$\nabla^{(m)}R_{ijkl},S_{*}\nabla^{(\nu)}R_{ijkl}$;\footnote{See the discussion in 
subsection \ref{sectioninduction}.} assume also that 
 none of the $\mu$-tensor fields have special free indices,
and  that
 $L_\mu^{*}\bigcup L_\mu^{+}=\emptyset$. 
 We claim that there exists a linear combination of
acceptable $(\mu+2)$-tensor fields, $\Sum_{p\in P} a_p C^{p,i_1\dots i_{\mu
+2}}_{g}(\Omega_1,\dots ,\Omega_p, \phi_1,\dots ,\phi_u)$ with simple
 character $\vec{\kappa}_{simp}$, so that:

\begin{equation}
\label{dog2}
\begin{split}
&\Sum_{l\in L''_{+}} a_l Xdiv_{i_1}\dots
Xdiv_{i_{\mu+1}} C^{l,i_1\dots i_{\mu+1}}_{g}(\Omega_1,\dots
,\Omega_p, \phi_1,\dots ,\phi_u)+
\\& \Sum_{p\in P} a_p Xdiv_{i_1}\dots Xdiv_{i_{\mu+2}}
C^{p,i_1\dots i_{\mu +2}}_{g}(\Omega_1,\dots ,\Omega_p,
\phi_1,\dots ,\phi_u)=
\\& \Sum_{j\in J} a_j C^{j}_{g}(\Omega_1,\dots ,\Omega_p,
\phi_1,\dots ,\phi_u)
\\&+ \Sum_{l\in \tilde{L}'} a_l Xdiv_{i_1}\dots
Xdiv_{i_{\mu+1}} C^{l,i_1\dots i_{\mu+1}}_{g} (\Omega_1,\dots
,\Omega_p, \phi_1,\dots ,\phi_u),
\end{split}
\end{equation}
modulo complete contractions of length $\ge\sigma +u+1$.
$\Sum_{j\in J}\dots$ stands for  a linear combination of complete
contractions that are simply subsequent to $\vec{\kappa}_{simp}$.
\end{lemma}

\par We observe that if we can show the above, then replacing into (\ref{assumpcion}) we
are reduced to proving Proposition \ref{giade} under the extra
 assumptions that 
\\$L_\mu^{*}\bigcup L^{+}_\mu\bigcup  L''_{+}=\emptyset$. 
So for the rest of this section we will
be making that assumption. The proof of
these two Lemmas is given in the paper \cite{alexakis6} in this series.
\newline

{\it Notation and language conventions for Lemma \ref{pskovb}:}
 Recall the two cases A, B. We will
 first formulate our claim in case  A (where $M'\ge 2$). We
introduce some notation.

\par  We define $Z'_{Max}\subset Z_{Max}$  as follows:\footnote{Recall
that $\bigcup _{z\in Z_{Max}} L^z\subset L_\mu$ stands for the
index set of the $\mu$-tensor fields.} $z\in Z'_{Max}$
if and only if $C^{l,i_1\dots i_\mu}_{g}, l\in L^z$ has $M'$ free indices in
the second critical factor (see definition \ref{seccrit}).

\par Now, we first consider the case where there is a unique second critical factor in
(\ref{assumpcion}). For each $l\in L^z$, $z\in Z'_{Max}$, we
assume with no loss of  generality that the indices
${}_{i_1},{}_{i_2}$ belong to the second
 critical factor, and that the index ${}_{i_1}$ is a derivative
 index (the second assumption can be made since all $\mu$-tensor fields in (\ref{assumpcion})
  have no special free indices now; hence if two free indices belong
   to the same factor, one of them must be a derivative index). We then denote by
$\dot{C}^{l,i_2\dots i_\mu,i_{*}}_{g}(\Omega_1,\dots,
\Omega_p,\phi_1,\dots ,\phi_u)$ the tensor field that formally  arises from
$C^{l,i_1\dots i_\mu}_{g} (\Omega_1,\dots,\Omega_p,\phi_1,\dots
,\phi_u)$ by erasing
 the free index ${}_{i_1}$ from the
 second critical factor and adding a derivative index
$\nabla_{i_{*}}$ onto the critical factor, and making the
 index ${}_{i_{*}}$ free.
We denote by $\vec{L}^{z,\sharp}$ the $(u+1,\mu-1)$-refined double
character of these $\dot{C}^{l,i_2\dots i_\mu,i_{*}}_{g}
(\Omega_1,\dots,\Omega_p,\phi_1,\dots
,\phi_u)\nabla_{i_2}\phi_{u+1}$, $l\in L^z$, $z\in Z'_{Max}$.

\par Now, the case where there are multiple second critical factors:
 If there
are $k>1$ second critical factors $T_1,\dots T_k$ in (\ref{assumpcion}), then for
each $C^{l,i_1\dots i_\mu}_{g}$, $l\in L^z, z\in Z'_{Max}$ we
denote by $\{i_1,\dots i_\alpha\}, \{ i_{\alpha+1}, \dots
i_{2\alpha}\},\dots$,$\{i_{(k-1)\alpha+1},\dots ,i_{k\alpha}\}$
the set of free indices that belong to $T_1,\dots ,T_k$ respectively. We will be
making the assumption (with no loss of generality, for the reason 
explained in the previous paragraph) that the index $i_{t\alpha+1}$ is a
derivative index for every $t=0,1,\dots ,k-1$. We then denote by
$\dot{C}^{l,i_1\dots \hat{i}_{t\alpha +1}\dots i_\mu,i_{*}}_{g}$
the tensor field that arises from $C^{l,i_1\dots i_\mu}_{g}$
 by erasing the index $i_{t\alpha+1}$ and adding a free derivative
 index ${}_{i_{*}}$ onto the (a) critical factor (and adding,
if there are multiple critical factors).

In both cases above we define $\vec{\kappa}^{+}_{simp}=
Simp(\vec{L}^{z,\sharp})$, for some $z\in Z'_{Max}$ (notice that
the definition is independent of the element $z\in Z'_{Max}$).

\par A note is needed regarding this definition: In the case where the
set of critical and second critical factors coincide, then when we
``add a free derivative
 index ${}_{i_{*}}$ onto (a) critical factor'', we will be adding
 it on any critical factor other than the one to which
 ${}_{i_{t\alpha+2}}$ belongs.
Observe that for any $l=0,1,\dots , k$ the tensor fields
 $\dot{C}^{l,i_1\dots \hat{i}_{t\alpha +1}\dots i_\mu,i_{*}}_{g}\nabla_{t\alpha+2}\phi_{u+1}$
 have the same $(u+1,\mu-1)$-refined double character, which we again denote by $\vec{L^{z,\sharp}}$,
  $z\in Z'_{Max}$ (as in the case of a unique second critical factor).

One last language convention: For uniformity, in case A
of Lemma \ref{pskovb} we will call the (set of) second critical factor(s) the {\it
(set of) crucial factor(s)}; in case B of Lemma \ref{pskovb} we will
call the (set of) critical factor(s) the {\it (set of) crucial
factor(s)}.

\par Our claim is then the following:

\begin{lemma}
\label{pskovb}
\par Assume (\ref{assumpcion}) with weight $-n$, real length $\sigma$,
 $u=\Phi$ and $\sigma_1+\sigma_2$ factors
$\nabla^{(m)}R_{ijkl},S_{*}\nabla^{(\nu)}R_{ijkl}$, and additionally assume
 that no $\mu$-tensor field in (\ref{assumpcion}) has special free indices;
assume also that $L^{*}_\mu\bigcup L^{+}_\mu\bigcup
L''_{+}=\emptyset$ (in the notation of the extra claims above).
 Recall the cases A, B that we have distinguished above.

\par Consider case A: Recall that $k$ stands
 for the (universal) number of second critical factors
among the tensor fields indexed in $\bigcup_{z\in Z'_{Max}}L^z$.
Recall also that for each $z\in Z'_{Max}$ $\alpha$ is the number
of free indices in the (each) second critical factor. We claim that:

\begin{equation}
\label{esmen}
\begin{split}
&{\alpha\choose{2}}\Sum_{z\in Z'_{Max}}\Sum_{l\in L^z} a_l
\Sum_{r=0}^{k-1}Xdiv_{i_2}\dots Xdiv_{i_{*}}\dot{C}^{l,i_1\dots
\hat{i}_{r\alpha+1}\dots i_\mu,i_{*}}_{g}
(\Omega_1,\dots,\Omega_p,\phi_1,\dots ,\phi_u)
\\&\nabla_{i_{r\alpha+2}}\phi_{u+1}
+\Sum_{\nu\in N} a_\nu Xdiv_{i_2}\dots
Xdiv_{i_\mu}C^{\nu,i_1\dots
,i_\mu}_{g}(\Omega_1,\dots,\Omega_p,\phi_1,\dots
,\phi_u)\nabla_{i_1} \phi_{u+1}+
\\& \Sum_{t\in T_1} a_t Xdiv_{i_1}\dots Xdiv_{i_{z_t}}
C^{t, i_1\dots i_{z_t}}_{g}(\Omega_1,\dots ,\Omega_p,\phi_1,\dots
,\phi_{u+1})+
\\& \Sum_{t\in T_2}  a_t Xdiv_{i_2}\dots Xdiv_{i_{z_t}}
C^{t, i_1\dots i_{z_t}}_{g}(\Omega_1,\dots ,\Omega_p,\phi_1,\dots
,\phi_u)\nabla_{i_1}\phi_{u+1}+
\\& \Sum_{t\in T_3}  a_t Xdiv_{i_1}\dots Xdiv_{i_{z_t}}
C^{t, i_1\dots i_{z_t}}_{g}(\Omega_1,\dots ,\Omega_p,\phi_1,\dots
,\phi_{u+1})
\\& \big{(}+\Sum_{t\in T_4}  a_t Xdiv_{i_1}\dots Xdiv_{i_{z_t}}
C^{t, i_1\dots i_{z_t}}_{g}(\Omega_1,\dots ,\Omega_p,\phi_1,\dots
,\phi_{u+1})\big{)}=
\\& \Sum_{j\in J}  a_j C^j_{g}(\Omega_1,\dots
,\Omega_p,\phi_1,\dots ,\phi_{u+1})=0,
\end{split}
\end{equation}
modulo complete contractions of length $\ge\sigma +u+2$. Here each
$C^{\nu,i_1\dots i_\mu}_{g}$ is acceptable and has a simple
character $\vec{\kappa}^{+}_{simp}$ and a double character that is
 doubly subsequent to each $\vec{L}^{z,\sharp}, z\in Z'_{Max}$.
$$\Sum_{t\in T_1} a_t
C^{t, i_1\dots i_{z_t}}_{g}(\Omega_1,\dots ,\Omega_p,\phi_1,\dots
,\phi_{u+1})$$
 is a generic linear combination of acceptable tensor fields
 with  a $(u+1)$-simple character $\vec{\kappa}^{+}_{simp})$,
 and with $z_t\ge \mu$.

$$\Sum_{t\in T_2} a_t
C^{t, i_1\dots i_{z_t}}_{g}(\Omega_1,\dots ,\Omega_p,\phi_1,\dots
,\phi_{u})$$
 ($z_t\ge \mu +1$) is a generic linear combination of acceptable
 tensor fields with a $u$-simple character $\vec{\kappa}_{simp}$, with the additional
restriction that the free index ${}_{i_1}$ that belongs to the
(a) crucial factor\footnote{i.e. the second critical factor, in this
case} is a special free index.\footnote{Recall that a special free
index is either an index ${}_k,{}_l$ in a factor
$S_{*}\nabla^{(\nu)}R_{ijkl}$ or an internal index in a factor
$\nabla^{(m)}R_{ijkl}$.}

\par Now, $\Sum_{t\in T_2} a_t
C^{t, i_1\dots i_{z_t}}_{g}(\Omega_1,\dots ,\Omega_p,\phi_1,\dots
,\phi_{u+1})$ is a generic linear combination of acceptable
 tensor fields with $(u+1)$-simple character $\vec{\kappa}^{+}_{simp}$ and
$z_t\ge \mu$,\footnote{If $z_t=\mu$ then we additionally claim that
$\nabla\phi_{u+1}$ is contracting against a derivative index, and
if it is contracting against a factor $\nabla^{(B)}\Omega_h$ then
$B\ge 3$; moreover, in this case $C^{t,i_1\dots i_\mu}_g$ will
contain no special free indices.}
and moreover one unacceptable factor
$\nabla\Omega_h$ which does not contract against any factor
$\nabla\phi_t$.

\par The sublinear combination $\Sum_{t\in T_4}\dots$ appears
only if the second critical factor is of the form
$\nabla^{(B)}\Omega_k$, for some $k$. In that case, $t\in T_4$
means that there is one unacceptable factor $\nabla\Omega_k$, and
it is contracting against a factor $\nabla\phi_r$:
$\nabla_i\Omega_k\nabla^i\phi_r$, and moreover if $z_t=\mu$ then
one of the free indices ${}_{i_1},\dots,{}_{i_{\mu}}$ is a derivative index,
 and if it belongs to $\nabla^{(B)}\Omega_h$ then $B\ge 3$.

Finally,
$$\Sum_{j\in J}  a_j C^j_{g}(\Omega_1,\dots
,\Omega_p,\phi_1,\dots ,\phi_{u+1})$$ stands for a generic linear
combination of complete contractions that are $u$-simply subsequent to
$\vec{\kappa}_{simp}$.
\newline

\par In case B, we just claim the whole of Proposition \ref{giade}.
\end{lemma}

{\it Note:} 
 Lemmas \ref{zetajones}, \ref{pool2}, \ref{pskovb} (and also Lemmas \ref{oui}, \ref{oui2}) 
 will be proven in the final paper \cite{alexakis6} 
in this series. In the remainder of the present paper 
we will show that these three Lemmas imply the inductive step of
Proposition \ref{giade}.

\section{Proof that Proposition \ref{giade} follows from Lemmas
 \ref{zetajones}, \ref{pool2}, \ref{pskovb} (and Lemmas \ref{oui}, \ref{oui2}).}
\label{laderivation}

\subsection{Introduction}

{\bf General Discussion:} In this section we will show how the inductive step 
 of  Proposition
\ref{giade} (see the discussion in the beginning of the last
section) follows from Lemmas \ref{zetajones}--\ref{pskovb} (apart from certain {\it special
cases} where we will prove the inductive step of Proposition
\ref{giade} {\it directly}, without using Lemmas \ref{zetajones},
\ref{pool2} and \ref{pskovb}). We stress that
 in this derivation, we {\it will be} using the inductive
 assumption on Proposition \ref{giade}. We also repeat that when we  {\it prove}
 the Lemmas \ref{zetajones}--\ref{pskovb}, we will
 be using the inductive assumptions on Proposition \ref{giade}.

\par More precisely, we will show that Lemmas
\ref{zetajones}, \ref{pool2}, \ref{pskovb} imply
the inductive step of Proposition \ref{giade} by distinguishing three cases regarding
 the {\it assumption} of Proposition \ref{giade} (recall
that the assumption is equation (\ref{assumpcion})). The cases we
distinguish are based on the maximal refined double characters
among the $\mu$-tensor fields in (\ref{assumpcion}):

\par Recall that the $(u,\mu)$-refined double characters $\vec{L}^z, z\in Z'_{Max}$ are among the maximal
 $(u,\mu)$-refined double characters in (\ref{assumpcion}). We
have then distinguished cases I,II,III as follows:
If for any $\vec{L}^z, z\in Z'_{Max}$
 there is a special free index in some factor $S_{*}\nabla^{(\nu)}R_{ijkl}$,
then we declare that (\ref{assumpcion}) falls under case I of
Proposition \ref{giade}.\footnote{Observe that if this property
holds for one of the maximal refined double characters $\vec{L}^z,
z\in Z'_{Max}$, it will then hold for all of them.} If for
$\vec{L}^z, z\in Z'_{Max}$ there are no special free indices in
any factor of the form $S_{*}\nabla^{(\nu)}R_{ijkl}$
 but there are special free indices in
 factors of the form $\nabla^{(m)}R_{ijkl}$,
then we declare that (\ref{assumpcion}) falls under case II of
Proposition \ref{giade}.\footnote{The observation of the above
footnote still holds.} Finally, if there are no special free
indices at all in any $\vec{L}^z, z\in Z'_{Max}$, then we declare
that (\ref{assumpcion}) falls under case III of Proposition
\ref{giade}.\footnote{The observation of the above footnote still
holds.}

\par In the remainder of this paper we will show that in case I, Lemma \ref{zetajones}
implies Proposition \ref{giade}. In case II, Lemma \ref{pool2}
implies Proposition \ref{giade}, while in case III 
Lemma \ref{pskovb} (and Lemmas \ref{oui}, \ref{oui2}) implies Proposition \ref{giade}.
\newline

\par More precisely, we will show that in the setting of
each of the Lemmas \ref{zetajones}, \ref{pool2}, \ref{pskovb}
it follows that for each $z\in Z'_{Max}$ there
is a linear combination of acceptable $(\mu+1)$-tensor fields
(indexed in $P$ below) with a $(u,\mu)$-double character
$\vec{L}^z$ 
so that:

\begin{equation}
\label{gindy}
\begin{split}
&\Sum_{l\in L^z} a_l C^{l,i_1\dots
i_\mu}_{g}(\Omega_1,\dots\Omega_p,\phi_1,\dots
,\phi_u)\nabla_{i_1}\upsilon\dots \nabla_{i_\mu}\upsilon
\\&-Xdiv_{i_{\mu+1}}\Sum_{p\in P} a_p C^{p,i_1\dots
i_{\mu+1}}_{g}(\Omega_1,\dots\Omega_p,\phi_1,\dots
,\phi_u)\nabla_{i_1}\upsilon\dots \nabla_{i_\mu}\upsilon=
\\&\Sum_{t\in T} a_t C^{t,i_1\dots i_\mu}_{g}(\Omega_1,\dots\Omega_p,\phi_1,\dots
,\phi_u)\nabla_{i_1}\upsilon\dots \nabla_{i_\mu}\upsilon,
\end{split}
\end{equation}
where each $C^{t,i_1\dots i_\mu}_{g}$ is subsequent (simply or
doubly) to $\vec{L}^z$.

\par Let us just observe how (\ref{gindy}) implies 
Proposition \ref{giade}: Firstly, (\ref{gindy}) shows us that the
conclusion of Proposition \ref{giade} holds for the sublinear
combination indexed in $\bigcup_{z\in Z'_{Max}}L^z\subset \bigcup
_{z\in Z_{Max}}L^z$. But then we only have to make the
$\nabla\upsilon$'s into $Xdiv$'s in the above\footnote{See the 
last Lemma in the Appendix of \cite{alexakis1}.} and substitute back
into (\ref{hypothese2}) and we will be reduced to proving our
Proposition \ref{giade} assuming an equation:

\begin{equation}
\label{niaaal}
\begin{split}
 &\sum_{z\in Z_{Max}\setminus Z'_{Max}}\Sum_{l\in L^z} a_l Xdiv_{i_1}\dots
Xdiv_{i_\mu} C^{l,i_1\dots i_\mu}_{g}(\Omega_1,\dots
,\Omega_p,\phi_1,\dots ,\phi_u)+
\\&\Sum_{t\in T} a_t Xdiv_{i_1}\dots
Xdiv_{i_\mu} C^{t,i_1\dots i_\mu}_{g}(\Omega_1,\dots
,\Omega_p,\phi_1,\dots ,\phi_u)+
\\& \Sum_{l\in L_{\beta>\mu}} a_l Xdiv_{i_1}\dots
Xdiv_{i_\beta} C^{l,i_1\dots i_\beta}_{g}(\Omega_1,\dots
,\Omega_p,\phi_1,\dots ,\phi_u)=
\\&\Sum_{j\in J} a_j C^j_{g}(\Omega_1,\dots
,\Omega_p,\phi_1,\dots ,\phi_u),
\end{split}
\end{equation}
where the tensor fields indexed in $T$ are acceptable $\mu$-tensor
fields which are (doubly) subsequent to the tensor fields in the
first line. But then our Proposition \ref{giade} follows by
induction, since there are finitely many $(u,\mu)$-refined double
characters.
\newline

{\bf Technical discussion of the difficulties in deriving
(\ref{gindy}) from the Lemmas \ref{zetajones}, \ref{pool2}, \ref{pskovb}:} 
Let us observe that {\it at a rough level}, it
would seem that the conclusions of Lemmas \ref{zetajones},
\ref{pool2} and \ref{pskovb} would fit into the inductive
assumption of Proposition \ref{giade} because the wegiht 
is $-n$, the real length of the terms in $\sigma$ but we have {\it increased} the number $\Phi$ of 
factors $\nabla\phi,\nabla\tilde{\phi},\nabla\phi'$. Hence, if that were true,
one could hope that a direct application of Corollary
\ref{corgiade}  to the conclusions of these Lemmas would
 imply the equation (\ref{gindy}).
Unfortunately, this is not quite the case, for the reasons we will explain in the next three paragraphs.
Therefore, there is some manipulation to be done with the conclusions of Lemmas
\ref{zetajones}, \ref{pool2} and \ref{pskovb} in order to be able to apply to them  the
inductive assumption of Proposition \ref{giade} (and hence also of Corollary \ref{corgiade}),
 and this manipulation will be done in the remainder of this paper.
 The obstacles to directly applying the inductive assumption of
Proposition \ref{giade} to the conclusions of Lemmas
 \ref{zetajones}, \ref{pool2} and \ref{pskovb} are as follows:
\newline

{\it Lemma \ref{zetajones}:} Here the $(\mu-1)$-tensor fields
$C^{l,i_1\dots i_\mu}_{g}\nabla_{i_1}\phi_{u+1}$ in equation
(\ref{narod}) have the factor $\nabla\phi_{u+1}$ contracting
against the index ${}_k$ of a factor
$S_{*}\nabla^{(\nu)}R_{ijkl}$. Thus, they are {\it not} of the
form (\ref{form2}). Therefore, the inductive assumption of
Proposition \ref{giade} cannot be directly applied to
(\ref{narod}).
\newline

{\it Lemma \ref{pool2}:} Here the $(\mu-1)$-tensor fields in (\ref{vecFskillb}) {\it
are} acceptable in the form (\ref{form2}), and the inductive assumption of Proposition
\ref{giade} {\it can} be applied to (\ref{vecFskillb}).
Nonetheless, if we directly apply the inductive assumption of
Proposition \ref{giade} to to (\ref{vecFskillb}), we will obtain
an equation similar to (\ref{bengreen}), but involving a linear
combination

$$\Sum_{z\in Z'_{Max}}\Sum_{l\in L^z} a_l \Sum_{i_h\in I_{*,l}}\tilde{C}^{l,i_1\dots
i_\mu}_{g}(\Omega_1,\dots ,\Omega_p,\phi_1,\dots
,\phi_u)\nabla_{i_h}\phi_{u+1}\nabla_{i_1}\upsilon\dots\hat{\nabla}_{i_h}\upsilon
\dots\nabla_{i_\mu}\upsilon$$ rather than a linear
combination

$$\Sum_{z\in Z'_{Max}}\Sum_{l\in L^z} a_l \Sum_{i_h\in I_{*,l}}C^{l,i_1\dots
i_\mu}_{g}(\Omega_1,\dots ,\Omega_p,\phi_1,\dots
,\phi_u)\nabla_{i_1}\upsilon\dots\nabla_{i_\mu}\upsilon$$ as
required (the important difference here is the symbol $\tilde{ }$,
which stands for an $S_{*}$-symmetrization). In other words,  in
$\tilde{C}^{l,i_1\dots i_\mu}_{g}(\Omega_1,\dots
,\Omega_p,\phi_1,\dots ,\phi_u)\nabla_{i_h}\phi_{u+1}$ some factor
 $\nabla^{(m)}R_{ijkl}$ has been $S_{*}$-symmetrized.
It is then not obvious how to manipulate this equation to obtain (\ref{bengreen}).
\newline

{\it Lemma \ref{pskovb}:} In this case there are numerous obstacles to deriving the inductive
 step of Proposition \ref{giade} from (\ref{esmen}).
 Firstly, the tensor fields indexed in $T_3, T_4$ are {\it not}
 acceptable. Secondly, even if these index sets
  were empty, the tensor fields indexed in $T_2$ {\it do not}
have the $(u+1)$-simple character $\vec{\kappa}_{simp}^{+}$ of the
tensor fields in the first line of (\ref{esmen}). Lastly, even if
this index set $T_2$ were also empty, and we directly applied the
inductive
 assumption of Proposition \ref{giade} to (\ref{esmen}), we
  would obtain a statement involving the expression:
\begin{equation}
\begin{split} 
&{\alpha\choose{2}}\Sum_{z\in Z'_{Max}}\Sum_{l\in L^z} a_l
\Sum_{r=0}^{k-1}\dot{C}^{l,i_1\dots \hat{i}_{r\alpha+1}\dots
i_\mu,i_{*}}_{g} (\Omega_1,\dots,\Omega_p,\phi_1,\dots ,\phi_u)
\\&\nabla_{i_{r\alpha+2}}\phi_{u+1}\nabla_{i_2}\upsilon\dots\nabla_{i_{*}}\upsilon,
\end{split}
\end{equation}
 and this expression is quite different from the expression we need in (\ref{gindy}).

\subsection{Derivation of Proposition \ref{giade} in case I from Lemma
\ref{zetajones}.}

\par We start this subsection with a technical Lemma that will be needed here, but will
also be used on multiple occasions throughout this
work:

\begin{lemma}
\label{petermichel}
 Let $\sum_{x\in X} a_x C^{x,i_1\dots , i_\beta}_{g}(\Omega_1,\dots ,
\Omega_p,\phi_1,\dots ,\phi_b)$ stand for  a linear combination of
tensor fields, each with rank $ \beta$,
 for some given number $\beta$, and with a given simple character
 $\vec{\kappa}^{*}_{simp}$,
 with real length $\sigma\ge 4$ and weight $-n$.
 We also assume that there is a  given factor
$\nabla^{(y)}\Omega_c$, $y\ge 1$, ($c$ independent of $x$) in each $C^{x,i_1\dots i_\beta}_{g}$
 all of whose indices are contracting against factors $\nabla\phi$.\footnote{If $y\ge 2$
 then our tensor fields are assumed to be acceptable.
 If $y=1$ then we assume that $\nabla\Omega_c$ is the only unacceptable factor.}
We assume an equation:

\begin{equation}
\label{pollhdoul} 
\begin{split}
&\sum_{x\in X} a_x X_{*}div_{i_1}\dots
X_{*}div_{i_\beta}C^{x,i_1\dots  i_\beta}_{g}(\Omega_1,\dots ,
\Omega_p,\phi_1,\dots ,\phi_b)
\\&+\sum_{j\in J} a_j
C^j_{g}(\Omega_1,\dots , \Omega_p,\phi_1,\dots ,\phi_b)=0,
\end{split}
\end{equation}
where $X_{*}div_i$ stands for the sublinear combination in
$Xdiv_i$ for which $\nabla_i$ is not allowed to hit the factor
$\nabla^{(y)}\Omega_c$. The complete contractions in $J$ are
simply subsequent to $\vec{\kappa}^{*}_{simp}$. We
additionally assume that if we formally erase the factor
$\nabla^{(y)}\Omega_c$ along with the factors $\nabla\phi_h$ that
it is contracting against, then none of the resulting
$\beta$-tensor fields is ``forbidden'' in the sense of Definition \ref{forbidden}.
\newline

\par We claim that we can then write:

\begin{equation}
\label{pollhdoul'} \begin{split} & \sum_{x\in X} a_x
Xdiv_{i_1}\dots Xdiv_{i_\beta}C^{x,i_1\dots
i_\beta}_{g}(\Omega_1,\dots , \Omega_p,\phi_1,\dots ,\phi_b)=
\\&\sum_{x\in X'} a_x
Xdiv_{i_1}\dots Xdiv_{i_\beta}C^{x,i_1\dots
i_\beta}_{g}(\Omega_1,\dots , \Omega_p,\phi_1,\dots ,\phi_b)+
\\&\sum_{j\in J} a_j
C^j_{g}(\Omega_1,\dots , \Omega_p,\phi_1,\dots ,\phi_b),
\end{split}
\end{equation}
where the tensor fields indexed in $X'$ are exactly like the ones
indexed in $X$, only the chosen factor $\nabla^{(y)}\Omega_c$ has
at least one index that is not contracting against a factor
$\nabla\phi$. The complete contractions in $J$ are simply
subsequent to
 $\vec{\kappa}^{*}_{simp}$.
\end{lemma}

{\it (Sketch of the) Proof of Lemma \ref{petermichel}:} We just apply the eraser to the
factor $\nabla^{(y)}\Omega_c$ and the factors $\nabla\phi_h$ that
it is contracting against in (\ref{pollhdoul}), obtaining a new
true equation. We can then iteratively apply Corollary
\ref{corgiade} to this new true equation, multiplying by
$\nabla^{(B)}_{r_1\dots
r_B}\Omega_c\nabla^{r_1}\upsilon\dots\nabla^{r_B}\upsilon$ and
making all the factors $\nabla\upsilon$ into $Xdiv$'s at each
stage. This would show our claim except for the  caveat that
in the last step of the above iteration, we might not be able to
apply Corollary \ref{corgiade} if the tensor fields of maximal
refined double character are in one of the forbidden forms of
Corollary \ref{corgiade} with rank $>\mu$. In that case, in the last
step we use Lemma \ref{appendix} below (setting
$\Phi=\nabla^{(y)}_{r_1\dots r_c}\Omega_c\nabla^{r_1}
\phi_{h_1}\dots\nabla^{r_y}\phi_{h_y}$). That concludes the proof
of our claim in this case. $\Box$
\newline

Furthermore, we have a weaker form of Lemma \ref{petermichel} when $\sigma=3$.
We firstly introduce a definition that will be used on a number of
occasions below:

\begin{definition}
\label{proextremovable} Consider any tensor field in the form
(\ref{form2}). We consider any set of indices,
$\{{}_{x_1},\dots,{}_{x_s}\}$ belonging to a factor $T$ (here $T$ is not in the form $\nabla\phi$). We
assume that these indices are neither free nor are contracting
against a factor $\nabla\phi_h$.

If the indices belong to a factor $T$ in the form
$\nabla^{(B)}\Omega_1$ then $\{{}_{x_1},\dots,{}_{x_s}\}$ are
removable provided $B\ge s+2$.

\par Now, we consider indices that belong to a factor
$\nabla^{(m)}R_{ijkl}$ (and are neither free nor are contracting
against a factor $\nabla\phi_h$). Any such index ${}_x$ which is a
derivative index will be removable. Furthermore, if $T$ has at
least two free derivative indices, then if neither of the indices
${}_i,{}_j$ are free then we will say one of ${}_i,{}_j$ is removable;
accordingly, if neither of ${}_k,{}_l$ is free then we will say
that one of ${}_k,{}_l$ is removable. Moreover, if $T$ has one free
derivative index then: if none of the indices ${}_i,{}_j$ are free
then we will say that one of the indices ${}_i,{}_j$ is removable;
on the other hand if one of the indices ${}_i,{}_j$ is also free
and none of the indices ${}_k,{}_l$ are free then we will say that
one of the indices ${}_k,{}_l$ is removable.

\par Now, we consider a set of indices $\{{}_{x_1},\dots,{}_{x_s}\}$ that belong to
a factor $T=S_{*}\nabla^{(\nu)}R_{ijkl}$ and are not special, and
are not free and are not contracting against any $\nabla\phi$. We
will say this set of indices is removable if $s\le\nu$.
Furthermore, if none of the indices ${}_k,{}_l$ are free and
$\nu>0$ and at least one of the other indices in $T$ is free, we
will say that one of the indices ${}_k,{}_l$ is removable.
\end{definition}

\par Weaker version of Lemma \ref{petermichel} when $\sigma=3$:

\begin{lemma}
\label{petermichel3} Assume the equation (\ref{pollhdoul}) when
$\sigma=3$ and assume additionally that every tensor field indexed
in $X$ has a removable index. Then (\ref{pollhdoul'}) still holds.
\end{lemma}

{\it Proof:} The argument essentially follows the ideas developed
in the paper \cite{alexakis3}. Firstly, we observe that (possibly
 after applying the second Bianchi identity) we can
explicitly write:

\begin{equation}
\label{grafeligo} \begin{split} & \sum_{x\in X} a_x
Xdiv_{i_1}\dots Xdiv_{i_\beta}C^{x,i_1\dots 
i_a}_{g}(\Omega_1,\dots , \Omega_p,\phi_1,\dots ,\phi_b)=
\\&\sum_{x\in X'} a_x
Xdiv_{i_1}\dots Xdiv_{i_\beta}C^{x,i_1\dots 
i_\beta}_{g}(\Omega_1,\dots , \Omega_p,\phi_1,\dots ,\phi_b)+
\\&\sum_{x\in \overline{X}} a_x
Xdiv_{i_1}\dots Xdiv_{i_\gamma}C^{x,i_1\dots 
i_{\gamma}}_{g}(\Omega_1,\dots , \Omega_p,\phi_1,\dots ,\phi_b)+
\\&\sum_{j\in J} a_j
C^j_{g}(\Omega_1,\dots , \Omega_p,\phi_1,\dots ,\phi_b),
\end{split}
\end{equation}
where the tensor fields indexed in $\overline{X}$ have all the
properties of the ones indexed in $X$, only they all have rank
$\gamma\ge\beta+1$, and they also have no removable indices. The
sublinear combination $\sum_{x\in X'}\dots$ (here and below, when
it appears on the RHS) stands for a {\it generic} linear
combination as described in Lemma \ref{petermichel}.

\par We then observe that we can write:

\begin{equation}
\label{grafeligo'} \begin{split} & \sum_{x\in \overline{X}} a_x
Xdiv_{i_1}\dots Xdiv_{i_\gamma}C^{x,i_1\dots 
i_\gamma}_{g}(\Omega_1,\dots , \Omega_p,\phi_1,\dots ,\phi_b)=
\\&(Const)_{*} Xdiv_{i_1}\dots Xdiv_{i_\gamma} C^{*,i_1\dots
i_\gamma}_g(\Omega_1,\dots , \Omega_p,\phi_1,\dots ,\phi_b)+
\\&\sum_{x\in X'} a_x
Xdiv_{i_1}\dots Xdiv_{i_\beta}C^{x,i_1\dots 
i_\beta}_{g}(\Omega_1,\dots , \Omega_p,\phi_1,\dots ,\phi_b)+
\\&\sum_{j\in J} a_j
C^j_{g}(\Omega_1,\dots , \Omega_p,\phi_1,\dots ,\phi_b)
\end{split}
\end{equation}
where the tensor field $C^{*,i_1\dots i_\gamma}_g$ is zero
 unless $\sigma_1=\sigma_2=0$ or $\sigma_1=2$
 or $\sigma_2=2$. In those cases, the tensor field
 $C^{*,i_1\dots
i_\gamma}_g(\Omega_1,\dots , \Omega_p,\phi_1,\dots ,\phi_b)$ is, respectively:

\begin{equation}
\begin{split}
\label{foula1}
&pcontr(\nabla^{(X)}_{i_1\dots i_a u_1\dots u_t}\Omega_1\otimes
\nabla^{(B)}_{j_1\dots j_by_1\dots
y_r}\Omega_2\otimes\nabla^{u_1}\phi_1\otimes\dots
\nabla^{j_b}\phi_f \otimes\nabla^{(B)}_{z_1\dots z_q}\Omega_3
\\&\otimes\nabla^{z_1}\phi_{f+1}\otimes\nabla^{z_q}\phi_{u+1}),
\end{split}
\end{equation}
(here if $y\ge 2$ then $b=0$; if $y\le 1$ then $y=2-b$),

\begin{equation}
\begin{split}
\label{foula2}
&pcontr(\nabla^{(X)}_{i_1\dots i_a u_1\dots u_t}R_{i_{a+1}ji_{a+2}l}\otimes
\nabla^{(r)}_{y_1\dots
y_r}{{{R_{i_{a+3}}}^j}_{i_{a+4}}}^l\otimes\nabla^{u_1}\phi_1\otimes\dots
\nabla^{j_b}\phi_f\otimes
\\&\nabla^{(B)}_{z_1\dots z_q}\Omega_3\otimes
\nabla^{z_1}\phi_{f+1}\otimes\nabla^{z_q}\phi_{u+1}),
\end{split}
\end{equation}

\begin{equation}
\begin{split}
\label{foula3}
&pcontr(S_{*}\nabla^{(X)}_{i_1\dots i_a u_1\dots u_t}R_{ii_{a+1}i_{a+2}l}\otimes
\nabla^{(r)}_{y_1\dots
y_r}{{R_{i'i_{a+3}}}_{i_{a+4}}}^l\otimes\nabla^i\tilde{\phi}_1\otimes\nabla^{i'}\tilde{\phi}_2\otimes
\\&\nabla^{u_1}\phi_3\otimes\dots
\nabla^{j_b}\phi_f
\otimes\nabla^{(B)}_{z_1\dots z_q}\Omega_3
\nabla^{z_1}\phi_{f+1}\otimes\nabla^{z_q}\phi_{u+1}).
\end{split}
\end{equation}

\par Then, picking out the sublinear combination in
(\ref{grafeligo'}) with only factors $\nabla\phi$ contracting
against $\nabla^{(B)}\Omega_c$ we derive that $(Const)_{*}=0$.
$\Box$
\newline

{\bf Derivation of Proposition \ref{giade} (in case I) from Lemma \ref{zetajones}:}
\newline

{\it Special cases etc:} Now, we return to the derivation of Proposition \ref{giade}
(in case I) from Lemma \ref{zetajones}. We will be singling out a further 
case, which we will call ``delicate''. In this ``delicate'' case we will
derive Proposition \ref{giade} from Lemma \ref{zetajones} by using an extra 
Lemma (see Lemma \ref{delicate} below). The proof of Lemma \ref{delicate} 
will be provided in the Appendix to this paper.  

{\it The ``delicate case'':} The delicate case is when the 
$\mu$-tensor fields of maximum refined double character in our 
Lemma assumption have no removable free indices, and moreover 
their critical factor is in the form: 
$S_{*}\nabla^{(\nu)}_{r_1\dots r_\nu}R_{i r_{\nu+1}i_1l}$, 
where {\it all} indices ${}_{r_1},\dots ,{}_{r_\nu},{}_{r_{\nu+1}}$ 
are either free or contracting against a factor $\nabla\phi'_h$.\footnote{Notice 
that by weight considerations and by the definition of maximal refined double 
character, if one tensor field of maximal refined double 
character has this property then all of them will.} 

In that case we have an extra claim, which we will prove in the Appendix:
\begin{lemma}
\label{delicate} 
For each $z\in Z'_{Max}$, we let $L^z_*\subset L^z$
stand for the index set of tensor fields $C^{l,i_1\dots i_\mu}_g, l\in L^z$
for which the index ${}_l$ in the critical factor contracts against a 
special index in some factor $S_{*}R_{ijkl}$.  

We claim that for each $z\in Z'_{Max}$,we can write:

\begin{equation}
\label{delicateeqn} 
\begin{split}
&\sum_{l\in L^z_*} a_l C^{l,i_1\dots i_\mu}_g\nabla_{i_1}\upsilon\dots
 \nabla_{i_\mu}\upsilon =
\sum_{l\in L'^z} a_l C^{l,i_1\dots i_\mu}_g\nabla_{i_1}\upsilon\dots
 \nabla_{i_\mu}\upsilon.
\end{split}
\end{equation}
Here the terms indexed in $L'^z$ have all the 
properties of the terms indexed in $L^z$,
 but in addition the index ${}_l$ in the critical factor does 
not contract against a special index in a factor $S_{*}R_{ijkl}$. 
\end{lemma}

{\it The derivation of Proposition \ref{giade} (in case I):}
\newline

Recall the conclusion of Lemma \ref{zetajones}, equation (\ref{narod}).
Recall that for each tensor field and each complete
 contraction in (\ref{narod}),
$\nabla\phi_{u+1}$ is contracting against the crucial factor,
which was defined to
 be the factor $S_{*}\nabla^{(\nu)}R_{ijkl}$
  in $\vec{\kappa}_{simp}$ whose index ${}_i$  is contracting against a
chosen factor $\nabla\tilde{\phi}_{Min}$.
For notational convenience, we will assume that $Min=1$, i.e. that
that index ${}_i$ in the crucial factor is contracting against a
factor $\nabla\tilde{\phi}_1$.

\par Define the set $Stan$ to stand for the set of numbers $o$ for which the factor
 $\nabla\phi'_o$ is contracting against one of the indices
${}_{r_1},\dots, ,{}_{r_{\nu+1}}$ in the crucial factor
$S_{*}\nabla^{(\nu)}_{r_1\dots r_\nu}R_{ir_{\nu+1}kl}$.
 With no loss of generality, we assume that $Stan=\{2,\dots ,q\}$ or
$Stan=\emptyset$ (which is equivalent to saying
 $q=1$--we will be using that convention below).

\par For convenience, we will assume that for
each of the tensor fields appearing in (\ref{narod})
the factors $\nabla\phi'_2,\dots,\nabla\phi'_q$ are contracting against the
 indices ${}_{r_1},\dots ,{}_{r_{q-1}}$ in the crucial factor
$S_{*}\nabla^{(\nu)}_{r_1\dots r_\nu}R_{ir_{\nu+1}kl}$.

\par With this convention, we introduce a new definition:

\begin{definition}
\label{ypsilanti}
 For each tensor
field  $C^{l,i_1\dots i_\mu}_{g}(\Omega_1,\dots
,\Omega_p,\phi_1,\dots ,\phi_u)\nabla_{i_1}\phi_{u+1}$,
$C^{\nu,i_1\dots i_\mu}_{g}(\Omega_1,\dots ,\Omega_p,\phi_1,\dots
,\phi_u)\nabla_{i_1}\phi_{u+1}$, $C^{p,i_1,\dots
,i_{\mu+1}}_{g}(\Omega_1,\dots ,\Omega_p,\phi_1,\dots
,\phi_u,\phi_{u+1})$ in (\ref{narod}), we define
$C^{l,i_2\dots i_\mu}_{g}(\Omega_1,\dots ,\Omega_p,Y,\phi_2,\dots
,\phi_u)$, \\$C^{\nu,i_2\dots i_\mu}_{g}(\Omega_1,\dots
,\Omega_p,Y,\phi_2,\dots ,\phi_u)$,
$C^{p,i_2,\dots ,i_{\mu+1}}_{g}(\Omega_1,\dots
,\Omega_p,Y,\phi_2,\dots ,\phi_u)$ to stand for the
 tensor fields that arise by formally replacing the  expression
\\$S_{*}\nabla^{(\nu)}_{r_1\dots
r_{\nu}}R_{ir_{\nu+1}kl}\nabla^i\tilde{\phi}_1\nabla^k
\phi_{u+1}$ by a factor $\nabla^{(\nu+2)}_{r_{1}\dots r_\nu
r_{\nu+1}l}Y$ ($Y$ is a scalar function).
\end{definition}

\par We observe that the tensor fields we are left with have
 length $\sigma+u-1$, a factor $\nabla^{(B)} Y$ with
$B\ge 2$, and are acceptable if we set $Y=\Omega_{p+1}$.
 We observe that all these tensor fields have the same $(u-1)$-simple character
(where we treat the function $Y$ as a function $\Omega_{p+1}$),
which we denote by $\tilde{\kappa}_{simp}$.

 We also note that each of the tensor fields
$C^{l,i_2\dots i_\mu}_{g}(\Omega_1,\dots ,\Omega_p,Y,\phi_2,\dots
,\phi_u)$, $C^{\nu,i_2\dots i_\mu}_{g}(\Omega_1,\dots
,\Omega_p,Y,\phi_2,\dots ,\phi_u)$, $C^{p,i_2\dots
i_b}_{g}(\Omega_1,\dots ,\Omega_p,Y,\phi_2,\dots ,\phi_u)$
 will have the property that
the last index ${}_l$ in $\nabla^{(B)}Y$ (in the tensor 
field \\$C^{l,i_2\dots i_\mu}_{g}(\Omega_1,\dots ,\Omega_p,Y,\phi_2,\dots
,\phi_u)$) is neither free nor
contracting against a factor $\nabla\phi_h$: This is because the
last index ${}_l$ in $\nabla^{(B)}Y$ corresponds to the index
${}_l$ in the crucial factor $S_{*}\nabla^{(\nu)}_{r_1\dots
r_{\nu}}R_{ir_{\nu+1}kl}\nabla^i\tilde{\phi}_1\nabla^k
\phi_{u+1}$ of the tensor field $C^{l,i_1\dots i_\mu}_{g}(\Omega_1,\dots ,\Omega_p,\phi_1,\dots
,\phi_u)$, and this index is neither
 free nor contracting against any factor $\nabla\phi'_o$ by hypothesis.

We claim an equation:

\begin{equation}
\label{Ig1}
\begin{split}
&\Sum_{z\in Z'_{Max}}\Sum_{l\in L^z} a_l Xdiv_{i_2}\dots
Xdiv_{i_\mu}C^{l,i_2\dots i_\mu}_{g}(\Omega_1,\dots
,\Omega_p,Y,\phi_2,\dots ,\phi_u)+ \\&\Sum_{\nu\in N} a_\nu
Xdiv_{i_2}\dots Xdiv_{i_\mu}C^{l,i_2\dots
i_\mu}_{g}(\Omega_1,\dots ,\Omega_p,Y,\phi_2,\dots ,\phi_u)-
\\&
\Sum_{p\in P} a_p Xdiv_{i_2}\dots Xdiv_{i_{\mu+1}}
C^{p,i_2\dots i_{\mu+1}}_{g}(\Omega_1,\dots
,\Omega_p,Y,\phi_2,\dots ,\phi_u)=
\\&\Sum_{t\in T'} a_t C^t_{g}
(\Omega_1,\dots ,\Omega_p,Y,\phi_2,\dots ,\phi_u),
\end{split}
\end{equation}
which will hold modulo complete contractions of length $\ge \sigma
+u$. Here the right hand side stands for a generic linear
combination of complete contractions that are
simply subsequent to $\tilde{\kappa}_{simp}$.

{\it Proof of (\ref{Ig1}):} Since the argument by which we derive this equation will
be used frequently in this series of papers, we codify
this claim in a Lemma:

\begin{lemma}
\label{technical1}
Consider a linear combination of acceptable $(\gamma+1)$-tensor fields,
 $\Sum_{x\in X} a_x C^{x,i_1\dots ,i_{\gamma+1}}_{g}(\Omega_1,\dots,\Omega_p,\phi_1,\dots,\phi_u)$,
all in the form (\ref{form2}) with weight $-n$ and with
 a given simple character $\overline{\kappa}_{simp}$.
 Assume that for each of the tensor fields then the index ${}_{i_1}$ is the index ${}_k$ in a given factor
$S_{*}\nabla^{(\nu)}R_{ijkl}$, for which the index ${}_i$ is contracting against a chosen factor
 $\nabla\tilde{\phi}_w$ (wlog we will assume $w=1$). Assume that
 $\Sum_{z\in Z} a_z C^{z,i_1\dots ,i_{\epsilon_z+1}}_{g}(\Omega_1,\dots,\Omega_p,\phi_1,\dots,\phi_u)$
is a linear combination with all the features of the
tensor fields indexed in $X$, only now each $\epsilon_z>\gamma$.
Assume an equation:

\begin{equation}
\label{mporwdenm}
\begin{split}
&\Sum_{x\in X} a_x Xdiv_{i_2}\dots Xdiv_{i_{\gamma+1}}
C^{x,i_1\dots
,i_{\gamma+1}}_{g}(\Omega_1,\dots,\Omega_p,\phi_1,\dots,\phi_u)\nabla_{i_1}\phi_{u+1}
\\&+\Sum_{z\in Z} a_z Xdiv_{i_2}\dots Xdiv_{i_{\epsilon_z}}
C^{z,i_1\dots ,i_{\epsilon_z+1}}_{g}(\Omega_1,\dots,\Omega_p,
\phi_1,\dots,\phi_u)\nabla_{i_1}\phi_{u+1}=
\\& \Sum_{j\in J} a_j C^{j,i_1}_{g}(\Omega_1,\dots,\Omega_p,
\phi_1,\dots,\phi_u)\nabla_{i_1}\phi_{u+1};
\end{split}
\end{equation}
here the vector fields in the RHS
have a $u$-weak character
$Weak(\overline{\kappa}_{simp})$ and are either simply
 subsequent to $\overline{\kappa}_{simp}$ or have one of the
two factors $\nabla\phi_1$, $\nabla\phi_{u+1}$ contracting
 against a derivative index, or both factors
$\nabla\phi_w$, $\nabla\phi_{u+1}$ are contracting against
 anti-symmetric indices ${}_i,{}_j$ or ${}_k,{}_l$ in some curvature factor.
\newline

Denote by $\overline{C}^{x,i_2\dots
,i_{\gamma+1}}_{g}(\Omega_1,\dots,\Omega_{p+1},\phi_2,\dots,\phi_u)$,
$C^{z,i_1\dots ,i_{\epsilon_z+1}}_{g}(\Omega_1,\dots,\Omega_{p+1},
\phi_2,\dots,\phi_u)$ the tensor fields that arise from
$C^{x,i_1\dots
,i_{\gamma+1}}_{g}(\Omega_1,\dots,\Omega_{p+1},\phi_2,\dots,\phi_u)$,
\\$C^{z,i_1\dots ,i_{\epsilon_z+1}}_{g}(\Omega_1,\dots,\Omega_{p+1},
\phi_2,\dots,\phi_u)$ by formally replacing the expression
\\$S_{*}\nabla^{(\nu)}_{r_1\dots r_\nu}R_{ijkl}
\nabla^i\tilde{\phi}_1\nabla^k\phi_{u+1}$ by an expression
$\nabla^{(\nu+2)}_{r_1\dots r_\nu jl}\Omega_{p+1}$.
 Denote by $\tilde{\kappa}_{simp}$ the $(u-1)$-simple character of each of
 the resulting tensor fields
(they have length $\sigma+u-1$). We then claim
that modulo complete contractions of length $\ge\sigma+u$:

\begin{equation}
\label{mporwdenm2}
\begin{split}
&\Sum_{x\in X} a_x Xdiv_{i_2}\dots Xdiv_{i_{\gamma+1}}
\overline{C}^{x,i_1\dots
,i_{\gamma+1}}_{g}(\Omega_1,\dots,\Omega_{p+1},\phi_2,\dots,\phi_u)
\\&+\Sum_{z\in Z} a_z Xdiv_{i_2}\dots Xdiv_{i_{\epsilon_z}}
\overline{C}^{z,i_1\dots
,i_{\epsilon_z+1}}_{g}(\Omega_1,\dots,\Omega_{p+1},
\phi_2,\dots,\phi_u)=
\\& \Sum_{j\in J'} a_j C^{j}_{g}(\Omega_1,\dots,\Omega_{p+1},
\phi_1,\dots,\phi_u),
\end{split}
\end{equation}
 where the complete contractions indexed in $J'$ are simply subsequent to $\tilde{\kappa}_{simp}$.
We note that the proof of this Lemma will be independent of Proposition \ref{giade}.
\end{lemma}

{\it Note:} Before we prove this Lemma we remark that by applying it
to (\ref{narod}) we derive (\ref{Ig1}).
\newline

{\it Proof of Lemma \ref{technical1}:} Denote the left hand side
of (\ref{mporwdenm}) by $F_{g}$. We then denote by $F'_{g}$ the
linear combination that arises from $F_{g}$ by formally replacing
the factors $\nabla_a\phi_1,\nabla_b\phi_{u+1}$ by $g_{ab}$ (the
uncontracted metric tensor).
 Notice that $F'_{g}$ then consists of complete contractions with one internal
contraction in a curvature factor, and with weight $-n+2$ and
length $\sigma+u$. Since $F_{g}=0$ modulo longer complete
contractions, and since this equation holds formally, we derive
that $F'_{g}=0$ modulo longer
 complete contractions. Now, apply the operation $Ricto\Omega_{p+1}$ to $F'_{g}$ 
 (see the relevant Lemma  in the Appendix of \cite{alexakis1}).
 Denote the resulting
 linear combination by $F''_{g}$. By definition of $Ricto\Omega_{p+1}$, the minimum
length of the complete contractions in $F''_{g}$ is $\sigma+u-1$.
If we denote this sublinear combination
 by $F''^{\sigma+u-1}_{g}$ then (virtue of the afformentioned Lemma) we will have $F''^{\sigma+u-1}_{g}=0$
modulo longer complete contractions. By following
all the operations we have performed we observe
that this equation is precisely (\ref{mporwdenm2}). $\Box$
\newline

\par Now, observe that (\ref{Ig1}) falls under our inductive assumption
of Proposition \ref{giade}:\footnote{Notice
that since our assumption (\ref{assumpcion}) does not include
tensor fields in any of the ``forbidden forms'', it follows that
the tensor fields of minimum rank in (\ref{Ig1}) are also {\it
not} in any of the ``forbidden forms''.} All the tensor fields
are acceptable, and they all have
 a given simple character $\tilde{\kappa}_{simp}$;
 furthermore, the weight of the complete contractions in (\ref{Ig1})
 is $-n+2>-n$. Lastly, recall that we have noted that the last index in $\nabla^{(B)}Y$ is neither
 free nor contracting against any  factor $\nabla\phi_{u+1}$.

 We observe that the sublinear
 combinations of $(\mu -1)$-tensor fields on the left hand side of
 (\ref{Ig1}) with maximal double characters are the
 sublinear combinations:

$$\Sum_{z\in Z'_{Max}} \Sum_{l\in L^z} a_l C^{l,i_2\dots
i_\mu}_{g}(\Omega_1,\dots ,\Omega_p,Y,\phi_2,\dots ,\phi_u).$$
(This follows directly from the definition of the maximal refined
double characters).
 We denote the respective refined double characters for
 the complete contractions of this form by $\vec{L^z}', z\in Z'_{Max}$.
Applying the inductive hypothesis of Corollary
\ref{corgiade} to (\ref{Ig1}) and picking out the sublinear
 combination with a $(u-1,\mu-1)$-double character $Doub(\vec{L^z}')$,
  we deduce that there is a
linear combination of acceptable $\mu$-tensor fields with a refined
double character $\vec{L^z}'$,

$$\Sum_{r\in R^z} a_r C^{r,i_2\dots
i_{\mu +1}}_{g}(\Omega_1,\dots ,\Omega_p,Y,\phi_2,\dots ,\phi_u),$$
such that:

\begin{equation}
\label{ajax}
\begin{split}
&\Sum_{l\in L^z} a_l C^{l,i_2\dots i_\mu}_{g}(\Omega_1,\dots
,\Omega_p,Y,\phi_2,\dots ,\phi_u)\nabla_{i_2}\upsilon\dots
\nabla_{i_\mu}\upsilon
\\&- \Sum_{r\in R^z} a_r Xdiv_{i_{\mu +1}}C^{r,i_2\dots
i_{\mu+1}}_{g}(\Omega_1,\dots ,\Omega_p,Y,\phi_2,\dots ,\phi_u)
\nabla_{i_2}\upsilon\dots\nabla_{i_\mu}\upsilon=
\\&\Sum_{t\in T} a_t C^{t,i_2\dots i_\mu}_{g}
(\Omega_1,\dots ,\Omega_p,Y,\phi_2,\dots
,\phi_u)\nabla_{i_2}\upsilon\dots \nabla_{i_\mu}\upsilon,
\end{split}
\end{equation}
 modulo complete contractions of greater length.
(Recall that since we are considering the factor $Y$ to be of the
form $\Omega_{p+1}$, each of the tensor fields above
has a factor $\nabla^{(b)}Y, b\ge 2$).

Here,
$$\Sum_{t\in T} a_t C^{t,i_2\dots i_\mu}_{g}
(\Omega_1,\dots ,\Omega_p,Y,\phi_1,\dots ,\phi_u)$$ stands for a
generic linear combination of tensor fields whose
 refined double character is (simply or doubly) subsequent to the refined double character
$\vec{L^z}'$. Moreover, if this tensor field $C^{t,i_2\dots i_\mu}_{g}$
is only doubly subsequent to $\vec{L^z}'$,
 then at least one of the indices in the factor $\nabla^{(B)}Y$
 is not contracting against a factor of the
 form $\nabla\upsilon$, $\nabla\phi$. (This follows since those
 tensor fields arise from the tensor fields indexed in $N$ in (\ref{Ig1}),
  which have that property by construction).

Now, we refer to (\ref{ajax}) and we make
all the factors $\nabla\upsilon$ that are {\it not} contracting
against the critical factor into an $Xdiv$ (see the last 
Lemma in the Appendix of \cite{alexakis1}). We thus obtain a new
true equation:

\begin{equation}
\label{ajaxnew}
\begin{split}
&\Sum_{l\in L^z} a_l Xdiv_{i_{M+1}}\dots Xdiv_{i_\mu}C^{l,i_2\dots
i_\mu}_{g}(\Omega_1,\dots ,\Omega_p,Y,\phi_2,\dots
,\phi_u)\nabla_{i_M}\upsilon\dots \nabla_{i_{M-1}}\upsilon=
\\& \Sum_{r\in R^z} a_r Xdiv_{i_{M+1}}\dots Xdiv_{i_{\mu +1}}C^{r,i_2\dots
i_{\mu+1}}_{g}(\Omega_1,\dots ,\Omega_p,Y,\phi_2,\dots ,\phi_u)
\nabla_{i_2}\upsilon\dots\nabla_{i_M}\upsilon
\\&+\Sum_{t\in T} a_t Xdiv_{i_{M+1}}\dots Xdiv_{i_\mu} C^{t,i_2\dots i_\mu}_{g}
(\Omega_1,\dots ,\Omega_p,Y,\phi_2,\dots
,\phi_u)\nabla_{i_2}\upsilon\dots \nabla_{i_M}\upsilon.
\end{split}
\end{equation}

\par Now, before we can proceed with our argument we will
need to employ a technical Lemma 
whose proof we will present in the end of this subsection. 

\begin{lemma}
\label{technlem1}
We claim that we may assume with no loss of
generality that for each tensor field indexed in $R^z$, one index
in $\nabla^{(B)}Y$ (the last one, wlog) is not contracting against
any factor $\nabla\phi_h$ or $\nabla\upsilon$.
\end{lemma}

\par Under this extra assumption, we can now derive the claim of Proposition \ref{giade}:

\par Now, since all the tensor fields in $L^z, R^z$ have the same $(u-1,\mu-1)$-double character,
$\vec{L^z}'$, it follows that for each tensor field appearing in
(\ref{ajax}) there is a fixed number $\tau$ of factors
$\nabla\upsilon$ contracting against
 the factor $\nabla^{(B)}Y$, and also a fixed number ($q-1$)
 of factors $\nabla\phi'_o$ contracting against $\nabla^{(B)}Y$. We may assume with no
loss of generality that the $\tau$ factors $\nabla\upsilon$  are contracting against
the first $\tau$ indices in $\nabla^{(B)}Y$ and the $q-1$ factors
$\nabla\phi_o$ are contracting against the next $q-1$ indices in
$\nabla^{(B)}Y$, in a decreasing rearrangement according to the
numbers $o$. By using the Eraser,\footnote{See the relevant 
Lemma in the Appendix of \cite{alexakis1}.} we can see that under these
assumptions, (\ref{ajax}) will hold formally, subject to the
additional feature that for each complete contraction in
(\ref{ajax}), the first $\tau+q-1$ indices are not permuted. (Call
this the extra feature).

\par In this setting, we define an operation $Replace$ that acts
on the tensor fields in (\ref{ajax}) by replacing the factor
$\nabla^{(b)}_{r_1\dots r_b}Y$ (recall $b\ge 2$) by an expression
$S_{*}\nabla^{(b-2)}_{r_1\dots r_{b-2}}R_{ir_{b-1}kr_b}\nabla^i\phi_1
\nabla^k\upsilon$. We denote the resulting tensor fields by:

$$\tilde{C}^{l,i_1\dots i_\mu}_{g}(\Omega_1,\dots ,\Omega_p,
\phi_1,\dots ,\phi_u)\nabla_{i_1}\upsilon,$$

$$\tilde{C}^{r,i_1\dots i_{\mu+1}}_{g}(\Omega_1,\dots ,\Omega_p,
\phi_1,\dots ,\phi_u)\nabla_{i_1}\upsilon,$$

$$\tilde{C}^{t,i_1\dots i_\mu}_{g}(\Omega_1,\dots ,\Omega_p,
\phi_1,\dots ,\phi_u)\nabla_{i_1}\upsilon.$$

\par We immediately observe that for each $l\in L^z$:

\begin{equation}
\label{tsuna}\begin{split} &Replace[\tilde{C}^{l,i_1\dots i_\mu}_{g}(\Omega_1,\dots
,\Omega_p, \phi_1,\dots ,\phi_u)\nabla_{i_1}\phi_{u+1}]
\\&=C^{l,i_1,i_2\dots i_\mu}_{g}(\Omega_1,\dots ,\Omega_p,
\phi_1,\dots ,\phi_u)\nabla_{i_1}\upsilon.
\end{split}
\end{equation}

We then claim that the vector field

$$\Sum_{r\in R} a_r \tilde{C}^{r,i_1\dots i_{\mu+1}}_{g}
(\Omega_1,\dots ,\Omega_p,\phi_1,\dots ,
\phi_u)\nabla_{i_1}\upsilon\dots \nabla_{i_\mu}\upsilon$$
is the one we need for Proposition \ref{giade}. In order to see this, we only have to recall that
(\ref{ajax}) holds formally. We then ``memorize'' the sequence of
 formal applications of the identities in Definition 6 in \cite{a:dgciI}, by
which we can make the linearization of (\ref{ajax}) formally equal
to the linearization of the right hand side.
 We recall that an application of the identities in Definition 6 in \cite{a:dgciI} to
 the factor $\nabla^{(b)}_{r_1\dots r_b}Y$ (subject to the extra feature)
means that we may freely {\it permute} the
 indices ${}_{r_{\tau+q}},\dots ,{}_{r_b}$.

\par Now, we will perform the same sequence of applications of the
identities in Definition 6 in \cite{a:dgciI} to the linear combination:

\begin{equation}
\label{lafforgue}
\begin{split}
&\Sum_{l\in L^z} a_l \tilde{C}^{l,i_1\dots
i_\mu}_{g}(\Omega_1,\dots ,\Omega_p,\phi_1,\dots
,\phi_u)\nabla_{i_1}\upsilon\nabla_{i_2}\upsilon\dots
\nabla_{i_\mu}\upsilon
\\&- \Sum_{r\in R^z} a_r Xdiv_{i_{\mu +1}}\tilde{C}^{r,i_1\dots
i_{\mu+1}}_{g}(\Omega_1,\dots ,\Omega_p,\phi_1,\dots
,\phi_u)\nabla_{i_1}\upsilon
\nabla_{i_2}\upsilon\dots\nabla_{i_\mu}\upsilon
\end{split}
\end{equation}

 We impose one restriction: When we had permuted the
indices ${}_{r_1},\dots ,{}_{r_b}$ in a factor
$\nabla^{(b)}_{r_1\dots r_b}Y$ in (\ref{ajax}), we now freely permute them again,
but also introduce correction terms, by virtue of the Bianchi
identities:

\begin{equation}
\label{koichi1} \nabla^{(m)}_{r_1\dots
r_m}R_{ijkl}\nabla^i\phi_1\nabla^k\upsilon-
\nabla^{(m)}_{r_1\dots
r_m}R_{ilkj}\nabla^i\phi_1\nabla^k\upsilon=
\nabla^{(m)}_{r_1\dots
r_m}R_{ikjl}\nabla^i\phi_1\nabla^k\upsilon,
\end{equation}

\begin{equation}
\label{koichi2} \nabla^{(m)}_{r_1\dots
r_m}R_{ijkl}\nabla^i\phi_1\nabla^k\upsilon-
\nabla^{(m)}_{r_1\dots
j}R_{ir_mkl}\nabla^i\phi_1\nabla^k\upsilon=
\nabla^{(m)}_{r_1\dots
i}R_{r_mjkl}\nabla^i\phi_1\nabla^k\upsilon,
\end{equation}

\begin{equation}
\label{koichi3} \nabla^{(m)}_{r_1\dots
r_m}R_{ijkl}\nabla^i\phi_1\nabla^k\upsilon-
\nabla^{(m)}_{r_1\dots
l}R_{ijkr_m}\nabla^i\phi_1\nabla^k\upsilon=
\nabla^{(m)}_{r_1\dots
k}R_{r_mlij}\nabla^i\phi_1\nabla^k\upsilon.
\end{equation}

Hence, we derive that modulo complete contractions of length
$\ge\sigma +u+\mu +1$:

\begin{equation}
\label{sarnak}
\begin{split}
&\Sum_{l\in L^z} a_l \tilde{C}^{l,i_1\dots
i_\mu}_{g}(\Omega_1,\dots ,\Omega_p,\phi_1,\dots
,\phi_u)\nabla_{i_1}\upsilon\dots \nabla_{i_\mu}\upsilon
\\&- \Sum_{r\in R^z} a_r Xdiv_{i_{\mu +1}}\tilde{C}^{r,i_1\dots
i_{\mu+1}}_{g}(\Omega_1,\dots ,\Omega_p,\phi_1,\dots ,\phi_u)
\nabla_{i_1}\upsilon\dots \nabla_{i_\mu}\upsilon=
\\&\Sum_{t\in T} a_t \tilde{C}^{t,i_1\dots i_\mu}_{g}
(\Omega_1,\dots ,\Omega_p,\phi_1,\dots ,\phi_u)
\nabla_{i_1}\upsilon\dots\nabla_{i_\alpha}\upsilon
\\& +\Sum_{t\in T'} a_t \tilde{C}^{t,i_1\dots i_\mu}_{g}
(\Omega_1,\dots ,\Omega_p,\phi_1,\dots ,\phi_u)
\nabla_{i_1}\upsilon\dots\nabla_{i_\mu}\upsilon,
\end{split}
\end{equation}
where

\begin{equation}
\label{tutu} \Sum_{t\in T'} a_t \tilde{C}^{t,i_1\dots i_\mu}_{g}
(\Omega_1,\dots ,\Omega_p,\phi_1,\dots ,\phi_u)
\nabla_{i_1}\upsilon\dots\nabla_{i_\mu}\upsilon
\end{equation}
stands for the linear combination of correction terms that arises
 by virtue of the identities (\ref{koichi1}), (\ref{koichi2}), (\ref{koichi3}).

Specifically, the linear combination of correction terms arises by
replacing the crucial factor $ \nabla^{(m)}_{r_1\dots
r_m}R_{ijkl}\nabla^i\phi_1\nabla^k\phi_{u+1}$ by one of the
expressions on the right hand
 sides of (\ref{koichi1}), (\ref{koichi2}), (\ref{koichi3}).
We observe that all the correction terms are acceptable
tensor fields that are (simply or doubly) {\it subsequent} to
$\vec{L^z}'$.

\par Thus we have proven that under the assumptions of Lemma
\ref{zetajones}, the claim of Proposition \ref{giade} follows in
this  case I. $\Box$
\newline

{\bf Proof of Lemma \ref{technlem1}:}

 We refer to (\ref{ajaxnew}). 
We pick out the sublinear combination of terms in that equation 
where all indices in the function $\nabla^{(B)}Y$ are contracting
against a factor $\nabla\upsilon$ or $\nabla\phi'_h$
(denote the index set of tensor fields in $R^z$
with that property by $\overline{R}^z$). We thus
obtain a new equation:

\begin{equation}
\label{juliocesar} \begin{split}
 &\Sum_{r\in \overline{R}^z} a_r
X_{*}div_{i_{M+1}}\dots X_{*}div_{i_{\mu
+1}}C^{r,i_2\dots i_{\mu+1}}_{g}(\Omega_1,\dots
,\Omega_p,Y,\phi_2,\dots ,\phi_u)
\nabla_{i_2}\upsilon\dots\nabla_{i_M}\upsilon
\\&=\sum_{j\in J} a_j C^{j,i_2\dots i_M}_g(\Omega_1,\dots
,\Omega_p,Y,\phi_2,\dots
,\phi_u)\nabla_{i_2}\upsilon\dots\nabla_{i_M}\upsilon,
\end{split}
\end{equation}
(where the terms indexed in $J$ are simply subsequent to the
simple character of the terms in the first line). If $\sigma\ge 4$ and {\it if} 
 our Lemma assumption (\ref{assumpcion}) 
{\it does not} fall under the delicate case,\footnote{We will 
explain what to do if (\ref{juliocesar}) {\it does} fall under 
the delicate case below. The fact that (\ref{assumpcion}) does not 
fall under the delicate case ensures that (\ref{juliocesar})
 does not fall under a forbidden case of Lemma \ref{petermichel}.} 
we then apply
Lemma \ref{petermichel}; if $\sigma=3$  and {\it if} 
 our Lemma assumption (\ref{assumpcion}) 
{\it does not} fall under the delicate case we apply   
Lemma \ref{petermichel3}.\footnote{We observe
that Lemma \ref{petermichel3} can be applied, since we are in a
non-delicate case hence the tensor fields indexed in $R^z$ which
have all indices in $\nabla^{(B)}Y$ contracting against factors
$\nabla\upsilon$ must have a removable index.} We thus derive:

\begin{equation}
\label{juliocesar'}
\begin{split}
 &\Sum_{r\in \overline{R}^z} a_r
Xdiv_{i_{M+1}}\dots Xdiv_{i_{\mu +1}}C^{r,i_2\dots
i_{\mu+1}}_{g}(\Omega_1,\dots ,\Omega_p,Y,\phi_2,\dots ,\phi_u)
\nabla_{i_2}\upsilon\dots\nabla_{i_M}\upsilon
\\& =\Sum_{r\in \tilde{R}^z} a_r
Xdiv_{i_{M+1}}\dots Xdiv_{i_{\mu+1}}C^{r,i_2\dots
i_{\mu+1}}_{g}(\Omega_1,\dots ,\Omega_p,Y,\phi_2,\dots ,\phi_u)
\nabla_{i_2}\upsilon\dots\nabla_{i_M}\upsilon
\\&+\sum_{j\in J} a_j C^{j,i_2\dots i_M}_g(\Omega_1,\dots
,\Omega_p,Y,\phi_2,\dots
,\phi_u)\nabla_{i_2}\upsilon\dots\nabla_{i_M}\upsilon,
\end{split}
\end{equation}
where the terms indexed in $\tilde{R}^z$ have all the features of
the terms in $R^z$ but in addition have at least one index in
$\nabla^{(B)}Y$ {\it not} contracting against a factor
$\nabla\phi_h$ or $\nabla\upsilon$. Replacing the above into
(\ref{ajaxnew}), we may assume that all terms in $R^z$ have at
least one index in $\nabla^{(B)}Y$ {\it not} contracting against a
factor $\nabla\phi_h$ or $\nabla\upsilon$.
\newline

{\it Proof of (\ref{juliocesar'}) when 
(\ref{assumpcion}) falls under the delicate case:}
 All the arguments above can be repeated
{\it except for the application of Lemma \ref{petermichel}}, because 
in this setting (\ref{juliocesar}) might fall under the 
forbidden cases of Lemma \ref{petermichel}.  On
the other hand, we have also imposed a ``delicate assumption'' on
(\ref{assumpcion}), which we will utilize now:
\newline

{\it Proof that the technical claim follows from (\ref{ajax}) when 
(\ref{assumpcion}) falls under the
 ``delicate case'':}

Refer to (\ref{ajax}). We denote by $R^z_{Bad}\subset R^z$ the
index set of tensor fields that have the free index
${}_{i_{\mu+1}}$ being a special index in some simple factor
$S_{*}\nabla^{(\nu)}R_{ijkl}$, {\it and} with all indices in
$\nabla^{(B)}Y$ contracting against a factor $\nabla\phi_h$ or
$\nabla\upsilon$ (we will call tensor fields with that property
``bad''). We observe that (\ref{juliocesar'}) can again be
derived, {\it provided} that $R^z_{Bad}=\emptyset$ in
(\ref{juliocesar}); (this is because when $R^z_{Bad}=\emptyset$ 
there is no danger of falling under a ``forbidden case'' of 
Lemma \ref{petermichel}, by weight considerations). 
We will now show that we can write:

\begin{equation}
\label{lastrada}
\begin{split}
&\Sum_{r\in R^z_{Bad}} a_r C^{r,i_2\dots
i_{\mu+1}}_{g}(\Omega_1,\dots ,\Omega_p,Y,\phi_2,\dots ,\phi_u)
\nabla_{i_2}\upsilon\dots\nabla_{i_\mu}\upsilon\nabla_{i_{\mu+1}}\Phi=
\\&\Sum_{r\in R^z_{NotBad}} a_r
C^{r,i_2\dots i_{\mu+1}}_{g}(\Omega_1,\dots ,\Omega_p,Y,\phi_2,\dots ,\phi_u)
\nabla_{i_2}\upsilon\dots\nabla_{i_\mu}\upsilon\nabla_{i_{\mu+1}}\Phi+
\\&\sum_{j\in J} a_j C^{j,i_2\dots i_{\mu+1}}_{g}(\Omega_1,\dots ,\Omega_p,Y,\phi_2,\dots ,\phi_u)
\nabla_{i_2}\upsilon\dots\nabla_{i_\mu}\upsilon\nabla_{i_{\mu+1}}\Phi,
\end{split}
\end{equation}
where the terms in $R^z_{NotBad}$ are in the form described in
(\ref{ajax}) and in addition are not bad. The terms indexed in $J$
are simply subsequent to $\vec{\kappa}_{simp}$.
 If we can show the above, then by making the
$\nabla\Phi$ into an $Xdiv$ and replacing into (\ref{ajax}),
we are reduced to the case $R^z_{Max}=\emptyset$. We have then
noted that the proof above goes through. 
Thus matters are reduced to showing (\ref{lastrada}).

{\it Proof of (\ref{lastrada}):} We may assume with no loss of generality that the free index
${}_{i_{\mu+1}}$ is the index ${}_k$ in a factor
$S^\sharp\nabla^{(\rho)}_{r_1\dots r_\rho}R_{ijkl}$, where
$S^\sharp$ stands for the symmetrization over the index ${}_l$ and
all the indices in the above that are not contracting against a
factor $\nabla\upsilon$ or $\nabla\phi'_h$ (the correction terms
that we obtain from this $S^\sharp$-symmetrization would be tensor
fields which are not ``bad''--as allowed in (\ref{lastrada})).
 We then pick out the sublinear combination of terms in (\ref{ajax}) with a factor $\nabla^{(B)}Y$
that has all its indices {\it except} one (say the index ${}_s$)
contracting against a factor $\nabla\phi_h$ or $\nabla\upsilon$,
and the index ${}_s$ contracting against a special index in a
factor $S_{*}\nabla^{(\rho)}R_{ijkl}$. By virtue of the ``delicate
assumption'', this sublinear combination will be of the form:

\begin{equation}
\label{nektarios} \begin{split} &\sum_{r\in R^z_{BAD}} a_r
Hit_Ydiv_{i_{\mu+1}} C^{r,i_2\dots i_{\mu+1}}_{g}(\Omega_1,\dots
,\Omega_p,Y,\phi_2,\dots ,\phi_u)
\nabla_{i_2}\upsilon\dots\nabla_{i_\mu}\upsilon+
\\&\sum_{j\in J} a_j
C^{j,i_2\dots i_\mu}_g(\Omega_1,\dots ,\Omega_p,Y,\phi_2,\dots
,\phi_u) \nabla_{i_2}\upsilon\dots\nabla_{i_\mu}\upsilon.
\end{split}
\end{equation}
Then, we consider the first conformal variation $Image^1_X[F_g]$
of (\ref{ajax}) and we pick out the terms where one of the factors
$\nabla\phi_h, h\in Def(\vec{\kappa}_{simp})$ is contracting
against the factor $\nabla^{(B)}Y$, {\it for which we additionally
require that all other indices contract against $\nabla\upsilon$'s
or $\nabla\phi$'s}. The above sublinear combination must vanish
separately. We thus derive a new equation:

\begin{equation}
\label{nektarios} \begin{split} &\sum_{r\in R^z_{BAD}} a_r
Hit_Ydiv_{i_{\mu+1}} Op[C]^{r,i_2\dots
i_{\mu+1}}_{g}(\Omega_1,\dots ,\Omega_p,Y,\phi_2,\dots ,\phi_u)
\nabla_{i_2}\upsilon\dots\nabla_{i_\mu}\upsilon+
\\&\sum_{j\in J} a_j
C^{j,i_2\dots i_\mu}_g(\Omega_1,\dots ,\Omega_p,Y,\phi_2,\dots
,\phi_u) \nabla_{i_2}\upsilon\dots\nabla_{i_\mu}\upsilon=0,
\end{split}
\end{equation}
where $Op[C]^{r,i_2\dots i_{\mu+1}}_{g}$ formally arises from by
replacing the factor \\$S^\sharp \nabla^{(\rho)}_{r_1\dots
r_\rho}R_{iji_{\mu+1}l}\nabla^i\tilde{\phi}_q$ by
$\nabla^{(\rho+2)}_{r_1\dots r_\rho jl}Y\nabla_{i_{\mu+1}}\phi_q$
(and $Hit_Ydiv_{i_{\mu+1}}$ still means that $\nabla^{i_{\mu+1}}$
is forced to hit the factor $\nabla^{(B)}Y$). Then, formally
replacing the expression
$$\nabla^{(C)}_{y_1\dots
y_C}X\nabla^{y_1}\upsilon\dots\nabla^{y_a}\upsilon\nabla^{y_{a+1}}\phi_{w_1}\dots\nabla^{y_b}\phi_{w_f}$$
by an expression $$S^\sharp\nabla^{(C-2)}_{y_1\dots
y_{C-2}}R_{iy_{C-1}i_{\mu+1}y_C}\nabla^{y_1}
\upsilon\dots\nabla^{y_a}\upsilon\nabla^{y_{a+1}}
\phi_{w_1}\dots\nabla^{y_b}\phi_{w_f}\nabla^i\tilde{\phi}_q,$$ 
and repeating the formal identities by which (\ref{nektarios}) 
is proven ``formally'',\footnote{See the argument that proves (\ref{sarnak}).} we
derive (\ref{lastrada}). $\Box$
\newline

\subsection{Derivation of Proposition \ref{giade} in case II from Lemma
\ref{pool2}.}
\label{stellara1}

\par Recall the notation of Lemma \ref{pool2}. Our point of departure will be the 
Lemma's conclusion, equation (\ref{vecFskillb}).
Recall that in that equation, all {\it tensor fields} have a given $(u+1)$-simple character,
which we have denoted by $\vec{\kappa}'_{simp}$. We also recall that the
 $(\mu-1)$-tensor fields $\tilde{C}^{l,i_1\dots i_\mu}_{g}(\Omega_1,\dots ,
\Omega_p,\phi_1,\dots ,\phi_u)\nabla_{i_h}\phi_{u+1}$
in the first line of (\ref{vecFskillb})
have {\it maximal refined double characters among 
all $(\mu-1)$-refined double characters
appearing in (\ref{vecFskillb})};
 we have denoted
the maximal $(\mu-1)$-refined double characters by 
$\vec{\kappa}^z_{ref-doub}, z\in Z'_{Max}$.

\par For all tensor fields in (\ref{vecFskillb}) the factor
$\nabla\phi_{u+1}$ is contracting against an index ${}_i$ in some
chosen factor $S_{*}\nabla^{(\nu)}R_{ijkl}$.
For this subsection, we will be calling that factor the A-crucial factor.
 We recall that all tensor fields in the
first line in (\ref{vecFskillb}) have a given number of special free indices in the
 A-crucial factor. This follows from the definition of $Z'_{Max}$.

\par We further distinguish two subcases of case II: Observe that either all
refined $(u+1,\mu-1)$-double characters $\vec{\kappa}^z, z\in Z'_{Max}$  one
 internal free index ${}_k$ or ${}_l$ in the A-crucial factor, or
 have no such free index.\footnote{Notice
that this dichotomy corresponds to the following dichotomy
regarding the tensor fields in $\bigcup_{z\in Z'_{Max}}L^z$ in
(\ref{assumpcion}): Either for those tensor fields the critical
factor $\nabla^{(m)}R_{ijkl}$ contains two internal free indices,
or it contains one internal free index.} We accordingly call
these subcases A and B,\footnote{Sometimes we will
 refer to subcases IIA, IIB, to stress that these are
subcases of case II.}
 and we will prove our claim separately for these two subcases.
\newline

\par We here prove our assertion for all cases apart from
certain {\it special cases}  where we will derive the claim of
Proposition \ref{giade} directly from (\ref{assumpcion}) 
(in the paper \cite{alexakis5} in this series). 
\newline

{\bf The special cases:}
\newline

{\it Case A:} The special case here is when for each
tensor field $C^{l,i_1\dots i_\mu}_g, l\in L^z,z\in Z'_{Max}$ in
(\ref{assumpcion}) there are no removable free indices
 among any of its factors.\footnote{Observe that if this 
is true of one of the tensor fields $C^{l,i_1\dots i_\mu}_g, l\in L^z,z\in Z'_{Max}$, 
it will be true of all of them, by weight considerations.}

{\it Case B:} The special case here is when for each
tensor field $C^{l,i_1\dots i_\mu}_g, l\in L^z,z\in Z'_{Max}$ in
(\ref{assumpcion}) there are no removable indices {\it other than} 
(possibly) the indices ${}_k,{}_l$ in the (one of the) crucial factor
  $\nabla^{(m)}_{v_1\dots v_xi_1\dots i_b}R_{i_{b+1}jkl}$.\footnote{The 
remark in the previous footnote still holds.}
\newline

\par Throughout the rest of this subsection we will be
 assuming that the $\mu$-tensor fields of maximal refined double
 character in (\ref{assumpcion}) are {\it not} special. In the special cases, Proposition \ref{giade}
will be derived directly (without recourse to Lemma \ref{pool2}) in \cite{alexakis5}.  
\newline

{\it Derivation of Proposition \ref{giade} in case II from Lemma
\ref{pool2} (subcase A):}  Recall the conclusion of Lemma
\ref{pool2}:

\begin{equation}
\label{rewrite1}
\begin{split}
&\Sum_{z\in Z'_{Max}}\Sum_{l\in L^z} a_l \Sum_{i_h\in
I_{*,l}}Xdiv_{i_1}\dots\hat{Xdiv}_{i_h}\dots
 Xdiv_{i_\mu}\tilde{C}^{l,i_1\dots
i_\mu}_{g}(\Omega_1,\dots ,\Omega_p,\phi_1,\dots
,\phi_u)\\&\nabla_{i_h}\phi_{u+1}
+\Sum_{\nu\in N} a_\nu Xdiv_{i_2}\dots Xdiv_{i_\mu}
C^{\nu,i_1\dots i_{\mu}}_{g}(\Omega_1,\dots
,\Omega_p,\phi_1,\dots ,\phi_u)\nabla_{i_1}\phi_{u+1} +
\\& \Sum_{d\in D} a_d Xdiv_{i_1}\dots Xdiv_{i_\mu} C^{d,i_1\dots
i_\mu}_{g}(\Omega_1,\dots ,\Omega_p,\phi_1,\dots
,\phi_u,\phi_{u+1})=
\\& \Sum_{t\in T} a_t C^{t,i_{*}}_{g}
(\Omega_1,\dots ,\Omega_p,\phi_1,\dots
,\phi_u,\phi_u)\nabla_{i_{*}}\phi_{u+1}.
\end{split}
\end{equation}

\par We now apply our inductive assumption of Corollary \ref{corgiade}
to the above.\footnote{The fact that we are assuming 
that (\ref{hypothese2}) (the assumption of Lemma \ref{pool2}) 
does not fall under a ``special case'' of subcase A ensures that 
(\ref{rewrite1}) satisifies the requirements  of Corollary \ref{corgiade}.
Also, since we have introduced a new factor $\nabla\phi_{u+1}$, 
the above falls under the inductive assumption of Corollary \ref{corgiade}.} 
We derive that there is a linear combination of acceptable $(\mu-1)$-tensor
fields with a $(u+1)$-simple character $\vec{\kappa}'_{simp}$ (indexed in $H$ below), so that:

\begin{equation}
\label{mayrh1/2}
\begin{split}
&\Sum_{z\in Z'_{Max}}\Sum_{l\in L^z} a_l \Sum_{i_h\in I_{*,l}}
\tilde{C}^{l,i_1\dots i_\mu}_{g} (\Omega_1,\dots
,\Omega_p,\phi_1,\dots
,\phi_u)\nabla_{i_h}\phi_{u+1}\nabla_{i_1}\upsilon\dots\hat{\nabla}_{i_h}\upsilon
\dots\nabla_{i_\mu}\upsilon
\\&+\Sum_{\nu\in N} a_\nu C^{\nu,i_1\dots i_\mu}_{g}
(\Omega_1,\dots ,\Omega_p,\phi_1,\dots
,\phi_u)\nabla_{i_1}\phi_{u+1}\nabla_{i_2}\upsilon
\dots\nabla_{i_\mu}\upsilon+
\\&\Sum_{h\in H} a_h X div_{i_{\mu +1}}
C^{h,i_1\dots i_{\mu +1}}_{g}(\Omega_1,\dots ,\Omega_p,
\phi_1,\dots,\phi_u)\nabla_{i_1}\phi_{u+1}
\nabla_{i_2}\upsilon\dots\nabla_{i_\mu}\upsilon+
\\&\Sum_{j\in J} a_j C^{j,i_1\dots i_\mu}_{g}
(\Omega_1,\dots ,\Omega_p,\phi_1,\dots
,\phi_u)\nabla_{i_1}\phi_{u+1}\nabla_{i_2}\upsilon
\dots\nabla_{i_\mu}\upsilon,
\end{split}
\end{equation}
where each $C^{j,i_1\dots i_\mu}_{g} (\Omega_1,\dots
,\Omega_p,\phi_1,\dots ,\phi_u)\nabla_{i_1}\phi_{u+1}$ is simply
subsequent to \\$Simp(\vec{\kappa}^z_{ref-doub})$.
{\it Note:} We will now be calling the factor against which
$\nabla\phi_{u+1}$ is contracting the A-crucial factor for {\it
any} contraction appearing in the above.

\par  Now, we observe that in view of the claim in Lemma \ref{pool2}, for each
$\tilde{C}^{l,i_1\dots i_\mu}_{g}$, $C^{\nu,i_1\dots i_\mu}_{g}$
we must have that at least one
 of the indices ${}_{r_1},\dots ,{}_{r_\nu},{}_j$ in the A-crucial
 factor $S_{*}\nabla^{(\nu)}_{r_1\dots r_\nu}R_{ijkl}$ is not free and
 not contracting against a factor $\nabla\phi_f,f\le u$.
Furthermore, it follows from the definition of $\bigcup_{z\in Z'_{Max}} L^z$ that all $\tilde{C}^{l,i_1\dots
i_\mu}_{g}$ tensor fields are contracting against a given number
$c$ of factors $\nabla\upsilon$. We denote by $\overline{N}\subset
N, \overline{H}\subset H, \overline{J}\subset J$ the index sets of
the contractions above for which the $A$-crucial factors is
contracting against $c$ factors $\nabla\upsilon$. Then, since
(\ref{mayrh1/2}) holds formally, we derive that:

\begin{equation}
\label{mayrh3/4}
\begin{split}
&\Sum_{z\in Z'_{Max}}\Sum_{l\in L^z} a_l \Sum_{i_r\in I_{*,l}}
\tilde{C}^{l,i_1\dots i_\mu}_{g} (\Omega_1,\dots
,\Omega_p,\phi_1,\dots
,\phi_u)\nabla_{i_r}\phi_{u+1}\nabla_{i_1}\upsilon\dots\hat{\nabla}_{i_r}\upsilon
\dots\nabla_{i_\mu}\upsilon
\\&+\Sum_{\nu\in \overline{N}} a_\nu C^{\nu,i_1\dots i_\mu}_{g}
(\Omega_1,\dots ,\Omega_p,\phi_1,\dots
,\phi_u)\nabla_{i_1}\phi_{u+1}\nabla_{i_2}\upsilon
\dots\nabla_{i_\mu}\upsilon=
\\&\Sum_{h\in \overline{H}} a_h X div_{i_{\mu +1}}
C^{h,i_1\dots i_{\mu +1}}_{g}(\Omega_1,\dots ,\Omega_p,
\phi_1,\dots,\phi_u)\nabla_{i_1}\phi_{u+1}
\nabla_{i_2}\upsilon\dots\nabla_{i_\mu}\upsilon+
\\&\Sum_{j\in \overline{J}} a_j C^{j,i_1\dots i_\mu}_{g}
(\Omega_1,\dots ,\Omega_p,\phi_1,\dots
,\phi_u)\nabla_{i_1}\phi_{u+1}\nabla_{i_2}\upsilon
\dots\nabla_{i_\mu}\upsilon.
\end{split}
\end{equation}

\par In fact, we will be using a weaker version of
 this equation:
Consider (\ref{mayrh3/4}): We assume with no loss of generality, only for
 notational convenience, that in each
$\tilde{C}^{l,i_1\dots i_{\mu}}$, $C^{\nu,i_1\dots i_\mu}_{g}$
 and each $C^{h,i_1\dots i_{\mu +1}}$
the free indices ${}_{i_1},\dots ,{}_{\hat{i_r}},\dots
,{}_{i_{c+1}}$ belong to the A-crucial factor, while the indices
${}_{i_{c+2}}, \dots ,{}_{i_\mu}$ do not. Now, recall that
(\ref{mayrh3/4}) holds formally. Then, define an operation that
formally erases the factors $\nabla\upsilon$ that are {\it not}
contracting against the A-crucial factor and then takes $Xdiv$s of
the resulting free indices that we obtain.  If we denote
  the expression that we (formally) thus obtain by $F'$,
it follows that $F'=0$ (modulo longer complete 
contractions) by virtue of the last Lemma in
 the Appendix of \cite{alexakis1}. We have thus derived:\footnote{We
will revert to writing $N,H,J$ instead of $\overline{N}$,
$\overline{H}$, $\overline{J}$ for notational simplicity.}

\begin{equation}
\label{weakmayrh}
\begin{split}
&\Sum_{z\in Z'_{Max}} \Sum_{l\in L^z} a_l \Sum_{i_h\in
I_{*,l}}Xdiv_{i_{c+2}}\dots Xdiv_{i_\mu}
\\&\tilde{C}^{l,i_1\dots i_\mu}_{g} (\Omega_1,\dots,\Omega_p,
\phi_1,\dots,\phi_u)\nabla_{i_{1}}\upsilon\dots\nabla_{i_h}\phi_{u+1}
\dots\nabla_{i_{c+1}}\upsilon+
\\&\Sum_{\nu\in N} a_\nu Xdiv_{i_{c+2}}\dots Xdiv_{i_\mu} C^{\nu,i_1\dots i_\mu}_{g}
(\Omega_1,\dots,\Omega_p,\phi_1,\dots,\phi_u)\nabla_{i_1}
\phi_{u+1}\nabla_{i_{2}}\upsilon\dots\nabla_{i_{c+1}}\upsilon
\\&=\Sum_{h\in H} a_h Xdiv_{i_{c+1}}\dots 
Xdiv_{i_{\mu +1}} C^{h,i_1\dots i_{\mu +1}}_{g}(\Omega_1,\dots
,\Omega_p, \phi_1,\dots,\phi_u)\nabla_{i_1}\phi_{u+1}\\&
\cdot\nabla_{i_2}\upsilon\dots\nabla_{i_{c+1}}\upsilon
+\Sum_{j\in J'} a_j
C^{j,i_1\dots i_{c+1}}_{g}(\Omega_1,\dots ,\Omega_p,
\phi_1,\dots,\phi_u)\nabla_{i_1}\phi_{u+1}
\nabla_{i_2}\upsilon\dots\nabla_{i_{c+1}}\upsilon.
\end{split}
\end{equation}
Here the complete contractions indexed in $J'$ are simply
subsequent to $\vec{\kappa}'_{simp}$.

Now, we break the index set $H$ into two subsets: $h\in H_1$
 if and only if the A-crucial factor
$S_{*}\nabla^{(\nu)}_{r_1\dots r_\nu}R_{ijkl}$ in
$$C^{h,i_1\dots i_{\mu +1}}_{g}(\Omega_1,\dots ,\Omega_p,
\phi_1,\dots,\phi_u)\nabla_{i_1}\phi_{u+1}
\nabla_{i_2}\upsilon\dots\nabla_{i_\pi}\upsilon$$ has the property
that at least one of the indices ${}_{r_1},\dots ,{}_{r_\nu},{}_j$
is neither contracting against a factor $\nabla\upsilon$ nor a
factor $\nabla\phi'_h$. We say that $h\in H_2$ if and only
 if all the indices ${}_{r_1},\dots ,{}_{r_\nu},{}_j$ in the A-crucial factor
  are contracting against a
 factor $\nabla\phi'_h$ or $\nabla\upsilon$. We 
complete our proof in two steps: Step 1 involves getting rid of the terms indexed in $H_2$:
\newline

{\it Step 1:} We introduce some notation:

\begin{definition} 
\label{Xstar} We denote by $X_{*} div_{i}[\dots]$ the sublinear combination in
each $X div_{i}[\dots]$ where we impose the additional restriction that
$\nabla_i$ is not allowed to hit the A-crucial factor (in the
 form $S_{*}\nabla^{(\nu)}R_{ijkl}$). 
\end{definition} 
 Then, since
(\ref{weakmayrh}) holds formally, we deduce that:

\begin{equation}
\label{guillemin}
\begin{split}
& \Sum_{h\in H_2} a_h X_{*} div_{i_{c+2}}\dots X_{*} div_{i_{\mu
+1}} C^{h,i_1\dots i_{\mu+1}}_{g}(\Omega_1,\dots ,\Omega_p,
\phi_1,\dots,\phi_u)\nabla_{i_1}\phi_{u+1}
\nabla_{i_2}\upsilon\\&\dots\nabla_{i_{c+1}}\upsilon
 +\Sum_{j\in J} a_j
C^{j,i_2\dots i_{c+2}}_{g}(\Omega_1,\dots ,\Omega_p,
\phi_1,\dots,\phi_u)\nabla_{i_1}\phi_{u+1}
\nabla_{i_2}\upsilon\dots\nabla_{i_{c+2}}\upsilon=0,
\end{split}
\end{equation}
(modlo longer terms), where each $C^{j,i_2\dots i_\pi}_{g}$ is simply subsequent to
$\vec{\kappa}'_{simp}$.

\par A notational convention that can be made with no loss of
 generality is that in $\vec{\kappa}^z_{simp}$ the $b$
  factors $\nabla\phi_h'$ that are contracting against
 indices ${}_{r_1},\dots ,{}_{r_\nu},{}_j$ in the A-crucial factor
 $S_{*}\nabla^{(\nu)}R_{ijkl}$ are precisely the factors
$\nabla\phi_1,\dots ,\nabla\phi_b$.

\par We then claim that we can write:

\begin{equation}
\label{guilleminberk}
\begin{split}
& \Sum_{h\in H_2} a_h X div_{i_{c+2}}\dots X div_{i_{\mu +1}}
C^{h,i_1\dots i_{\mu+1}}_{g}(\Omega_1,\dots ,\Omega_p,
\phi_1,\dots,\phi_u)\nabla_{i_1}\phi_{u+1}
\nabla_{i_2}\upsilon\dots\\&\nabla_{i_{c+1}}\upsilon
 =\Sum_{h\in \tilde{H}} a_h X div_{i_{c+2}}\dots X
div_{i_{\mu +1}} C^{h,i_1\dots i_{\mu+1}}_{g}(\Omega_1,\dots
,\Omega_p, \phi_1,\dots,\phi_u)\nabla_{i_1}\phi_{u+1}
\nabla_{i_2}\upsilon\\&\dots\nabla_{i_{c+1}}\upsilon
+\Sum_{j\in J} a_j
C^{j,i_2\dots i_{c+2}}_{g}(\Omega_1,\dots ,\Omega_p,
\phi_1,\dots,\phi_u)\nabla_{i_1}\phi_{u+1}
\nabla_{i_2}\upsilon\dots\nabla_{i_{c+2}}\upsilon,
\end{split}
\end{equation}
where the sublinear combination indexed in $\tilde{H}$
in the RHS stands for a generic linear combination
 of tensor fields with all the properties of the tensor fields indexed in $H_1$ above.

\par We will show below that  (\ref{guilleminberk}) follows by
applying Lemma \ref{obote} or Lemma \ref{obote3} to
(\ref{guillemin}) (and we will prove Lemmas
\ref{obote}, \ref{obote3} in \cite{alexakis5}). For now, let us observe
how  (\ref{guilleminberk}) implies Proposition \ref{giade}
in this case A:
\newline

{\it Step 2: Proposition \ref{giade} follows from (\ref{guilleminberk}).}
By replacing (\ref{guilleminberk}) into (\ref{weakmayrh}), we are
reduced to showing
 our claim when $H_2=\emptyset$. Now, for each of the tensor fields in (\ref{weakmayrh}) we  denote by
$C^{l,i_{c+2}\dots i_\mu}_{g} (\Omega_1,\dots ,\Omega_p,
Y,\phi_{b+1},\dots,\phi_u)$, \\$C^{\nu,i_{c+2}\dots i_\mu}_{g}
(\Omega_1, \dots ,\Omega_p, Y,\phi_1,\dots,\phi_u)$,
$C^{h,i_{c+2}\dots i_{\mu +1}}_{g} (\Omega_1,\dots ,\Omega_p,
Y,\phi_{b+1},\dots ,\phi_u)$ the tensor fields that
 arise from
\\ $C^{l,i_1\dots i_\mu}_{g} (\Omega_1,\dots ,\Omega_p,
\phi_1,\dots ,\phi_u)\nabla_{i_h}\phi_{u+1}
\nabla_{i_2}\upsilon\dots\nabla_{i_{c+1}}\upsilon$,
$$C^{h,i_1\dots i_{\mu+1}}_{g} (\Omega_1,\dots ,\Omega_p,
 \phi_1,\dots,\phi_u)\nabla_{i_1}\phi_{u+1}
\nabla_{i_2}\upsilon\dots\nabla_{i_{c+1}}\upsilon$$
 by replacing the expressions
$$S_{*}\nabla^{(A+c+b-1)}_{u_1\dots u_{A-1}s_1\dots s_br_1\dots
 r_{c}}R_{iu_Akl}\nabla^{r_1}\upsilon\dots
 \nabla^{r_c}\upsilon\nabla^k\upsilon\nabla^{s_1}
\phi_1\dots \nabla^{s_b}\phi_b\nabla^i\phi_{u+1}$$ by a factor
$\nabla^{(A+1)}_{u_1\dots u_Al}Y$.

\par Now, by polarizing the function $\upsilon$ in (\ref{weakmayrh}) and applying Lemma \ref{technical1},
  we deduce that:

\begin{equation}
\label{weakmayrhb}
\begin{split}
& \Sum_{z\in Z'_{Max}} \Sum_{l\in L^z} a_l \Sum_{i_h\in I_{*,l}}
Xdiv_{i_{c+2}}\dots Xdiv_{i_\mu} \tilde{C}^{l,i_{c+2}\dots
i_\mu}_{g} (\Omega_1,\dots, \Omega_p,Y,\phi_{b+1},\dots,\phi_u)+
\\&\Sum_{\nu\in N} a_\nu
Xdiv_{i_{c+2}}\dots Xdiv_{i_\mu} C^{\nu,i_{c+2}\dots i_\mu}_{g}
(\Omega_1,\dots,\Omega_p,Y,\phi_{b+1},\dots,\phi_u)=
\\&\Sum_{h\in H} a_h Xdiv_{i_{c+1}}\dots Xdiv_{i_{\mu +1}}
C^{h,i_{c+2}\dots i_{\mu +1}}_{g}(\Omega_1,\dots ,\Omega_p,
Y,\phi_1,\dots,\phi_u)+
\\&\Sum_{j\in J} a_j
C^{j}_{g}(\Omega_1,\dots ,\Omega_p,Y, \phi_{b+1},\dots,\phi_u).
\end{split}
\end{equation}

{\it Notice that all the tensor fields in the above have at least
2 derivatives on the factor $\nabla^{(b)}Y$}. For the tensor
fields indexed in $H$, this follows since $H_2=\emptyset$; for the
tensor fields indexed in $(\bigcup_{z\in Z'_{Max}}L^z)\bigcup N$, it
follows by our observation after (\ref{mayrh1/2}). Thus, if we
treat
 $Y$ as a function $\Omega_{p+1}$, all tensor fields in the above have a given simple character,
 which we will denote by $\tilde{\kappa}_{simp}$.

  We denote the refined double character of the tensor fields
\\$C^{l,i_{c+2}\dots i_{\mu +1}}_{g} (\Omega_1,\dots
,\Omega_p,Y,\phi_{b+1},\dots,\phi_u)$, $l\in L^z,z\in Z'_{Max}$ by
$\vec{L}'_z$ (observe
 that they are the maximal refined double characters among all the
 $(\mu-c-1)$-tensor fields appearing in (\ref{weakmayrhb})).
By virtue of our inductive assumption of Proposition
\ref{giade},\footnote{Notice that  {\it since we
are assuming that the $\mu$-tensor fields of maximal refined
double character in (\ref{assumpcion}) do not have special free
indices in any factor $S_{*}\nabla^{(\nu)}R_{ijkl}$} then it
follows that the tensor fields of minimum rank in
(\ref{weakmayrhb}) satisfy the requirememnts of 
Proposition \ref{giade}.} we
 derive that for each $z\in Z'_{Max}$ there is a linear combination of acceptable
$(\mu-c)$-tensor fields \\ $\Sum_{h\in H'} a_h C^{h,i_{c+1}\dots i_{\mu +1}}_{g} (\Omega_1,\dots
,\Omega_p,Y,\phi_{b+1}, \dots,\phi_u)$ with an $(\mu-c-1)$-refined
double character  $\vec{L}'_z$ (so in particular they have a
factor $\nabla^{(B)}Y$, $B\ge 2$, not contracting against any
factor $\nabla\upsilon$, $\nabla\phi_h$), so that for each $z\in
Z'_{Max}$:

\begin{equation}
\label{mayrh2}
\begin{split}
&\Sum_{l\in L^z} a_l \Sum_{i_h\in I_{*,l}} C^{l,i_1i_{c+2}\dots
i_\mu}_{g} (\Omega_1,\dots ,\Omega_p,
Y,\phi_{b+1},\dots,\phi_u)\nabla_{i_{c+2}}\upsilon\dots
 \nabla_{i_\mu}\upsilon
\\&-\Sum_{h\in H} a_h Xdiv_{i_{\mu +1}}
C^{h,i_{c+2}\dots i_{\mu +1}}_{g}(\Omega_1,\dots
,\Omega_p,Y,\phi_{b+1},\dots
,\phi_u)\nabla_{i_{c+2}}\upsilon\dots\nabla_{i_\mu}\upsilon
\\& =\Sum_{t\in T} a_t C^{t,i_{c+2}\dots i_\mu}_{g}(\Omega_1,\dots ,\Omega_p,Y,
\phi_{b+1},\dots,\phi_u)\nabla_{i_{c+2}}\upsilon\dots
\nabla_{i_\mu}\upsilon,
\end{split}
\end{equation}
where each $C^{t,i_{c+2}\dots i_\mu}_{g}(\Omega_1, \dots
,\Omega_p,Y,\phi_{b+1},\dots,\phi_u)$  is acceptable and simply or
doubly subsequent to $\vec{L}'_z$.

\par We use the fact that (\ref{mayrh2}) holds formally. 
 We then define an operation $Op[\dots ]$
that replaces each factor $\nabla^{(B)}_{t_1\dots t_B}Y$ ($B\ge
2$) by an expression
$$\nabla^{(m)}_{s_1\dots s_br_1\dots r_{c-1}t_1\dots
t_{k-2}}R_{it_{B-1}kt_B}\nabla^i\upsilon\nabla^{r_1}
\upsilon\dots \nabla^{r_{c-1}}\upsilon \nabla^k
\upsilon\nabla^{s_1}\phi'_1\dots \nabla^{s_b}\phi'_b.$$
 We observe that for each $z\in Z'_{Max}$:

\begin{equation}
\label{church-hill1}
\begin{split}
&Op\{\Sum_{l\in L^z} a_l \Sum_{i_h\in I_{*,l}} \tilde{C}^{l,i_{c+2}\dots
i_\mu}_{g} (\Omega_1,\dots ,\Omega_p,
Y,\phi_{b+1},\dots,\phi_u)\nabla_{i_{c+2}}\upsilon\dots\nabla_{i_{\mu}}\upsilon\}=
\\& |I_{*,l}| \Sum_{l\in L^z} a_l C^{l,i_1\dots i_\mu}_{g}
(\Omega_1,\dots ,\Omega_p,\phi_1,\dots,\phi_u)
\nabla_{i_1}\upsilon\dots\nabla_{i_\mu}\upsilon,
\end{split}
\end{equation}
where we have noted that $|I_{*,l}|$ is universal, i.e.
independent of the element $l\in L^z$ (in most cases
$|I_{*,l}|=1$).

\par Hence, since (\ref{mayrh2}) holds formally, we deduce that:

\begin{equation}
\begin{split}
\label{vaskonik2} &\Sum_{l\in L^z} a_l C^{l,i_1\dots
i_\mu}_{g}(\Omega_1,\dots ,\Omega_p,\phi_1,\dots,\phi_u)
\nabla_{i_1}\upsilon\dots\nabla_{i_\mu}\upsilon-
\\&\Sum_{h\in H} a_h Xdiv_{i_{\mu +1}}Op[
C^{h,i_{\pi+1}\dots i_{\mu +1}}_{g}(\Omega_1,\dots
,\Omega_p,Y,\phi_1,\dots
,\phi_u)\nabla_{i_{\pi+1}}\upsilon\dots\nabla_{i_\mu}\upsilon]
\\&= \Sum_{t\in T'} a_t C^{t,i_1\dots i_\mu}_{g}(\Omega_1,\dots ,\Omega_p,
\phi_1,\dots,\phi_u)\nabla_{i_1}\phi_{u+1}
\nabla_{i_2}\upsilon\dots \nabla_{i_\mu}\upsilon,
\end{split}
\end{equation}
where again each $C^{t,i_1\dots i_\mu}_{g}(\Omega_1,\dots ,
\Omega_p,\phi_1,\dots,\phi_u)$ is either simply or doubly
 subsequent to $\vec{L^z}$. The above is obtained from
 (\ref{mayrh2}) by the usual argument as in the derivation of
 Proposition \ref{giade} from Lemma \ref{zetajones} (see the
 argument above equations (\ref{koichi1}), (\ref{koichi2}),
 (\ref{koichi3})): We
 may repeat the permutations by which we make (\ref{mayrh2})
 formally zero, modulo introducing corrections terms that are
 simply or doubly subsequent by virtue of the Bianchi identities.

\par Therefore, we have shown that Lemma \ref{pool2}
 implies case II of Proposition \ref{giade}
in subcase A (in the non-special cases), provided we 
can prove (\ref{guilleminberk}). We reduce (\ref{vaskonik2}) 
to certain other Lemmas, which will be proven in \cite{alexakis5} 
in the next  subsection. $\Box$
\newline

{\bf Derivation of case II of Proposition \ref{giade} from Lemma
\ref{pool2} in case B:}
\newline

\par Our point of departure is again equation (\ref{vecFskillb}).

\par Recall that in this second case, 
for each $z\in Z'_{Max}$  {\it none} of the free indices in the
A-crucial factor $S_{*}\nabla^{(\nu)}R_{ijkl}$ in any
$\tilde{C}^{l,i_1\dots i_\mu}_{g}\nabla_{i_1} \phi_{u+1}$, $l\in
L^z$ are special.

\par In that case, we again have equation (\ref{mayrh3/4}).
We will re-write the equation in a somewhat more convenient form,
but first we recall some of the notational conventions from the
previous case. For
notational convenience, we have assumed that the $b$ factors
$\nabla\phi_h$, $h\le u$ that are contracting against the
A-crucial factor $S_{*}\nabla^{(\nu)}R_{ijkl}\nabla^i\phi_{u+1}$
are precisely the factors $\nabla\phi'_1,\dots ,\nabla\phi'_b$. We
also recall that $\vec{\kappa}^z_{ref-doub}$ stands for the
$(u+1,\mu-1)$-refined double character of the contractions in
$\Sum_{i_h\in I_{*,l}} \tilde{C}^{l,i_1\dots i_\mu}_{g}
(\Omega_1,\dots ,\Omega_p,\phi_1,\dots
,\phi_u)\nabla_{i_1}\phi_{u+1}$. Equation (\ref{mayrh3/4})
 can then be re-written in the form:

\begin{equation}
\label{mayrh3}
\begin{split}
&\Sum_{z\in Z'_{Max}}\Sum_{l\in L^z} a_l \Sum_{i_r\in I_{*,l}}
\tilde{C}^{l,i_1\dots i_\mu}_{g} (\Omega_1,\dots
,\Omega_p,\phi_1,\dots
,\phi_u)\nabla_{i_r}\phi_{u+1}\nabla_{i_1}\upsilon
\dots\hat{\nabla}_{i_r}\upsilon\dots
\\&\nabla_{i_\mu}\upsilon
+\Sum_{h\in H} a_h Xdiv_{i_{\mu +1}}
C^{h,i_1\dots i_{\mu +1}}_{g}(\Omega_1,\dots ,\Omega_p,
\phi_1,\dots,\phi_u)\nabla_{i_1}\phi_{u+1}
\nabla_{i_2}\upsilon\dots\nabla_{i_\mu}\upsilon
\\&=\Sum_{t\in T} a_t C^{t,i_1\dots i_\mu}_{g}(\Omega_1,\dots ,\Omega_p,
\phi_1,\dots,\phi_u)\nabla_{i_1}\phi_{u+1}
\nabla_{i_2}\upsilon\dots\nabla_{i_\mu}\upsilon,
\end{split}
\end{equation}
where each $C^{t,i_1\dots i_\mu}_{g}(\Omega_1,\dots ,\Omega_p,
\phi_1,\dots,\phi_u)\nabla_{i_1}\phi_{u+1}$ is (simply or doubly)
subsequent to $\vec{\kappa}^z_{ref-doub}$. Moreover, if some
$C^{t,i_1\dots i_\mu}_{g}$ is doubly subsequent to
$\vec{\kappa}^z_{ref-doub}$ then at least one of the indices
${}_{r_1},\dots ,{}_{r_\nu},{}_j$ in the A-crucial factor
$S_{*}\nabla^{(\nu)}_{r_1\dots r_n}R_{ijkl}$ is neither
contracting against a factor $\nabla\upsilon$ nor against a factor
$\nabla\phi_h$. The complete contractions on the right hand side
arise by indexing together the contractions in
$\overline{N},\overline{J}$ in (\ref{mayrh3/4}).

 We then again assume with no loss of
 generality that in each tensor field in the first line above, the indices
${}_{i_1},\dots ,{}_{i_{c+1}}$ belong to the
 A-crucial factor and the indices ${}_{i_{c+2}},\dots ,{}_{i_\mu}$
 do not. We will then again use a weakened version of
(\ref{mayrh3}).

{\it Weakened version of (\ref{mayrh3}):} Now, we return to
(\ref{mayrh3}). We derive an equation:

\begin{equation}
\label{weakmayrh3}
\begin{split}
&\Sum_{z\in Z'_{Max}}Xdiv_{i_{c+2}}\dots Xdiv_{i_\mu} \Sum_{l\in
L^z} a_l \Sum_{i_r\in I_{*,l}}\tilde{C}^{l,i_1\dots i_\mu}_{g}
(\Omega_1,\dots ,\Omega_p,\phi_1,\dots
,\phi_u)
\\&\nabla_{i_r}\phi_{u+1}\nabla_{i_1}\upsilon
\dots\hat{\nabla}_{i_r}\upsilon\dots\nabla_{i_{c+1}}\upsilon+\Sum_{t\in T'} a_t Xdiv_{i_{c+2}}\dots
\\& Xdiv_{i_\mu}
C^{t,i_1\dots i_\mu}_{g}(\Omega_1,\dots ,\Omega_p,
\phi_1,\dots,\phi_u)\nabla_{i_1}\phi_{u+1}
\nabla_{i_2}\upsilon\dots\nabla_{i_{c+1}}\upsilon +\Sum_{h\in H} a_h Xdiv_{i_{c+2}}
\\&\dots Xdiv_{i_{\mu+1}} C^{h,i_1\dots i_{\mu +1}}_{g}(\Omega_1,\dots
,\Omega_p, \phi_1,\dots,\phi_u)\nabla_{i_1}\phi_{u+1}
\nabla_{i_2}\upsilon\dots\nabla_{i_{c+1}}\upsilon \\&+ \Sum_{j\in
J} a_j C^{j,i_1\dots i_\mu}_{g}(\Omega_1,\dots ,\Omega_p,
\phi_1,\dots,\phi_u)\nabla_{i_1}\phi_{u+1}
\nabla_{i_2}\upsilon\dots\nabla_{i_{c+1}}\upsilon=0;
\end{split}
\end{equation}
here each of the tensor fields indexed in $T'$ are doubly
subsequent to the maximal refined double characters
$\vec{\kappa}^z_{ref-doub}$, while each of the complete
contractions indexed in $J$ is simply subsequent to
$\vec{\kappa}'_{simp}$. This just follows from the previous
equation by making the factors $\nabla\upsilon$ that are not
contracting against the A-crucial factor into $Xdiv$'s (as in the
previous case A--we are applying the last
 Lemma in the Appendix of \cite{alexakis1} here).

\par Now, similarly to the previous case,
we complete our proof in two steps; we first introduce some notation: 
We divide the index set $H$ into two subsets:
We say $h\in H_1$ if at least one of the indices ${}_{r_1},\dots
,{}_{r_m},{}_j$ in the A-crucial factor
$S_{*}\nabla^{(m)}_{r_1\dots r_m}R_{ijkl}$ does not
 contract against a factor $\nabla\upsilon$ or $\nabla\phi_h,h\le u$. If all
 the indices ${}_{r_1},\dots ,{}_{r_m},{}_j$ in the
 A-crucial factor $S_{*}\nabla^{(m)}_{r_1\dots r_m}R_{ijkl}$
 contract against a factor $\nabla\upsilon$ or $\nabla\phi_f$, $f\le u$,
 we say that $h\in H_2$. Now, 
step 1 involves 
getting rid of the terms indexed in $H_2$. 
\newline

{\it Step 1:} For each $h\in H_2$, recall (from definition \ref{Xstar}) that 
 $X_{*} div_i[\dots]$ the sublinear combination in $Xdiv_i[\dots]$,
where we impose the extra restriction that $\nabla_i$ is not
allowed to hit the A-crucial factor $S_{*}\nabla^{(m)}R_{ijkl}$.
Then, since (\ref{mayrh3})  holds formally, we deduce that:

\begin{equation}
\label{weakmayrh7} \begin{split} &\Sum_{h\in H_2} a_h
X_{*}div_{i_{c+2}}\dots X_{*}div_{i_\mu} X_{*} div_{i_{\mu +1}}
C^{h,i_1\dots i_{\mu +1}}_{g}(\Omega_1,\dots ,\Omega_p,
\phi_1,\dots,\phi_u)\\&\cdot\nabla_{i_1}\phi_{u+1}
\nabla_{i_2}\upsilon\dots\nabla_{i_{c+1}}\upsilon +\Sum_{j\in J}
a_j C^{j}_{g}(\Omega_1,\dots ,\Omega_p, \phi_1,\dots,\phi_u)=0.
\end{split}
\end{equation}

\par We then claim that (\ref{weakmayrh7}) will imply a new equation,
for which we will need some more notation: Let us consider the
tensor fields
$$\tilde{C}^{l,i_1\dots i_\mu}_{g} (\Omega_1,\dots
,\Omega_p,\phi_1,\dots
,\phi_u)\nabla_{i_r}\phi_{u+1}\nabla_{i_1}\upsilon
\dots\hat{\nabla}_{i_r}\upsilon\dots\nabla_{i_{c+1}}\upsilon,$$
\\$C^{t,i_1\dots i_\mu}_{g} (\Omega_1,\dots ,\Omega_p,\phi_1,\dots
,\phi_u)\nabla_{i_r}\phi_{u+1}\nabla_{i_1}\upsilon
\dots\hat{\nabla}_{i_r}\upsilon\dots\nabla_{i_{c+1}}\upsilon$ in
(\ref{weakmayrh3}).

\par For each $l\in L^z$, we denote by
$\tilde{C}^{l,i_1i_{c+2}\dots i_\mu}_{g} (\Omega_1,\dots
,\Omega_p,\phi_{b+1},\dots,\phi_u)\nabla_{i_1}\phi_{u+1}$ the
tensor field that arises from $\tilde{C}^{l,i_1\dots
i_\mu}_{g}\nabla_{i_h}\phi_{u+1}$ (as it appears in
(\ref{weakmayrh3})) by replacing the A-crucial factor
\begin{equation}
\label{servia} S_{*}\nabla^{(\nu)}_{i_2\dots i_{c+1} l_1\dots l_b
y_{\pi+1}\dots y_\nu}R_{iy_{\nu+1}kl}
\nabla^i\phi_{u+1}\nabla^{l_1}\phi_2\dots\nabla^{l_b}\phi_b
\nabla^{i_2}\upsilon\dots\nabla^{i_{c+1}}\upsilon
\end{equation}
 (${}_{i_2},\dots ,{}_{i_\pi}$ are the free indices
that belong to that A-crucial factor) by
\\$S_{*}\nabla^{(\nu-\pi-b+1)}_{y_{\pi+1}\dots
y_\nu}R_{iy_{\nu+1}kl}\nabla^i \phi_{u+1}$. We analogously define
\\$C^{t,i_1i_{c+2}\dots i_\mu}_{g} (\Omega_1,\dots
,\Omega_p,\phi_{b+1},\dots,\phi_u)\nabla_{i_1}\phi_{u+1}$. Notice
these constructions
 are well-defined, since we know that at least one of the indices
${}_{i_2},\dots ,{}_{r_{\nu+1}}$ in the left hand side of
(\ref{servia})
 are not contracting against any factor $\nabla\phi$ or $\nabla\upsilon$.

\par Furthermore, observe that the tensor fields constructed above are acceptable, and have a given
$(u-b)$-simple character, which we denote by
$\vec{\kappa}''_{simp}$.

\par Now, our claim is that assuming (\ref{weakmayrh7}), there is a linear combination
of acceptable tensor fields (indexed in $\tilde{H}$ below) with a
simple character $\vec{\kappa}''_{simp}$, and each with rank
$\mu-c>\mu-c-1$ so that:

\begin{equation}
\label{weakmayrh7berk}
\begin{split}
&\Sum_{z\in Z'_{Max}} \Sum_{l\in
L^z} a_l \Sum_{i_r\in I_{*,l}}Xdiv_{i_{c+2}}\dots X div_{i_\mu} \tilde{C}^{l,i_{c+2}\dots
i_{\mu}}_{g} (\Omega_1,\dots ,\Omega_p,\phi_{b+1},\dots
,\phi_u)\\&\nabla_{i_r}\phi_{u+1}+
\Sum_{t\in
T'} a_t  Xdiv_{i_{c+2}}\dots X div_{i_\mu}C^{t,i_1i_{c+2}\dots
i_{\mu}}_{g} (\Omega_1,\dots ,\Omega_p,\phi_{b+1},\dots ,\phi_u)+
\\&\Sum_{h\in \tilde{H}} a_h Xdiv_{i_{c+2}}\dots Xdiv_{i_{\mu+1}}
C^{h,i_{c+2}\dots i_{\mu+1}}_{g}(\Omega_1,\dots ,\Omega_p,
\phi_{b+1},\dots,\phi_u)\nabla_{i_1}\phi_{u+1}
\\& +\Sum_{j\in J} a_j C^j_{g} (\Omega_1,\dots ,\Omega_p,\phi_{b+1},\dots
,\phi_u)=0;
\end{split}
\end{equation}
here the contractions indexed in $J$ are simply subsequent to $\vec{\kappa}''_{simp}$.
The above holds modulo contractions of length $\ge\sigma+u-b+1$.
This claim will be reduced to Lemma \ref{vanderbi} in the next subsection. We now take it for granted
and check how Proposition \ref{giade} in case IIB follows from (\ref{weakmayrh7berk}).
\newline

{\it Step 2: Derivation of Proposition \ref{giade} in case II subcase B
from (\ref{weakmayrh7berk}):} Denote the refined double
characters of the tensor fields in $\bigcup_{z\in Z'_{Max}}L^z$ by
$\vec{\kappa}'^z_{ref-doub}$; observe that the tensor fields
$C^{t,i_1i_{c+2}\dots i_\mu}_{g} (\Omega_1,\dots
,\Omega_p,\phi_{b+1},\dots,\phi_u)\nabla_{i_1}\phi_{u+1}$ are
doubly subsequent to  the refined double characters
$\vec{\kappa}'^z_{ref-doub}$. Moreover, the refined double
characters $\vec{\kappa}'^z_{ref-doub}$ are then all maximal.

\par Now, the above falls under the inductive assumption of
Proposition \ref{giade}:\footnote{Observe that the tensor fields of minimum rank in
(\ref{weakmayrh7berk}) will not contain special free indices in factors
$S_{*}\nabla^{(\nu)}R_{ijkl}$ (since
we are considering (\ref{assumpcion}) in the setting of case IIB).
Therefore there is no danger of falling under a ``forbidden case'' of
 Proposition \ref{giade}.} If $b+c>0$ then the weight of the above
complete contractions is $>-n$, and if $b+c=0$ then we have $u+1$
factors $\nabla\phi$.  Thus we derive that for each $z\in
Z'_{Max}$ there is a linear combination of acceptable tensor
fields with a refined double character
$\vec{\kappa}'^z_{ref-doub}$ (indexed in $H^z$ below) so that:

\begin{equation}
\label{weakmayrh4}
\begin{split}
&\Sum_{l\in L^z} a_l \Sum_{i_h\in
I_{*,l}}\tilde{C}^{l,i_hi_{c+2}\dots i_\mu}_{g} (\Omega_1,\dots
,\Omega_p,\phi_{b+1},\dots ,\phi_u)\nabla_{i_h}\phi_{u+1}
\nabla_{i_{c+2}}\upsilon\dots
\nabla_{i_\mu}\upsilon-
\\& Xdiv_{i_{\mu+1}}
\Sum_{h\in H} a_h C^{h,i_1i_{c+2}\dots i_{\mu
+1}}_{g}(\Omega_1,\dots ,\Omega_p,
\phi_{b+1},\dots,\phi_u)\nabla_{i_{1}}\phi_{u+1}
\nabla_{i_{c+2}}\upsilon\dots\nabla_{i_\mu}\upsilon\\&=\Sum_{t\in
T} a_t C^{t,i_{c+1}\dots i_\mu}_{g}(\Omega_1,\dots
,\Omega_p,\phi_{b+1},\dots ,\phi_u)\nabla_{i_{c+1}}\phi_{u+1}
\nabla_{i_{c+2}}\upsilon\dots \nabla_{i_\mu}\upsilon;
\end{split}
\end{equation}
(here each $C^{t,i_{c+1}\dots i_\mu}_{g}(\Omega_1,\dots
,\Omega_p,\phi_{b+1},\dots ,\phi_u)\nabla_{i_{c+1}}\phi_{u+1}$
is (simply or doubly) subsequent to
$\vec{\kappa}'^{z}_{ref-doub}$).

\par Now, we define an operation $Add[\dots ]$ that acts
 on the complete contractions and vector fields in the above
 by adding $c$
 derivative indices $\nabla_{g_1},\dots ,\nabla_{g_{c}}$ to the A-crucial
 factor and then contracting them against $c$ factors
$\nabla\upsilon$, and then adds $b$ derivative indices
$\nabla_{f_1},\dots \nabla_{f_b}$ onto the A-crucial factor and
contracts them against factors $\nabla^{f_1}\phi_1,\dots
,\nabla^{f_b}\phi_b$.

Since (\ref{weakmayrh4}) holds formally, we derive that for each
$z\in Z'_{Max}$:

\begin{equation}
\label{weakmayrh4'}
\begin{split}
&\Sum_{l\in L^z} a_l \Sum_{i_h\in
I_{*,l}}Add[C^{l,i_hi_{c+2}\dots i_{\mu }}_{g} (\Omega_1,\dots
,\Omega_p,\phi_{b+1},\dots ,\phi_u)\nabla_{i_h}\phi_{u+1}
\nabla_{i_{c+2}}\upsilon\dots \nabla_{i_\mu}\upsilon]
\\&=\{ Xdiv_{i_{\mu+1}}\Sum_{h\in H} a_h
Add[C^{h,i_{c+2}\dots i_{\mu +1}}_{g}(\Omega_1,\dots ,\Omega_p,
\phi_{b+1},\dots,\phi_u)\nabla_{i_1}\phi_{u+1}+
\\&\Sum_{t\in T} a_t Add[C^{t,i_{c+2}\dots i_{\mu-1}}_{g}(\Omega_1,\dots
,\Omega_p,\phi_{b+1},\dots,\phi_{u+1})\}
\nabla_{i_{c+2}}\upsilon\dots\nabla_{i_{\mu-1}}\upsilon],
\end{split}
\end{equation}
where each complete contraction indexed in  $T$ is simply or
doubly subsequent to $\tilde{\kappa}_{ref-doub}^z$. If we set
$\phi_{u+1}=\upsilon$ in the above, and we observe that:

\begin{equation}
\label{candles}
\begin{split}
&\Sum_{l\in L^z} a_l \Sum_{i_h\in
I_{*,l}}Add[\tilde{C}^{l,i_1i_{c+2}\dots i_{\mu}}_{g}
(\Omega_1,\dots ,\Omega_p,\phi_{b+1},\dots
,\phi_u)\nabla_{i_h}\upsilon \nabla_{i_{c+2}}\upsilon\dots
\nabla_{i_\mu}\upsilon]=
\\& |I_{*,l}|\cdot \Sum_{l\in L^z} a_l C^{l,i_1\dots i_\mu}_{g}
(\Omega_1,\dots ,\Omega_p,\phi_1,\dots ,\phi_u)
\nabla_{i_1}\upsilon\dots \nabla_{i_\mu}\upsilon+ \\&\Sum_{t\in T''}
a_t C^{t,i_1\dots
i_\mu}_{g}\nabla_{i_1}\upsilon\dots\nabla_{i_\mu}\upsilon,
\end{split}
\end{equation}
where the complete contractions indexed in $T''$ are acceptable
and doubly subsequent to $\vec{\kappa}^z_{ref-doub}$. We thus
derive that the vector field needed for Proposition \ref{giade} in
this case is precisely:

$$\frac{1}{|I_{*,l}|}\Sum_{l\in L^z} a_l \Sum_{i_h\in
I_{*,l}}Add[\tilde{C}^{l,i_1i_{\pi+1}\dots i_{\mu +1}}_{g}
(\Omega_1,\dots ,\Omega_p,\phi_{b+1},\dots
,\phi_u)\nabla_{i_h}\upsilon \nabla_{i_{c+2}}\upsilon\dots
\nabla_{i_\mu}\upsilon].$$

Therefore we have shown that Lemma \ref{pool2} implies 
Proposition \ref{giade} in case IIB, provided we can prove 
(\ref{weakmayrh7berk}).
\newline

We now show how  the claims (\ref{guilleminberk}) and
(\ref{weakmayrh7berk}) follows from four Lemmas, 
\ref{obote}, 
\ref{vanderbi} and \ref{obote3}, \ref{vanderbi3}, which we will 
state below. These four Lemmas will be derived in 
the paper \cite{alexakis5} in this series.

\subsection{Reduction of the claims (\ref{guilleminberk}) and
(\ref{weakmayrh7berk}) to the Lemmas \ref{obote}, 
\ref{vanderbi} and \ref{obote3}, \ref{vanderbi3} below.}

{\bf Reduction of claim (\ref{guilleminberk}) to Lemma
\ref{obote}:}

\par Since our Lemma \ref{obote} will also be used in other instances in this series of papers,
we will re-write our hypothesis (equation \ref{guillemin}) in
slightly more general notation:

 We will set
$c+1=\pi$ and write $\alpha$ instead of $\mu$, to stress that our
Lemma \ref{obote} below is independent of the specific values of
the parameters $\mu,c$. Furthermore, with no loss of generality,
we will assume further down that $b=0$ (in other words that there are no 
factors $\nabla\phi'_h$ contracting agianst the crucial factor--this can be done since
we can just re-name the
factors $\nabla\phi'_h$ that contract against the $A$-crucial
factor and make them into $\nabla\upsilon$s). Now, recall the
operation introduced in Step 2 after (\ref{guillemin}), where for each $h\in H_2$ we obtain
tensor fields $C^{h,i_{\pi+1}\dots i_{\mu +1}}_{g}(\Omega_1,\dots
, \Omega_p,Y, \phi_{b+1},\dots,\phi_u)$ by formally replacing the
expression $S_{*}\nabla^{(\nu)}_{r_1\dots
r_\nu}R_{ijkl}\nabla^{r_1}\upsilon\dots\nabla^j\upsilon\nabla^i\tilde{\phi}_1\nabla^k\phi_{u+1}$
by an expression $\nabla_lY$. As we noted after (\ref{guillemin}),
if we apply this operation to a true equation, we again obtain a
true equation. Thus, applying this operation to (\ref{guillemin})
we derive a new equation:

\begin{equation}
\label{guillemin2}
\begin{split}
& \Sum_{h\in H_2} a_h X_{*} div_{i_{\pi+1}}\dots X_{*} div_{i_{\mu
+1}} C^{h,i_{\pi+1}\dots i_{\mu +1}}_{g}(\Omega_1,\dots ,
\Omega_p,Y, \phi_{b+1},\dots,\phi_u)=
\\& \Sum_{j\in J} a_j
C^j_{g}(\Omega_1,\dots ,\Omega_p, \phi_{b+1},\dots,\phi_u),
\end{split}
\end{equation}
where all complete contractions and tensor fields in the above
have $\sigma+u-b-c$ factors, and are in the form:

\begin{equation}
\label{form2'}
\begin{split}
&pcontr(\nabla^{(m_1)}R_{ijkl}\otimes\dots\otimes\nabla^{m_{\sigma_1}}
R_{ijkl}\otimes
\\&S_{*}\nabla^{(\nu_1)}R_{ijkl}\otimes\dots\otimes
S_{*}\nabla^{(\nu_t)} R_{ijkl}\otimes \nabla Y\otimes
\\& \nabla^{(b_1)}\Omega_1\otimes\dots\otimes
\nabla^{(b_p)}\Omega_p\otimes
\\& \nabla\phi_{z_1}\dots \otimes\nabla\phi_{z_w}\otimes\nabla
\phi'_{z_{w+1}}\otimes
\dots\otimes\nabla\phi'_{z_{w+d}}\otimes\dots \otimes
\nabla\tilde{\phi}_{z_{w+d+1}}\otimes\dots\otimes\nabla\tilde{\phi}_{z_{w+d+y}}).
\end{split}
\end{equation}
(Notice this is the same as the form (\ref{form2}), but for the
fact that we have inserted a factor $\nabla Y$ in the second line).

\begin{definition}
 \label{XstarY}
\par In the setting of (\ref{guillemin})
$X_{*}div_i$ will stand for the sublinear combination in $Xdiv_i$
with the additional restriction that $\nabla_i$ is not allowed to
hit the factor $\nabla Y$. Moreover, we observe that the complete
contractions  in (\ref{guillemin2})
 have weight $-n +2(b+c)$.
\end{definition}

\par Some language conventions  before our next claim: We will be considering tensor
fields $C^{i_1\dots i_a}_{g}(\Omega_1,\dots ,\Omega_p,
Y,\phi_1,\dots, \phi_u)$ in the form (\ref{form2'}),
 and even more generally in the form:

\begin{equation}
\label{form2''}
\begin{split}
&pcontr(\nabla^{(m_1)}R_{ijkl}\otimes\dots\otimes\nabla^{(m_{\sigma_1})}
R_{ijkl}\otimes
\\&S_{*}\nabla^{(\nu_1)}R_{ijkl}\otimes\dots\otimes
S_{*}\nabla^{(\nu_t)} R_{ijkl}\otimes \nabla^{(B)} Y\otimes
\\& \nabla^{(b_1)}\Omega_1\otimes\dots\otimes
\nabla^{(b_p)}\Omega_p\otimes
\\& \nabla\phi_{z_1}\dots \otimes\nabla\phi_{z_w}\otimes\nabla
\phi'_{z_{w+1}}\otimes
\dots\otimes\nabla\phi'_{z_{w+d}}\otimes\dots \otimes
\nabla\tilde{\phi}_{z_{w+d+1}}\otimes\dots\otimes\nabla\tilde{\phi}_{z_{w+d+y}}).
\end{split}
\end{equation}
(Notice this only differs from (\ref{form2'}) by the fact that we
allow $B\ge 1$ derivatives on the function $Y$).

 We will say that the tensor field in the form
 (\ref{form2''}) is {\it acceptable} if all
its factors are acceptable when we {\it disregard} the factor
$\nabla^{(B)}Y$ (i.e.~ we may have $B=1$ derivatives on $Y$ but
the tensor fields will still be considered acceptable). Also, we
will still use the notion of a {\it simple character} for such tensor fields (where
 we again just disregard the factor $\nabla^{(B)}Y$). With
this convention, it follows that all the tensor fields in
(\ref{weakmayrhb}) have the same simple character, which
 we will denote by $\vec{\kappa}_{simp}'$. For such complete
 contractions $\sigma$ will stand for the number of factors
in one of the forms  $\nabla^{(m)}R_{ijkl}$,
$S_{*}\nabla^{(\nu)}R_{ijkl}$, $\nabla^{(p)}\Omega_h$,
$\nabla^{(B)}Y$.

We now state a technical Lemma, which will be proven in the next
paper in this series.

\begin{lemma}
\label{obote} Assume an equation:

\begin{equation}
\label{guillemin3}
\begin{split}
& \Sum_{h\in H_2} a_h X_{*} div_{i_{1}}\dots X_{*} div_{i_{a_h}}
C^{h,i_{1}\dots i_{a_h}}_{g}(\Omega_1,\dots , \Omega_p,Y,
\phi_{1},\dots,\phi_{u'})=
\\& \Sum_{j\in J} a_j
C^j_{g}(\Omega_1,\dots ,\Omega_p, \phi_{1},\dots,\phi_{u'}),
\end{split}
\end{equation}
where all tensor fields have rank $a_h\ge \alpha$. All tensor
fields have a given $u$-simple character $\vec{\kappa}'_{simp}$,
for which $\sigma\ge 4$. Moreover, we assume that if we formally
treat the factor $\nabla Y$ as a factor $\nabla\phi_{u'+1}$ in the
above equation, then the inductive assumption of Proposition
\ref{giade} can be applied.

The conclusion (under various assumptions which we will explain
below):  Denote by $H_{2,\alpha}$ the index set of tensor fields
with rank $\alpha$.

 We claim that there
is a linear combination of acceptable\footnote{``Acceptable'' in
the sense given after (\ref{form2''}).} tensor fields, $\Sum_{d\in
D} a_d C^{d,i_{1}\dots i_{\alpha +1}}_{g}(\Omega_1,\dots ,
\Omega_p,Y,\phi_{1},\dots, \phi_u)$, each with a simple character
$\vec{\kappa}'_{simp}$ so that:

\begin{equation}
\begin{split}
\label{vaskonik} &\Sum_{h\in H_{2,\alpha}} a_h C^{h,i_{1}\dots
i_{\alpha }}_{g}(\Omega_1,\dots ,\Omega_p, Y,\phi_{1},\dots,
\phi_{u'})\nabla_{i_{1}}\upsilon\dots \nabla_{i_{\alpha}}\upsilon-
\\& X_{*}div_{i_{\alpha+1}}\Sum_{d\in D} a_d
C^{d,i_{1}\dots i_{\alpha +1}}_{g}(\Omega_1,\dots ,\Omega_p,
Y,\phi_{1},\dots, \phi_{u'})\nabla_{i_{1}}\upsilon\dots
\nabla_{i_{\alpha}}\upsilon=
 \\&+\Sum_{t\in T}
a_t C^t_{g}(\Omega_1,\dots , \Omega_p
,Y,\phi_{1},\dots,\phi_{u'},\upsilon^{\alpha}).
\end{split}
\end{equation}
 The linear combination on the right hand side stands for a generic
linear combination of complete contractions in the form (\ref{form2'}) with a factor $\nabla
Y$ and with a simple character that is
 subsequent to $\vec{\kappa}'_{simp}$.
\newline

{\it The assumptions under which (\ref{vaskonik}) will hold:} The
assumption under which (\ref{vaskonik}) holds is that there should be no tensor fields 
 of rank $\alpha$ in (\ref{guillemin3}) which are ``bad''. Here ``bad'' means the following: 

If $\sigma_2=0$ in $\vec{\kappa}'_{simp}$ then a tensor field in the form (\ref{form2'}) 
is ``bad'' provided:
\begin{enumerate}
 \item{The factor $\nabla Y$ contains a free index.}
\item{If we formally erase the factor $\nabla Y$ (which contains 
a free index), then the resulting tensor field should have no removable indices,\footnote{Thus,
the tensor field should consist of factors $S_*R_{ijkl},\nabla^{(2)}\Omega_h$,
 and factors $\nabla^{(m)}_{r_1\dots r_m}R_{ijkl}$
 with all the indices ${}_{r_1},\dots ,{}_{r_m}$ 
contracting against factors $\nabla\phi_h$.} and no free 
indices.\footnote{I.e.~$\alpha=1$ in (\ref{guillemin3}).}
Moreover, any factors $S_{*}R_{ijkl}$ should be {\it simple}. }
\end{enumerate}

If $\sigma_2>0$ in $\vec{\kappa}'_{simp}$ then a tensor field in the form (\ref{form2'}) 
is ``bad'' provided:
\begin{enumerate}
 \item{The factor $\nabla Y$ should contain a free index.}
\item{If we formally erase the factor $\nabla Y$ (which contains 
a free index), then the resulting tensor field should have no removable indices, 
any factors $S_{*}R_{ijkl}$ should be {\it simple}, any
 factor $\nabla^{(2)}_{ab}\Omega_h$ should 
have at most one of the indices 
${}_a,{}_b$ free or contracting 
against  a factor $\nabla\phi_s$.}
\item{Any factor $\nabla^{(m)}R_{ijkl}$ can contain at most 
one (necessarily special, by virtue of 2.) free index. }
\end{enumerate}

  Furthermore, we claim that the proof of
 this Lemma will only rely on the inductive assumption of
 Proposition \ref{giade}. Moreover, we claim that if all the tensor
 fields indexed in $H_2$ (in (\ref{guillemin3})) do not have a free index in $\nabla Y$
 then we may assume that the tensor fields indexed in $D$ in (\ref{vaskonik}) have
 the same property.
\end{lemma}

{\it Note:} It follows (by weight considerations) 
that none of the tensor fields of minimum rank in (\ref{guillemin2}) 
is ``bad'' in the above sense, since our assumption
(\ref{assumpcion}) {\it does not} fall under one of the special
cases, as described in the beginning of this subsection. 
\newline

\par We also claim a corollary of Lemma \ref{obote}. Firstly, we
introduce some notation:

$$\Sum_{q\in Q} a_q C^{q,i_{\pi+1}\dots i_{\mu +1}}_{g}
(\Omega_1,\dots ,\Omega_p, Y, \phi_{b+1},\dots ,\phi_u)$$ will
stand for a generic linear combination of acceptable tensor fields
with a simple character $\vec{\kappa}'_{simp}$ and with a factor
$\nabla^{(B)}Y$ with $B\ge 2$ (and where this factor is not
 contracting against any factors $\nabla\phi_h$).

\begin{corollary}
\label{obotecor} Under the assumptions of Lemma \ref{obote},
 with $\sigma\ge 4$ we  can write:
\begin{equation}
\begin{split}
\label{vaskonik} &\Sum_{h\in H_2} a_h Xdiv_{i_{1}}\dots
Xdiv_{i_{\alpha}}C^{h,i_{1}\dots i_{\alpha}}_{g}(\Omega_1,\dots
,\Omega_p, Y,\phi_{1},\dots, \phi_{u'})=
\\&\Sum_{q\in Q} a_q Xdiv_{i_{1}}\dots
Xdiv_{i_{a'}} C^{q,i_{1}\dots i_{a'}}_{g} (\Omega_1,\dots
,\Omega_p, Y,\phi_{1},\dots,\phi_{u'})
\\&+\Sum_{t\in T} a_t C^t_{g}(\Omega_1,\dots , \Omega_p
,Y,\phi_{1},\dots,\phi_{u'}),
\end{split}
\end{equation}
where the linear combination 
$\Sum_{q\in Q} a_q C^{q,i_{1}\dots
i_{a'}}_{g}$ stands
 for a generic linear combination of tensor fields
  in the form (\ref{form2''}) with $B\ge 2$, with a
 simple character $\vec{\kappa}'_{simp}$ and with each $a'\ge \alpha$. The acceptable complete
  contractions \\$C^t_{g}(\Omega_1,\dots , \Omega_p
 ,Y,\phi_{1},\dots,\phi_{u'})$ are simply subsequent to
 $\vec{\kappa}'_{simp}$. $Xdiv_i$ here means that $\nabla_i$ is
 not allowed to hit the factors $\nabla\phi_h$ (but it is
 allowed to hit $\nabla^{(B)}Y$).
\end{corollary}

\par We have an analogue of the above corollary when $\sigma=3$
 (the next Lemma, \ref{obote3}, will also be proven in the 
paper \cite{alexakis5}).

\begin{lemma}
\label{obote3} We assume (\ref{guillemin3}), where $\sigma=3$. We
also assume that for each of the tensor fields in
$H_2^{\alpha,*}$\footnote{Recall that $H_2^{\alpha,*}$
 is the index set of tensor fields of rank $\alpha$ in (\ref{guillemin3})
 with a free index in the factor $\nabla Y$.} there is at least
 one removable index. We then have two
claims:

Firstly, the conclusion of Lemma \ref{obote}
 holds in this setting also. Secondly, the conclusion
of Corollary \ref{obotecor} is true in this setting.
\end{lemma}

Before we show that  Corollary \ref{obotecor} follows from Lemma
\ref{obote}, let us see how our desired equation
(\ref{guilleminberk}) follows from the above corollary:
\newline

{\it Corollary \ref{obotecor} (or Lemma \ref{obote3} when
$\sigma=3$) implies (\ref{guilleminberk}):} We introduce an
operation $Op\{\dots\}$ which acts on complete contractions and
tensor fields in the form (\ref{form2''}) by formally replacing
the factor $\nabla^{(B)}_{r_1\dots r_B}Y$ (recall $B\ge 1$) by
$$S_{*}\nabla^{(B+b+c-2)}_{y_1\dots y_b s_1\dots
s_cr_1\dots r_{B-2}}R_{ir_{B-1}s_{c+1}r_B}\nabla^{y_1}\phi_1\dots
\nabla^{y_b}\phi_b\nabla^{s_1}\upsilon\dots
\nabla^{s_{c+1}}\upsilon.$$ Then, if we apply this operation to
(\ref{vaskonik}) and we repeat the permutations by which we make
(\ref{vaskonik}) formally zero (modulo introducing correction
terms by the Bianchi identities (\ref{koichi1}), (\ref{koichi2}),
(\ref{koichi3})), we derive (\ref{guilleminberk}). So matters are
reduced to showing that Corollary \ref{obotecor} follows from
Lemma \ref{obote}. $\Box$
\newline

{\it Proof that Corollary \ref{obotecor} follows from Lemma
\ref{obote}:}
\newline

\par The proof is by induction. Firstly, we apply Lemma
\ref{obote} and we pick out the sublinear combination in the
conclusion of Lemma \ref{obote} where $\nabla Y$ is contracting
against a factor $\nabla\upsilon$. That sublinear combination must
vanish separately, thus we obtain an equation:

\begin{equation}
\begin{split}
\label{vaskonikaaa} &\Sum_{h\in H_{2}^{\alpha,*}} a_h
C^{h,i_{1}\dots i_{\alpha }}_{g}(\Omega_1,\dots ,\Omega_p,
Y,\phi_{1},\dots, \phi_{u'})\nabla_{i_{1}}\upsilon\dots
\nabla_{i_{\alpha}}\upsilon-
\\& X_{*}div_{i_{\alpha+1}}\Sum_{d\in D'} a_d
C^{d,i_{1}\dots i_{\alpha +1}}_{g}(\Omega_1,\dots ,\Omega_p,
Y,\phi_{1},\dots, \phi_{u'})\nabla_{i_{1}}\upsilon\dots
\nabla_{i_{\alpha}}\upsilon=
 \\&+\Sum_{t\in T}
a_t C^t_{g}(\Omega_1,\dots , \Omega_p
,Y,\phi_{1},\dots,\phi_{u'},\upsilon^{\alpha}).
\end{split}
\end{equation}
\par Now, we make the the factors $\nabla\upsilon$ into 
$Xdiv$s\footnote{See the last Lemma in the Appendix of \cite{alexakis1}.}
(which {\it are} allowed to hit the factor $\nabla Y$) and we
derive a new equation:

\begin{equation}
\begin{split}
\label{vaskonikbbb} &\Sum_{h\in H_{2}^{\alpha,*}}
a_hXdiv_{i_1}\dots Xdiv_{i_\alpha} C^{h,i_{1}\dots i_{\alpha
}}_{g}(\Omega_1,\dots ,\Omega_p, Y,\phi_{1},\dots, \phi_{u'})-
\\& Xdiv_{i_1}\dots Xdiv_{i_\alpha}Xdiv_{i_{\alpha+1}}\Sum_{d\in D'} a_d
C^{d,i_{1}\dots i_{\alpha +1}}_{g}(\Omega_1,\dots ,\Omega_p,
Y,\phi_{1},\dots, \phi_{u'})= \\&\Sum_{q\in Q} a_q
Xdiv_{i_{1}}\dots Xdiv_{i_{\alpha}} C^{q,i_{1}\dots
i_{\alpha}}_{g} (\Omega_1,\dots ,\Omega_p,
Y,\phi_{1},\dots,\phi_{u'})
 \\&+\Sum_{j\in J}
a_j C^j_{g}(\Omega_1,\dots , \Omega_p ,Y,\phi_{1},\dots,\phi_{u'}).
\end{split}
\end{equation}

\par In view of this equation, we are reduced to proving our claim
when $H_2^{\alpha,*}=\emptyset$. That is, we may then
additionally assume that in the hypothesis of Lemma \ref{obote} no
tensor fields contain a free index in $\nabla Y$ (if there are
such tensor fields with a factor $\nabla_{i_1} Y$, we just treat
$X_{*}div_{i_1} \nabla_{i_1}Y[\dots]$ as a sum of
$\beta$-tensor fields, $\beta\ge \alpha$). We will be making
this assumption until the end of this proof.

\par Then we proceed by induction. More precisely, our inductive statement is
the following: Suppose we know that for some number $f\ge 0$
we can write:

\begin{equation}
\begin{split}
\label{foroi} &\Sum_{h\in H_2} a_h Xdiv_{i_{1}}\dots
Xdiv_{i_{a}}C^{h,i_{1}\dots i_{a}}_{g}(\Omega_1,\dots ,\Omega_p,
Y,\phi_{1},\dots,\phi_{u'})=
\\&\Sum_{q\in Q} a_q Xdiv_{i_{1}}\dots
Xdiv_{i_{a}} C^{q,i_{1}\dots i_{a}}_{g}(\Omega_1,\dots ,\Omega_p,
Y,\phi_{1},\dots,\phi_{u'})+
\\&\Sum_{h\in H^f_2} a_h Xdiv_{i_{1}}\dots
Xdiv_{i_{a+f}}C^{h,i_{1}\dots i_{a +f}}_{g}(\Omega_1,\dots
,\Omega_p, Y,\phi_{1},\dots,\phi_{u'})+
 \\&+\Sum_{t\in T}
a_t C^t_{g}(\Omega_1,\dots , \Omega_p ,Y,\phi_{1},\dots,\phi_{u'}),
\end{split}
\end{equation}
where the tensor fields indexed in $H^f_2$ still have a factor
$\nabla Y$ (which {\it does not} contain a free index) but are
otherwise acceptable with simple character $\vec{\kappa}'_{simp}$
and have rank $\alpha+f$.

\par Our claim is that we can then write:

\begin{equation}
\begin{split}
\label{foroi2}  &\Sum_{h\in H_2} a_h Xdiv_{i_{1}}\dots
Xdiv_{i_{a}}C^{h,i_{1}\dots i_{a}}_{g}(\Omega_1,\dots ,\Omega_p,
Y,\phi_{1},\dots,\phi_{u'})=
\\&\Sum_{q\in Q} a_q Xdiv_{i_{1}}\dots
Xdiv_{i_{a}} C^{q,i_{1}\dots i_{a}}_{g}(\Omega_1,\dots ,\Omega_p,
Y,\phi_{1},\dots,\phi_{u'})+
\\&\Sum_{h\in H^{f+1}_2} a_h Xdiv_{i_{1}}\dots
Xdiv_{i_{a+f+1}}C^{h,i_{1}\dots i_{a +f+1}}_{g}(\Omega_1,\dots
,\Omega_p, Y,\phi_{1},\dots,\phi_{u'})+
 \\&+\Sum_{t\in T}
a_t C^t_{g}(\Omega_1,\dots , \Omega_p ,Y,\phi_{1},\dots,\phi_{u'});
\end{split}
\end{equation}
(with the same convention regarding $\sum_{h\in
H^{f+1}_2}\dots$--it is like the sublinear combination $\sum_{h\in
H^f_2}\dots$ only with rank $\ge f +1$).

\par Clearly, if we can show this inductive step then our Corollary
 will follow, since we are dealing with tensor fields of a
 {\it fixed} weight $-K, K\le n$.

\par This inductive step is not hard to deduce. Assuming
(\ref{foroi}), we pick out the sublinear combination that
 contains a  factor $\nabla Y$ (which vanishes separately)
and we replace it into (\ref{guillemin2}) to derive the equation:

\begin{equation}
\begin{split}
\label{foroi3} &\Sum_{h\in H^f_2} a_h X_{*}div_{i_{1}}\dots
X_{*}div_{i_{\alpha+f}}C^{h,i_{1}\dots i_{\alpha
+f}}_{g}(\Omega_1,\dots ,\Omega_p, Y,\phi_{1},\dots,\phi_{u'})=
 \\&\Sum_{t\in T}
a_t C^t_{g}(\Omega_1,\dots , \Omega_p ,Y,\phi_{1},\dots,\phi_{u'}).
\end{split}
\end{equation}

\par Now, applying Lemma \ref{obote} to
this equation,\footnote{Since we are assuming
 that the terms of maximal refined double character in the 
hypothesis of Proposition \ref{giade} are assumed not to be ``special'' 
(as defined in the begining of the previous subsection), it follows by weight 
considerations that no terms in (\ref{foroi3}) are ``bad'' in the language of Lemmas \ref{obote}.} 
 we derive that
there is a linear combination of acceptable $(\alpha+f+1)$-tensor
fields (indexed in $D^f$ below) with a factor $\nabla Y$ and a
simple character $\vec{\kappa}'_{simp}$ so that:

\begin{equation}
\begin{split}
\label{foroi4} &\Sum_{h\in H^f_2} a_h C^{h,i_{1}\dots i_{\alpha
+f}}_{g}(\Omega_1,\dots ,\Omega_p,
Y,\phi_{1},\dots,\phi_{u'})\nabla_{i_{1}}\upsilon\dots
\nabla_{i_{\alpha+f+1}}\upsilon-
\\& X_{*}div_{i_{\alpha+f+1}}\Sum_{d\in D^f} a_d
C^{d,i_{1}\dots i_{\alpha+f+1}}_{g}(\Omega_1,\dots ,
\Omega_p,Y,\phi_{1},\dots, \phi_{u'})\nabla_{i_{1}}\upsilon\dots
\nabla_{i_{\alpha+f}}\upsilon=
 \\&+\Sum_{t\in T}
a_t C^t_{g}(\Omega_1,\dots , \Omega_p
,Y,\phi_{1},\dots,\phi_{u'},\upsilon^{\alpha-\pi}).
\end{split}
\end{equation}

\par But observe that the above implies:

\begin{equation}
\begin{split}
\label{foroi5} &\Sum_{h\in H^f_2} a_h C^{h,i_{1}\dots i_{\alpha
+f}}_{g}(\Omega_1,\dots ,\Omega_p,
Y,\phi_{1},\dots,\phi_{u'})\nabla_{i_{1}}\upsilon\dots
\nabla_{i_{\alpha+f}}\upsilon-
\\& Xdiv_{i_{\alpha+f+1}}\Sum_{d\in D^f} a_d
C^{d,i_{1}\dots i_{\alpha+f+1}}_{g}(\Omega_1,\dots ,
\Omega_p,Y,\phi_{1},\dots, \phi_{u'})\nabla_{i_{1}}\upsilon\dots
\nabla_{i_{\alpha+f}}\upsilon=
\\& \Sum_{q\in Q} a_q
C^{q,i_{1}\dots i_{\alpha+f}}_{g}(\Omega_1,\dots ,
\Omega_p,Y,\phi_{b+1},\dots, \phi_{u'})\nabla_{i_{1}}\upsilon\dots
\nabla_{i_{\alpha+f}}\upsilon
\\&+\Sum_{t\in T}
a_t C^t_{g}(\Omega_1,\dots , \Omega_p
,Y,\phi_{1},\dots,\phi_{u'},\upsilon^{\alpha+f+1}),
\end{split}
\end{equation}
where the tensor fields indexed in $Q$ are acceptable with a
 simple character $\vec{\kappa}'_{simp}$ and with a factor
$\nabla^{(2)}Y$. But then just making the $\nabla\upsilon$'s into
 $Xdiv$'s in the above we obtain (\ref{foroi2}) with
$H^{f+1}=D^f$. $\Box$
\newline

{\it Reduction of equation (\ref{weakmayrh7berk}) to Lemmas
\ref{vanderbi}, \ref{vanderbi3} below:}
\newline

\par We define an operation that acts on the tensor fields in
(\ref{mayrh3})  and (\ref{weakmayrh7}) by replacing the expression
$$S_{*}\nabla^{(\nu+b)}_{s_1\dots s_br_1\dots
r_\nu}R_{ijkl}\nabla^i\phi_{u+1}\nabla^{r_1}\upsilon\dots
\nabla^{r_\tau}\upsilon\nabla^{s_1}\phi'_1\dots\nabla^{s_b}\phi'_b$$
by an expression $\nabla_{(r_{\tau+1}\dots
r_mj)l}\omega_1\nabla_k\omega_2-\nabla_{(r_{\tau+1}\dots
r_mj)k}\omega_1\nabla_l\omega_2$. We denote this operation by
$Repl\{\dots \}$. Thus, acting with the above operation we obtain
complete contractions and tensor fields in the form:

\begin{equation}
\label{tavuk}
\begin{split}
&contr(\nabla^{(m_1)}R_{ijkl}\otimes\dots \otimes
\nabla^{(m_s)}R_{ijkl}\otimes
\\&S_{*}\nabla^{(\nu_1)}R_{ijkl}\otimes\dots\otimes
S_{*}\nabla^{(\nu_b)}R_{ijkl}\otimes \nabla^{(B,+)}_{r_1\dots
r_B}(\nabla_a\omega_1\nabla_b\omega_2-\nabla_b\omega_1\nabla_a\omega_1)
\\&\otimes \nabla^{(d_1)}\Omega_p\otimes\dots\otimes\nabla^{(d_p)}
\Omega_p \otimes \nabla\phi_1\otimes\dots\otimes
\nabla\phi_u);
\end{split}
\end{equation}
here $\nabla^{(B,+)}_{r_1\dots r_B}(\dots )$ stands for the
sublinear combination in $\nabla^{(B)}_{r_1\dots r_B}(\dots )$
where each $\nabla$ is not allowed to hit the factor
$\nabla\omega_2$.

\begin{definition}
\label{newsimpchar}
 We define the simple character of a complete contraction
or tensor field in the above form to be the simple character of
the complete contraction or tensor fields that arises from it by
{\it disregarding} the two factors
$\nabla^{(B)}\omega_1,\nabla\omega_2$. For each tensor field in
the form (\ref{tavuk}),
 we will also define $\sigma$ to stand for the number of factors
$\nabla^{(m)}R_{ijkl},S_{*}\nabla^{(\nu)}R_{ijkl}$, $\nabla^{(p)}\Omega_h$ {\it plus one}.
(In other words, we are not counting the $\nabla\phi$'s
and we are counting the two factors $\omega_1,\omega_2$ as one).
\end{definition}

We then derive from (\ref{weakmayrh3}) that:

\begin{equation}
\label{weakmayrh9}
\begin{split}
&\Sum_{z\in Z'_{Max}}X^{+}div_{i_{c+2}}\dots X^{+}div_{i_\mu}
\Sum_{l\in L^z} a_l\Sum_{i_h\in I_{*,l}}
\\& Repl\{ C^{l,i_1\dots
i_{\mu}}_{g} (\Omega_1,\dots ,\Omega_p,\phi_1,\dots
,\phi_u)\nabla_{i_h}\phi_{u+1}\nabla_{i_1}\upsilon\dots\hat{\nabla}_{i_h}\upsilon
\dots\nabla_{i_{c+1}}\upsilon\}+
\\&\Sum_{t\in T'} a_t
X^{+}div_{i_{c+2}}\dots X^{+}div_{i_\mu} Repl\{ C^{t,i_1\dots
i_\mu}_{g}(\Omega_1,\dots ,\Omega_p,
\phi_1,\dots,\phi_u)\nabla_{i_1}\phi_{u+1}
\\&\nabla_{i_2}\upsilon\dots\nabla_{i_{c+2}}\upsilon\} -
X^{+}div_{i_{c+2}}\dots X^{+} div_{i_{\mu +1}}
 \\&\Sum_{h\in H} a_h Repl\{ C^{h,i_1\dots i_{\mu
+1}}_{g}(\Omega_1,\dots ,\Omega_p,
\phi_1,\dots,\phi_u)\nabla_{i_1}\phi_{u+1}
\nabla_{i_2}\upsilon\dots\nabla_{i_{c+1}}\upsilon\}
\\&=\Sum_{j\in J} a_j Repl\{ C^{j,i_1\dots i_{c+1}}_{g}(\Omega_1,\dots ,\Omega_p,
\phi_1,\dots,\phi_u)\nabla_{i_1}\phi_{u+1}
\nabla_{i_2}\upsilon\dots\nabla_{i_{c+1}}\upsilon\},
\end{split}
\end{equation}
where here $X^{+}div_i$ stands for the sublinear combination in
$Xdiv_i$ where $\nabla_i$ is not allowed to hit the factor
$\nabla\omega_2$. This equation follows from the {\it proof}
of the last Lemma in the Appendix in \cite{alexakis1}.\footnote{Be repeating 
exactly the same argument.} Thus, we derive that:

\begin{equation}
\label{weakmayrh8} \begin{split} &\Sum_{h\in H_2} a_h
X_{*}div_{i_{c+2}}\dots X_{*}div_{i_\mu} X_{*} div_{i_{\mu +1}}
Repl\{ C^{h,i_1\dots i_{\mu +1}}_{g}(\Omega_1,\dots ,\Omega_p,
\\&\phi_1,\dots,\phi_u)\nabla_{i_1}\phi_{u+1}
\nabla_{i_2}\upsilon\dots\nabla_{i_{c+1}}\upsilon\}=
\\&\Sum_{j\in J'} a_j Repl\{ C^{t,i_1\dots i_{c+1}}_{g}(\Omega_1,\dots ,\Omega_p,
\phi_1,\dots,\phi_u)\nabla_{i_1}\phi_{u+1}
\nabla_{i_2}\upsilon\dots\nabla_{i_{c+1}}\upsilon\},
\end{split}
\end{equation}
where here $X_{*}$ stands for the sublinear combination in
$Xdiv_i$ where $\nabla_i$ is not allowed to hit either of the
factors $\nabla\omega_1,\nabla\omega_2$. Also $J'\subset J$ stands
for the sublinear combination of complete contraction with two
factors $\nabla\omega_1,\nabla\omega_2$ (each with only one
derivative).

\par We then formulate Lemma \ref{vanderbi},  which we will 
show in \cite{alexakis5}. 
We introduce one further piece of notation
before stating this claim: 

\begin{definition}
\label{removable} 
Let $C^{x,i_1\dots
i_a}_{g}(\Omega_1,\dots
,\Omega_p,[\omega_1,\omega_2],\phi_1,\dots,\phi_{u'})$ stand for a
tensor field in the form (\ref{tavuk}) with $B=0$.
We will say that a derivative index in some
factor $\nabla^{(m)}R_{ijkl}$ or $S_{*}\nabla^{(\nu)}R_{ijkl}$ in
$C^{x,i_1\dots i_a}_{g}(\Omega_1,\dots
,\Omega_p,[\omega_1,\omega_2],\phi_1,\dots,\phi_{u'})$ is
``removable'' if it is neither free
not contracting against a factor $\nabla\phi_h$.
 \end{definition}

 Now, consider any
factor $\nabla^{(B)}_{r_1\dots r_B}\Omega_v$ in $C^{x,i_1\dots
i_a}_{g}$, where we make the normalizing requirement that all
indices that are either free or are contracting against a factor
$\nabla\phi_h$ or $\nabla\omega_f$ are pulled to the right. We
then say that an index in $\nabla^{(B)}_{r_1\dots r_B}\Omega_v$ is
``removable'' if it is one of the leftmost $B-2$ indices and it is
neither free, nor contracting against any factor
$\nabla\phi_h,\nabla\omega_f$.

\begin{lemma}
\label{vanderbi} Consider a linear combination of partial
contractions,
$$\Sum_{x\in X} a_x C^{x,i_1\dots
i_a}_{g}(\Omega_1,\dots
,\Omega_p,[\omega_1,\omega_2],\phi_1,\dots,\phi_{u'}),$$ where each
of the tensor fields $C^{x,i_1\dots
i_a}_{g}$ is in the form (\ref{tavuk}) with $B=0$ (and
is antisymmetric in the factors
$\nabla_a\omega_1,\nabla_b\omega_2$ by definition), with rank
$a\ge\alpha$ and length $\sigma\ge 4$.\footnote{ Recall we are
counting the two factors $\omega_1,\omega_2$ for one} We assume
all these tensor fields have a given simple character which we
denote by $\vec{\kappa}'_{simp}$ (we use $u'$ instead of $u$ to
stress that this Lemma holds in generality). We assume an
equation:

\begin{equation}
\label{para3enh}
\begin{split}
& \Sum_{x\in X} a_x X_{*}div_{i_1}\dots
X_{*}div_{i_a}C^{x,i_1\dots i_a}_{g}(\Omega_1,\dots
,\Omega_p,[\omega_1,\omega_2],\phi_1,\dots,\phi_u)+
\\& \Sum_{j\in J} a_j C^j_{g}(\Omega_1,\dots
,\Omega_p,[\omega_1,\omega_2],\phi_1,\dots,\phi_u)=0,
\end{split}
\end{equation}
where $X_{*}div_i$ stands for the sublinear combination in
$Xdiv_i$ where $\nabla_i$ is in addition not allowed to hit the
factors $\nabla\omega_1,\nabla\omega_2$. The contractions $C^j$
here are simply subsequent to $\vec{\kappa}'_{simp}$. We assume
that if we formally treat the factors
$\nabla\omega_1,\nabla\omega_2$ as factors
$\nabla\phi_{u+1},\nabla\phi_{u+2}$ (disregarding whether they are
contracting against special indices) in the above, then the
inductive assumption of Proposition \ref{giade} applies.

\par The conclusion we will draw (under various hypotheses that
we will explain below) is that we can write:

\begin{equation}
\label{para3enh2} \begin{split} &\Sum_{x\in X} a_x
X_{+}div_{i_1}\dots X_{+}div_{i_a}C^{x,i_1\dots
i_a}_{g}(\Omega_1,\dots
,\Omega_p,[\omega_1,\omega_2],\phi_1,\dots,\phi_u)=
\\&\Sum_{x\in X'} a_x
X_{+}div_{i_1}\dots X_{+}div_{i_a}C^{x,i_1\dots
i_a}_{g}(\Omega_1,\dots
,\Omega_p,[\omega_1,\omega_2],\phi_1,\dots,\phi_u)+
\\& \Sum_{j\in
J} a_j C^j_{g}(\Omega_1,\dots
,\Omega_p,[\omega_1,\omega_2],\phi_1,\dots,\phi_u),
\end{split}
\end{equation}
where the tensor fields indexed in $X'$ on the right hand side are
in the form (\ref{tavuk}) with $B>0$. All the other sublinear
combinations are as above. We recall that $X_{+}div_i$ stands for
the sublinear combination in $Xdiv_i$ where $\nabla_i$ is in
addition not allowed to hit the factor $\nabla\omega_2$ (it is
allowed to hit the factor $\nabla^{(B)}\omega_1$).
\newline

{\it Assumptions needed to derive (\ref{para3enh2}):} We claim
(\ref{para3enh2}) under certain assumptions on the $\alpha$-tensor
fields in (\ref{para3enh}) which have rank $\alpha$ and have a
free index in one of the factors $\nabla\omega_1,\nabla\omega_2$
(say to $\nabla\omega_1$ wlog)--we denote the index set of those
tensor fields by $X^{\alpha,*}\subset X$.

 The assumption we need in order for the claim to hold is that no
tensor field indexed in $X^{\alpha,*}$ should be ``bad''. A  tensor field 
is ``bad'' if it has the property that 
when we erase the expression 
$\nabla_{[a}\omega_1\nabla_{b]}\omega_2$ (and make 
the index that contracted against ${}_b$ into a free index) then the resulting tensor field 
will have no removable indices, and all factors $S_{*}R_{ijkl}$ will be simple.
\end{lemma}

\begin{lemma}
\label{vanderbi3} We assume (\ref{para3enh}), where now the tensor 
fields have length $\sigma=3$.
We also assume that for each of the tensor fields indexed in  $X$, there
is a removable index in each of the real factors. We then claim that the conclusion of
Lemma \ref{vanderbi} is still true in this setting.
\end{lemma}

\par We will show Lemmas \ref{vanderbi}, \ref{vanderbi3}
 in \cite{alexakis5}. For now, let us see how they
 imply (\ref{weakmayrh7berk}).

 {\it Note:} Observe
 that (\ref{weakmayrh8}) satisfies the requirements of 
 Lemma \ref{vanderbi} by  weight considerations, 
since we are assuming that (\ref{assumpcion}) {\it does not} fall
 under any of the ``special cases'' outlined in the beginning of the previous subsection. 
\newline

\par Thus, we will now apply Lemma \ref{vanderbi} (or \ref{vanderbi3})
  to (\ref{weakmayrh8}).

\par  Consider (\ref{weakmayrh8}).
 We denote by $\vec{\kappa}_{*}$ the simple character of the
tensor fields $Repl\{C^{h,i_1\dots i_{\alpha+1}}_{g}\}$. We then
observe that Lemmas \ref{vanderbi} (or \ref{vanderbi3}) imply:
\begin{equation}
\label{mayrh10}
\begin{split}
&\Sum_{h\in H_2} a_h X^{+} div_{i_{c+2}} X^{+} div_{i_{\mu +1}}
Repl\{ C^{h,i_1\dots i_{\mu +1}}_{g}(\Omega_1,\dots ,\Omega_p,
\phi_1,\dots,\phi_u)\\&\nabla_{i_1}\phi_{u+1}
\nabla_{i_2}\upsilon\dots\nabla_{i_\pi}\upsilon\}
\\&=\Sum_{q\in Q} a_q X^{+} div_{i_{c+2}} X^{+} div_{i_{\beta}}
  C^{q,i_{i_{c+2}}\dots i_\beta}_{g}(\Omega_1,\dots
,\Omega_p, \omega_1 ,\omega_2,\phi_1,\dots,\phi_u) +
\\&\Sum_{j\in J} a_j
  C^{j}_{g}(\Omega_1,\dots
,\Omega_p, \omega_1 ,\omega_2,\phi_1,\dots,\phi_u),
\end{split}
\end{equation}
where each $C^{q,i_{c+2}\dots i_\beta}_{g}(\Omega_1, \dots
,\Omega_p,\omega_1,\omega_2,\phi_1,\dots,\phi_u)$ ($\beta\ge \mu
+1$) is a generic acceptable tensor field in the
 form (\ref{tavuk}), with the additional restriction that it has an expression
$\nabla^{+}_\chi(\nabla_a\omega_1\nabla_b\omega_2-
\nabla_a\omega_2\nabla_b\omega_1)$.\footnote{Recall that
$\nabla^{+}_\chi$ stands for the sublinear combination in
$\nabla_\chi$ where $\nabla_\chi$ is not allowed to hit
$\nabla\omega_2$} Also, each \\$C^{j}_{g}(\Omega_1, \dots
,\Omega_p,\omega_1,\omega_2,\phi_1,\dots,\phi_u)$ is in the form
(\ref{tavuk}) with $(B=0)$ but is also simply  subsequent to
$\vec{\kappa}_{*}$.

\par We now define an operation $Op_{*}$ which formally acts on
the complete contractions (and linear combinations thereof) in
(\ref{weakmayrh9}) by replacing each expression
$\nabla^{(K)}_{(r_1\dots
r_K)}\omega_1\nabla_\gamma\omega_2$\footnote{${}_{(\dots)}$ stands
for symmetrization of the indices between parentheses.} by an
expression:

$$(K-1)\cdot \nabla^{(K+1)}_{\gamma (r_1\dots r_K)}\phi_{u+1}$$
 This operation can also be defined on the tensor
fields
 appearing in (\ref{weakmayrh3}). Before we proceed to explain how this
operation can act on true equations and produce true equations,
let us see what will be the outcome of formally
 applying $Op_{*}$ to the equation (\ref{weakmayrh9}):
\newline

{\it $Op_{*}$ acting on (\ref{weakmayrh9}) proves
(\ref{weakmayrh7berk}):} For each $l\in L^z$, we denote by

$$\tilde{C}^{l,i_1i_{c+2}\dots i_\mu}_{g} (\Omega_1,\dots
,\Omega_p,\phi_{b+1},\dots,\phi_u)\nabla_{i_1}\phi_{u+1}$$ the
tensor field that arises from $\tilde{C}^{l,i_1\dots
i_\mu}_{g}\nabla_{i_h}\phi_{u+1}$ (as it appears in
(\ref{weakmayrh3})) by replacing the A-crucial factor
\begin{equation}
\label{servia}
S_{*}\nabla^{(\nu)}_{i_2\dots i_{c+1} l_1\dots l_b y_{f}\dots y_\nu}R_{ir_{\nu+1}kl}
\nabla^i\phi_{u+1}\nabla^{l_1}\phi_2\dots\nabla^{l_b}\phi_b
\nabla^{i_2}\upsilon\dots\nabla^{i_{c+1}}\upsilon
\end{equation}
 (${}_{i_2},\dots ,{}_{i_\pi}$ are the free indices
that belong to that critical factor) by
\\$S_{*}\nabla^{(\nu-\pi-b+1)}_{y_{\pi+1}\dots
y_\nu}R_{ijkl}\nabla^i \phi_{u+1}$.

\par We analogously define the tensor fields $C^{h,i_1i_{c+2}\dots i_{\mu+1}}_{g}
(\Omega_1,\dots ,\Omega_p,\phi_{b+1},\dots,\phi_u)$,
$C^{t,i_1i_{c+2}\dots i_\mu}_{g} (\Omega_1,\dots
,\Omega_p,\phi_{b+1},\dots,\phi_u)\nabla_{i_1}\phi_{u+1}$. Observe
that for the tensor fields indexed in $H$, this is a well-defined
operation, since we are assuming that $H_2=\emptyset$ in
(\ref{weakmayrh9}) (thus for each of the tensor fields above we
will have that at least one of the indices ${}_{i_1},\dots
,{}_{r_{\nu+1}}$ in
 (\ref{servia}) is not contracting against a factor $\nabla\phi$ or $\nabla\upsilon$).

 We observe that for each $l\in L^z, z\in Z'_{Max}$:

\begin{equation}
\label{proxwra}
\begin{split}
&Op_{*}\{ Xdiv_{i_{\pi+1}}\dots Xdiv_{i_\mu}Repl\{
\tilde{C}^{l,i_1\dots i_\mu }_{g}
(\Omega_1,\dots,\Omega_p,\phi_1,\dots
,\phi_u)\nabla_{i_h}\phi_{u+1}\nabla_{i_2}\upsilon
\dots\nabla_{i_\pi}\upsilon\}\} \\&= Xdiv_{i_{\pi+1}}\dots
Xdiv_{i_\mu}C^{l,i_1i_{\pi+1}\dots i_\mu}_{g} (\Omega_1,\dots
,\Omega_p,\phi_{b+1},\dots,\phi_u)\nabla_{i_1}\phi_{u+1}+
\\&\Sum_{j\in J} a_j C^j_{g}(\Omega_1,\dots
,\Omega_p,\phi_{b+1},\dots,\phi_{u+1});
\end{split}
\end{equation}
(the tensor fields and complete contractions on the right hand
side have length $\sigma -b+u$. Here each $C^j_{g}$ has length
$\sigma -b+u+1$) and also has a
 factor $\nabla^{(s)}\phi_{u+1}, s>1$.

\par In the same way we derive that for each $h\in H$ (recall
 that $H_2=\emptyset$, hence the factor
$\nabla^{(K)}\omega_1$ in each $C^{h,i_1\dots i_{\mu+1}}_{g}$ has
$K>1$) and for each $t\in T'$:

\begin{equation}
\label{proxwra2}
\begin{split}
&Op_{*}\{ Xdiv_{i_{c+2}}\dots X div_{i_{\mu+1}}Repl\{
C^{h,i_1\dots i_{\mu+1} }_{g}
(\Omega_1,\dots,\Omega_p,\phi_1,\dots
,\phi_u)\nabla_{i_1}\phi_{u+1}\\&\nabla_{i_2}\upsilon
\dots\nabla_{i_\pi}\upsilon\}\}  =Xdiv_{i_{c+2}}\dots X
div_{i_{\mu+1}} C^{h,i_1i_{c+2}\dots i_{\mu+1}}_{g}
(\Omega_1,\dots ,\Omega_p,\phi_{b+1},\dots,\phi_u)
\\&+\Sum_{j\in
J} a_j C^j_{g}(\Omega_1,\dots ,\Omega_p,\phi_{b+1},\dots,\phi_u),
\end{split}
\end{equation}

\begin{equation}
\label{proxwra3}
\begin{split}
&Op_{*}\{ Xdiv_{i_{c+2}}\dots Xdiv_{i_\mu}Repl\{ C^{t,i_1\dots
i_\mu }_{g} (\Omega_1,\dots,\Omega_p,\phi_1,\dots
,\phi_u)\nabla_{i_1}\phi_{u+1}\nabla_{i_2}\upsilon
\dots\nabla_{i_\pi}\upsilon\}\} \\& =Xdiv_{i_{c+2}}\dots
Xdiv_{i_\mu}C^{t,i_1i_{c+2}\dots i_\mu}_{g} (\Omega_1,\dots
,\Omega_p,\phi_{b+1},\dots,\phi_u)\nabla_{i_1}\phi_{u+1}+
\\& \Sum_{j\in J} a_j C^j_{g}(\Omega_1,\dots
,\Omega_p,\phi_{b+1},\dots,\phi_{u+1}),
\end{split}
\end{equation}
where each $C^j_{g}$ has length $\sigma -b+u+1$ and also has a
 factor $\nabla^{(s)}\phi_{u+1}, s>1$.
\newline

{\it $Op_{*}$ produces a true equation:} Now, let us explain why
acting on (\ref{weakmayrh9}) by $Op_{*}$ will produce a true
equation: We break up (\ref{weakmayrh9}) (denote its left hand
side by $F$) into sublinear
 combinations $F^K$, according to the number $K$ of derivatives on
the factor $\nabla^{(K)}\omega_1$. Since (\ref{weakmayrh9}) holds
formally, it follows that $F^K=0$ formally (modulo longer
contractions).
 We then apply $Op_{*}$ to each equation $F^K=0$. This produces a true equation since
we may just repeat the permutations by which $F^K$ is made
formally zero to $Op_{*}\{F^K\}$. Adding over all equations
$Op_{*}\{F^K\}=0$, we derive our conclusion. $\Box$
\newline

\par For future reference, we formulate a corollary of Lemma \ref{vanderbi}:

\begin{corollary}
\label{vanderbicor} We consider a linear combination of
$\alpha$-tensor fields of weight $-n+\alpha$ and length $\sigma$
$$\Sum_{w\in W} a_w C^{w,i_1\dots i_\alpha}_{g}(\Omega_1,\dots \Omega_p,\phi_1,\dots ,\phi_b),$$
where each tensor field above has a given $b$-simple character
 $\vec{\kappa}_{*}$ and a given rank $\alpha$.
 We assume that for a given factor
 $T=S_{*}\nabla^{(\nu)}_{r_1\dots r_\nu} R_{ijkl}$ (for which the index ${}_i$
 is contracting against a given factor $\nabla\tilde{\phi}_k$)
  each tensor field indexed in $W$ has the feature
  that the factor $T$  has at least one of the indices
   ${}_{r_1},\dots ,{}_{r_\nu},{}_j$ {\it not}
   contracting against a factor $\nabla\phi$ (for this Lemma
    only we refer to this as the good property).
Assume an equation of the form:

\begin{equation}
\label{philon}
\begin{split}
&\Sum_{w\in W} a_w Xdiv_{i_1}\dots Xdiv_{i_\alpha} C^{w,i_1\dots
i_\alpha}_{g}(\Omega_1,\dots \Omega_p,\phi_1,\dots ,\phi_b)+
\\&\Sum_{h\in H} a_h Xdiv_{i_1}\dots Xdiv_{i_z}C^{h,i_1\dots i_z}_{g}
(\Omega_1,\dots \Omega_p,\phi_1,\dots ,\phi_b)
\\&= \Sum_{j\in J'} a_j
C^j_{g} (\Omega_1,\dots \Omega_p,\phi_1,\dots ,\phi_b),
\end{split}
\end{equation}
where each tensor field indexed in $H$ has a $b$-simple character
$\vec{\kappa}_{*}$ has rank $z\ge\beta$ (for some chosen $\beta$),
and does not satisfy the good property. Furthermore, we assume
that for these tensor fields of rank exactly $\beta$, if we
formally replace the expression $S_{*}\nabla^{(\nu)}_{r_1\dots
r_\nu}R_{ijkl}$ by $\nabla^{(\nu+1)}_{r_1\dots r_\nu j[k}\omega_1\nabla_{l]}\omega_2$,
 then the resulting tensor fields satisfy the
hypotheses of Lemma \ref{vanderbi} or \ref{vanderbi3}.
 Each complete contraction indexed in $J'$ is simply subsequent to $\vec{\kappa}_{*}$.

\par We then claim that:

\begin{equation}
\label{philon}
\begin{split}
&\Sum_{w\in W} a_w Xdiv_{i_1}\dots Xdiv_{i_\alpha} C^{w,i_1\dots
i_\alpha}_{g}(\Omega_1,\dots \Omega_p,\phi_1,\dots ,\phi_b)+
\\&\Sum_{h\in H'} a_h Xdiv_{i_1}\dots Xdiv_{i_z}C^{h,i_1\dots i_z}_{g}
(\Omega_1,\dots \Omega_p,\phi_1,\dots ,\phi_b)
\\&= \Sum_{j\in J} a_j
C^j_{g} (\Omega_1,\dots \Omega_p,\phi_1,\dots ,\phi_b),
\end{split}
\end{equation}
where each tensor field indexed in $H'$ are as the tensor field
indexed in $H$ above, only they now satisfy the good property.
Each complete contraction indexed in $J$ is simply subsequent to
$\vec{\kappa}_{*}$.
\end{corollary}

{\it Proof:} The proof just follows by reiterating the argument above:
We first use the operation
$Repl\{\dots\}$ as above, and then apply Lemma \ref{vanderbi}. We then
 use the operation $Op\{\dots\}$ as above and then the operation $Add$. $\Box$

\subsection{Derivation of Proposition \ref{giade} from Lemma
\ref{pskovb}.}

\par Firstly observe that we only have to show the above claim in case A,
since in case B the claims of Lemma \ref{pskovb} and Proposition \ref{giade} coincide. 

{\bf Proof in two steps:} We show that Lemma
\ref{pskovb} (in case A)  implies Proposition \ref{giade} in steps: Firstly,
using the conclusion of Lemma \ref{pskovb} we show that we
can derive a new equation:

\begin{equation}
\label{esmen'}
\begin{split}
&{\alpha\choose{2}}\Sum_{z\in Z'_{Max}}\Sum_{l\in L^z} a_l
\Sum_{r=0}^{k-1}Xdiv_{i_2}\dots Xdiv_{i_{*}}\dot{C}^{l,i_1\dots
\hat{i}_{r\alpha+1}\dots i_\mu,i_{*}}_{g}
(\Omega_1,\dots,\Omega_p,\phi_1,\dots ,\phi_u)
\\&\nabla_{i_{r\alpha+2}}\phi_{u+1}
+\Sum_{\nu\in N} a_\nu Xdiv_{i_1}\dots
Xdiv_{i_\mu}C^{\nu,i_1\dots
,i_\mu}_{g}(\Omega_1,\dots,\Omega_p,\phi_1,\dots
,\phi_u)\nabla_{i_1} \phi_{u+1}
\\& +\Sum_{t\in \tilde{T}_1} a_t Xdiv_{i_1}\dots Xdiv_{i_{z_t}}
C^{t, i_1\dots i_{z_t}}_{g}(\Omega_1,\dots ,\Omega_p,\phi_1,\dots
,\phi_{u+1})=
\\& \Sum_{j\in J}  a_j C^j_{g}(\Omega_1,\dots
,\Omega_p,\phi_1,\dots ,\phi_{u+1})=0;
\end{split}
\end{equation}
(notice the difference with (\ref{esmen}) is that we are not including
 $T_2,T_3,T_4$ and $T_1$ has been replaced by $\tilde{T}_1$).
Here the sublinear combination indexed in $\tilde{T}_1$ stands for a {\it generic}
 linear combination of the form $\sum_{t\in T_1}\dots$ 
described in the statement of Lemma \ref{pskovb}. This is step 1.

In step 2 we use (\ref{esmen'}) to derive Proposition \ref{giade}.
\newline

{\it Special cases:} There are two special cases which we will not
consider here, but treat in \cite{alexakis5}. The first  special case is
when $\sigma=3,p=3$,\footnote{Recall that $p$ stands for the number
of factors $\nabla^{(A)}\Omega_x$ in $\vec{\kappa}_{simp}$.} and
$n-2\mu-2u\le 2$. The second special case is when $\sigma=3$, $p=2,\sigma+2=1$ and $n=2\mu+2u$. 
 For the rest of this section, we will be
assuming that we do not fall under these special cases.
\newline

{\it Proof of Step 1:} Recall that $\sum_{t\in T_4}\dots$ appears 
only when the second critical factor is a simple factor of the from $\nabla^{(B)}\Omega_h$.  
In that case, we choose the factor $\nabla^{(B)}\Omega_x$ 
(referenced in the definition of the index set $L_\mu^*$)  to 
be the factor $\nabla^{(B)}\Omega_h$. (In other words we set $x=h$). 
In order to show (\ref{esmen'}), we recall the
hypothesis $L^{*}_\mu=\emptyset$ in Lemma \ref{pskovb}. In other
words, we are assuming that no
 tensor field $C^{l,i_1\dots i_\mu}_{g}$ in
(\ref{hypothese2}) has two free indices belonging to a
 factor $\nabla^{(2)}\Omega_h$. 

Now, refer to the conclusion of Lemma \ref{pskovb} (in case A). 
In view of the above remark, it follows that none of the $\mu$-tensor fields
$\dot{C}^{l,i_2\dots i_\mu,i_{*}}_{g} \nabla_{i_2}\phi_{u+1}$ or
$C^{\nu,i_1\dots i_\mu}_{g} \nabla_{i_1}\phi_{u+1}$ have an
expression $\nabla_i\Omega_h \nabla^i\phi_{u+1}$. Thus, all tensor
fields on the RHS of (\ref{esmen})
 with such an expression are indexed in $T_4$ (and have rank
$z_t\ge\mu$, by definition).

\par Now, we will firstly focus on the sublinear combination
$$\Sum_{t\in T_4}a_t Xdiv_{i_1}\dots Xdiv_{i_{z_t}}C^{t,i_1
\dots i_{z_t}}_{g}(\Omega_1,\dots ,\Omega_p,\phi_1,\dots ,
\phi_{u+1})$$ (recall $z_t\ge\mu$) {\it if it is nonzero}. (If it
is zero we move onto the next stage). We firstly seek to ``get
rid'' of the tensor fields indexed in $T_4$. More precisely, we
will show that we can write:

\begin{equation}
\label{pwpw}
\begin{split}
&\Sum_{t\in T_4}a_t Xdiv_{i_1}\dots Xdiv_{i_{z_t}}C^{t,i_1 \dots
i_{z_t}}_{g}(\Omega_1,\dots ,\Omega_p,\phi_1,\dots , \phi_{u+1})=
\\& \Sum_{t\in \tilde{T}_1}a_t Xdiv_{i_1}\dots Xdiv_{i_{z_t}}
C^{t,i_1 \dots i_{z_t}}_{g}(\Omega_1,\dots ,\Omega_p,\phi_1, \dots
,\phi_{u+1})+\\& \Sum_{j\in J} a_j C^j_{g} (\Omega_1,\dots
,\Omega_p,\phi_1,\dots, \phi_{u+1}),
\end{split}
\end{equation}
(with the same notational convention for the index
 set $\tilde{T}_1$ as above).
If we can prove this, we will be reduced to
 showing Step 1 under the assumption that $T_4=\emptyset$.

{\it Proof of \ref{pwpw}:} From (\ref{esmen}) we straightforwardly derive that:

\begin{equation}
\label{xaris} \Sum_{t\in T_4}a_t X_{*}div_{i_1}\dots
X_{*}div_{i_{z_t}} C^{t,i_1\dots i_{z_t}}_{g}(\Omega_1,\dots
,\Omega_p,\phi_1, \dots ,\phi_{u+1})=0,
\end{equation}
modulo complete contractions of length $\ge\sigma +u+2$. Here
 $X_{*}div_i$ stands for the sublinear combination in $Xdiv_i$
  where $\nabla_i$ is additionally not allowed to hit the
 expression $\nabla_k\Omega_h\nabla^k\phi_{u+1}$.

\par We observe that (\ref{pwpw}) follows
 from Lemma \ref{petermichel} if $\sigma>3$,\footnote{By the 
definition of $\sum_{t\in T_4}\dots$ in the statement of Proposition \ref{giade}, the assumptions of Lemma \ref{petermichel} are fulfilled.} and from
  Lemma \ref{petermichel3} if $\sigma=3$.\footnote{Notice that since
 we are assuming that (\ref{assumpcion})  does not fall under the special case
(described in the beginning of this subsection)
the requirements of Lemma \ref{petermichel3} are fulfilled.} $\Box$
\newline

 Thus, we may now prove our claim under the additionnal assumption that
$T_4=\emptyset$. Next, we want to ``get rid'' of the sublinear
combination:

$$\Sum_{t\in T_3}a_t Xdiv_{i_1}\dots Xdiv_{i_{z_t}}C^{t,i_1
\dots i_{z_t}}_{g}(\Omega_1,\dots ,\Omega_p,\phi_1,\dots ,
\phi_{u+1})$$ in (\ref{esmen}).

\par In particular, we will show that we can write:

\begin{equation}
\label{pwpw2}
\begin{split}
&\Sum_{t\in T_3}a_t Xdiv_{i_1}\dots Xdiv_{i_{z_t}}C^{t,i_1 \dots
i_{z_t}}_{g}(\Omega_1,\dots ,\Omega_p,\phi_1,\dots , \phi_{u+1})=
\\&\Sum_{t\in \tilde{T}_1}
a_t Xdiv_{i_1}\dots Xdiv_{i_{z_t}} C^{t,i_1 \dots
i_{z_t}}_{g}(\Omega_1,\dots ,\Omega_p,\phi_1, \dots
,\phi_{u+1})+\\& \Sum_{j\in J} a_j C^j_{g} (\Omega_1,\dots
,\Omega_p,\phi_1,\dots, \phi_{u+1});
\end{split}
\end{equation}
($\sum_{t\in \tilde{T}_1}\dots$ stands for a generic linear
combination as described in the conclusion of (\ref{esmen'})).

Thus, if we can show the above, we may additionaly assume that
$T_3=\emptyset$, in addition to our
assumption that $T_4=\emptyset$.
\newline

{\it Proof of (\ref{pwpw2}) when $\sigma>3$:} Break up $T_3$ into
subsets $\{T_3^h\}_{h=1,\dots ,p}$ according to the factor $\nabla\Omega_h$ that is
differentiated only once. We will then show (\ref{pwpw2})
 for each of the index sets $T_3^h$ separately.

\par To show this, we pick out the sublinear combination on (\ref{esmen}) with
a factor $\nabla\Omega_h$ (differentiated
 only once). This sublinear combination must vanish separately,
hence we derive an equation:

\begin{equation}
\label{mprei}
\begin{split}
& \Sum_{t\in T^h_3}a_t X_{*}div_{i_1}\dots
X_{*}div_{i_{z_t}}C^{t,i_1\dots i_{z_t}}_{g} (\Omega_1,\dots
,\Omega_p,\phi_1,\dots ,\phi_{u+1})
\\&+\sum_{j\in J} a_j C^j_g(\Omega_1,\dots
,\Omega_p,\phi_1,\dots ,\phi_{u+1})=0,
\end{split}
\end{equation}
modulo complete contractions of length $\ge\sigma +u+2$;
here as usual $X_{*}div_i$ stands for the subinear combination
 in $Xdiv_i$ where $\nabla_i$ is not allowed to hit the one
 factor $\nabla\Omega_h$.

Then, we see that (\ref{pwpw2}) follows from (\ref{mprei}) by
 applying  Corollary \ref{obotecor} above
(since $\mu\ge 4$ there are at least 2 derivative free indices 
for all maximal $\mu$-tensor fields in our Lemma assumption; therefore there exist 
at least two derivative free indices for each tensor field indexed in 
$T_3$, by weight considerations hence the requirements 
of Corollary \ref{corgiade} are fulfilled). $\Box$
\newline

{\it Proof of (\ref{pwpw2}) when $\sigma=3$:} We apply the
technique of the proof of Lemma \ref{petermichel3} (``manually''
constructing $Xdiv$'s) to write out:

\begin{equation}
\label{bravado}
\begin{split}
 &\sum_{t\in T_3} a_t Xdiv_{i_1}\dots Xdiv_{z_t}
a_t C^{t,i_1\dots i_{z_t}}_g(\Omega_1,\dots ,\Omega_p,\phi_1,\dots
,\phi_{u+1})=
\\&(Const)_{1,*} Xdiv_{i_1}\dots Xdiv_{i_{A+1}} C^{*1,i_1\dots i_{A+1}}_g
(\Omega_1,\dots ,\Omega_p,\phi_1,\dots ,\phi_{u+1})+
\\&(Const)_{2,*} Xdiv_{i_1}\dots Xdiv_{i_{A}} C^{*2,i_1\dots i_{A}}_g
(\Omega_1,\dots ,\Omega_p,\phi_1,\dots ,\phi_{u+1})+
\\&\sum_{t\in \tilde{T}_1} a_t Xdiv_{i_1}\dots Xdiv_{z_t}
a_t C^{t,i_1\dots i_{z_t}}_g(\Omega_1,\dots ,\Omega_p,\phi_1,\dots
,\phi_{u+1})+
\\&\sum_{j\in J} a_j C^j_g(\Omega_1,\dots ,\Omega_p,\phi_1,\dots
,\phi_{u+1}),
\end{split}
\end{equation}
where the tensor fields $C^{*1,i_1\dots i_{A+1}}_g$,
$C^{*2,i_1\dots i_A}_g$ are zero unless $p=3$ or $\sigma_1=2$ or
$\sigma_2=2$ or $\sigma_1=\sigma_2=1$ (in the last case there will
only be one tensor field $C^{*1,i_1\dots i_{A+1}}_g$ in the
above). In those cases, they stand for the following tensor
fields:

$$pcontr(\nabla^{(X)}_{i_1\dots i_a u_1\dots u_t}\Omega_1\otimes
\nabla^{(B)}_{j_1\dots j_by_1\dots
y_r}\Omega_2\otimes\nabla^{u_1}\phi_1\otimes\dots
\nabla^{y_r}\phi_{u+1}\otimes\nabla_{i_{A+1}}\Omega_3)$$

$$pcontr(\nabla^s\nabla^{(X)}_{i_1\dots i_a u_1\dots u_t}\Omega_1\otimes
\nabla^{(B)}_{j_1\dots j_by_1\dots
y_r}\Omega_2\otimes\nabla^{u_1}\phi_1\otimes\dots
\nabla^{y_r}\phi_{u+1}\otimes\nabla_{s}\Omega_3)$$ (here if $r\ge
2$ then $b=0$; if $r\le 1$ then $y=2-r$).

$$pcontr(\nabla^{(X)}_{i_1\dots i_a u_1\dots u_t}R_{i_{a+1}ji_{a+2}l}\otimes
\nabla^{(r)}_{y_1\dots
y_r}{{{R_{i_{a+3}}}^j}_{i_{a+4}}}^l\otimes\nabla^{u_1}\phi_1\dots
\nabla^{y_r}\phi_{u+1}\otimes\nabla_{i_{A+1}}\Omega_1),$$

$$pcontr(\nabla^s\nabla^{(X)}_{i_1\dots i_a u_1\dots u_t}R_{i_{a+1}ji_{a+2}l}\otimes
\nabla^{(r)}_{y_1\dots
y_r}{{{R_{i_{a+3}}}^j}_{i_{a+4}}}^l\otimes\nabla^{u_1}\phi_1\dots
\nabla^{y_r}\phi_{u+1}\otimes\nabla_{s}\Omega_1).$$ (In fact if
$t+r=0$ then there will be no $C^{*1,i_1\dots i_{A+1}}_g$ in
(\ref{bravado})).

\begin{equation}
\begin{split} 
&pcontr(S_{*}\nabla^{(X)}_{i_1\dots i_a u_1\dots u_t}R_{ii_{a+1}i_{a+2}l}\otimes
S_{*}\nabla^{(r)}_{y_1\dots
y_r}{{R_{i'i_{a+3}}}_{i_{a+4}}}^l\otimes\nabla^i\tilde{\phi}_1\otimes\nabla^{i'}\tilde{\phi}_2
\\&\otimes
\nabla^{u_1}\phi_3\otimes\dots
\nabla^{y_r}\phi_{u+1}\otimes\nabla_{i_{A+1}}\Omega_1)
\end{split}
\end{equation}

\begin{equation}
\begin{split} 
&pcontr(\nabla^sS_{*}\nabla^{(X)}_{i_1\dots i_a u_1\dots u_t}R_{ii_{a+1}i_{a+2}l}\otimes
S_{*}\nabla^{(r)}_{y_1\dots
y_r}{{R_{i'i_{a+3}}}_{i_{a+4}}}^l\otimes\nabla^i\tilde{\phi}_1\otimes\nabla^{i'}\tilde{\phi}_2
\\&\otimes
\nabla^{u_1}\phi_3\otimes\dots
\nabla^{y_r}\phi_{u+1}\otimes\nabla_{s}\Omega_1)
\end{split}
\end{equation}

\begin{equation}
\begin{split} 
&pcontr(S_{*}\nabla^{(X)}_{i_1\dots i_a u_1\dots u_t}R_{si_{a+1}i_{a+2}l}\otimes
\nabla^{(r)}_{y_1\dots
y_r}{{R_{i'i_{a+3}}}_{i_{a+4}}}^l\otimes\nabla^i\tilde{\phi}_1\otimes\nabla^{i'}\tilde{\phi}_2
\\&\otimes\nabla^{u_1}\phi_3\otimes\dots
\nabla^{y_r}\phi_{u+1}\otimes\nabla^{s}\Omega_1).
\end{split}
\end{equation}

\par As in the proof of Lemma \ref{petermichel3} we then derive
that $(Const)_{1,*}=(Const)_{2,*}=0$ in those cases; thus our
claim follows in this case also. $\Box$
\newline

\par Now, under the additional assumption that
 $T_3=T_4=\emptyset$,  we focus on the sublinear combination
$$\Sum_{t\in T_2} a_t Xdiv_{i_1}\dots Xdiv_{i_{z_t}}
C^{t,i_1,\dots ,i_{z_t}}_{g} (\Omega_1,\dots
,\Omega_p,\phi_1,\dots,\phi_{u+1})$$ in (\ref{esmen}). We will
show that we can write:

\begin{equation}
\label{adelfia}
\begin{split}
&\Sum_{t\in T_2} a_t Xdiv_{i_1}\dots Xdiv_{i_{z_t}} C^{t,i_1,\dots
,i_{z_t}}_{g} (\Omega_1,\dots ,\Omega_p,\phi_1,\dots, \phi_{u+1})=
\\&\Sum_{t\in \tilde{T}_1} a_t Xdiv_{i_1}\dots Xdiv_{i_{z_t}}
C^{t,i_1,\dots ,i_{z_t}}_{g} (\Omega_1,\dots
,\Omega_p,\phi_1,\dots, \phi_{u+1})+\\& \Sum_{j\in J} a_j C^j_{g}
(\Omega_1,\dots ,\Omega_p,\phi_1,\dots, \phi_{u+1}),
\end{split}
\end{equation}
where the notation is the same as in the statement of Lemma
\ref{pskovb}, and moreover $\Sum_{t\in \tilde{T}_1}\dots$ stands
 for a {\it generic} linear combination of the form described
 after (\ref{esmen'}). For each $t\in \tilde{T}_1$ we have
$z_t\ge \mu$.

\par We will show a more general statement, for
 future reference. 
\begin{lemma}
\label{addition}
 Consider a linear combination of
 acceptable tensor fields, $$\Sum_{l\in L_1} a_l
C^{l,i_1\dots i_{z_l}}_{g}(\Omega_1,\dots ,\Omega_p,\phi_1, \dots
,\phi_u)\nabla_{i_1}\phi_{u+1}$$ with a $u$-simple character
$\vec{\kappa}_{simp}$ ($\sigma\ge 3$) and with a $(u+1)$-simple
character $\vec{\kappa}^{+}_{simp}$, where in addition we are
assuming that if $\nabla_{i_1}\phi_{u+1}$ is contracting against a
factor $\nabla^{(m)}R_{ijkl}$ then it is contracting against a
derivative index, whereas if it is contracting against a factor
$S_{*}\nabla^{(\nu)}R_{ijkl}$ it must be contracting against
one of the indices ${}_{r_1},\dots ,{}_{r_\nu},{}_j$. (This is the defining property
of the $(u+1)$-simple character $\vec{\kappa}^{+}_{simp}$).

\par Consider another linear combination of
 acceptable tensor fields,  $$\Sum_{l\in L_2} a_l
C^{l,i_1\dots i_{z_l}}_{g}(\Omega_1,\dots ,\Omega_p,\phi_1, \dots
,\phi_u)\nabla_{i_1}\phi_{u+1}$$
 with a $u$-simple character
$\vec{\kappa}_{simp}$ and weak $(u+1)$-character equal to
$Weak(\vec{\kappa}^{+}_{simp})$, where in addition
$\nabla_{i_1}\phi_{u+1}$ is either contracting against an
 internal index in some factor $\nabla^{(m)} R_{ijkl}$ or an
index ${}_k$ or ${}_l$ in a factor $S_{*}\nabla^{(\nu)}R_{ijkl}$.
We moreover assume that each $l\in L_2$ we have $z_l\ge \gamma$,
for some number $\gamma$, and denote by $L^\gamma_2\subset L_2$
the index set of the tensor fields with order $\gamma$.

\par Assume that:

\begin{equation}
\label{terramomo}
\begin{split}
&Xdiv_{i_2}\dots Xdiv_{i_{z_l}}\Sum_{l\in L_1} a_l C^{l,i_1\dots
i_{z_l}}_{g}(\Omega_1,\dots ,\Omega_p,\phi_1, \dots
,\phi_u)\nabla_{i_1}\phi_{u+1}+
\\&Xdiv_{i_2}\dots Xdiv_{i_{z_l}}\Sum_{l\in L_2} a_l
C^{l,i_1\dots i_{z_l}}_{g}(\Omega_1,\dots ,\Omega_p,\phi_1, \dots
,\phi_u)\nabla_{i_1}\phi_{u+1}=
\\&\Sum_{j\in J} a_j C^{j,i_1}_{g}
(\Omega_1,\dots ,\Omega_p,\phi_1,
\dots ,\phi_u)\nabla_{i_1}\phi_{u+1},
\end{split}
\end{equation}
where each $C^{j,i_1}_{g}$ is $u$-subsequent to
$\vec{\kappa}_{simp}$. Furthermore assume that the above equation
falls under the inductive assumption of Proposition \ref{giade}
(with regard to the parameters weight, $\sigma$, $\Phi$, $p$).
Furthermore, we additionally assume that none
of the tensor fields $C^{l,i_1\dots i_{z_l}}_g$ 
of minimum rank in
(\ref{terramomo})\footnote{I.e. of rank $\gamma$.} are  ``forbidden''
in the sense of Proposition \ref{giade}.

\par Our first claim is then that there exists a linear combination
of $(\gamma+1)$-tensor fields,  $\Sum_{l\in L'_2} a_l
C^{l,i_1\dots i_{\gamma+1}}_{g} (\Omega_1,\dots ,\Omega_p,\phi_1,
\dots ,\phi_u)\nabla_{i_1}\phi_{u+1}$ with $u$-simple character
$\vec{\kappa}$ and weak $(u+1)$-character equal to
$Weak(\vec{\kappa}^{+}_{simp})$, where in addition
$\nabla_{i_1}\phi_{u+1}$ is either contracting against an
 internal index in some factor $\nabla^{(m)} R_{ijkl}$ or an
index ${}_k$ or ${}_l$ in a factor $S_{*}\nabla^{(\nu)}R_{ijkl}$,
so that:

\begin{equation}
\label{nerompogies}
\begin{split}
&\Sum_{l\in L^\gamma_2} a_l C^{l,i_1\dots
i_{\gamma}}_{g}(\Omega_1,\dots ,\Omega_p,\phi_1, \dots
,\phi_u)\nabla_{i_1}\phi_{u+1}\nabla_{i_2}\upsilon\dots
\nabla_{i_\gamma}\upsilon
\\& -Xdiv_{i_\gamma+1}\Sum_{l\in L'_2} a_l
C^{l,i_1\dots i_{\gamma+1}}_{g} (\Omega_1,\dots ,\Omega_p,\phi_1,
\dots ,\phi_u)\nabla_{i_1}\phi_{u+1}\nabla_{i_2}\upsilon\dots
\nabla_{i_\gamma}\upsilon
\\& +\Sum_{l\in \overline{L}_1} a_l C^{l,i_1\dots i_{\gamma}}_{g}(\Omega_1,\dots ,\Omega_p,\phi_1,
\dots ,\phi_u)\nabla_{i_1}\phi_{u+1}\nabla_{i_2}\upsilon\dots
\nabla_{i_\gamma}\upsilon+
\\& \Sum_{j\in J} a_j C^{j,i_1\dots i_\gamma}_{g}
(\Omega_1,\dots ,\Omega_p,\phi_1,\dots ,\phi_u)\nabla_{i_1}
\phi_{u+1}\nabla_{i_2}\upsilon\dots\nabla_{i_\gamma}\upsilon;
\end{split}
\end{equation}
here each $C^{j,i_1\dots i_\gamma}_{g}$ is $u$-subsequent to
$\vec{\kappa}_{simp}$. The tensor fields indexed in
$\overline{L}_1$ are like the ones indexed in $L_1$ in
(\ref{terramomo}), but in addition each $z_l\ge \gamma$.
\newline

\par Our second claim is that assuming (\ref{terramomo}) we can write:

\begin{equation}
\label{caputo}
\begin{split}
&Xdiv_{i_2}\dots Xdiv_{i_{z_l}}\Sum_{l\in L_2} a_l C^{l,i_1\dots
i_{z_l}}_{g}(\Omega_1,\dots ,\Omega_p,\phi_1, \dots
,\phi_u)\nabla_{i_1}\phi_{u+1}=
\\&\Sum_{l\in \overline{L}_1} a_lXdiv_{i_2}\dots Xdiv_{i_{z_l}}
C^{l,i_1\dots i_{z_l}}_{g}(\Omega_1,\dots ,\Omega_p,\phi_1, \dots
,\phi_u)\nabla_{i_1}\phi_{u+1}+
\\&\Sum_{j\in J} a_j C^{j,i_1}_{g}
(\Omega_1,\dots ,\Omega_p,\phi_1,
\dots ,\phi_u)\nabla_{i_1}\phi_{u+1},
\end{split}
\end{equation}
 where $C^{j,i_1}_{g}$ is $u$-subsequent to
$\vec{\kappa}_{simp}^{+}$.

\par Our third claim is that if $\gamma$ is the minimum rank among
all tensor fields in $L_1\bigcup L_2$ in our assumption and
$L^\gamma_1,L^\gamma_2$ their respective index sets, then there
exists a linear combination of $(\gamma+1)$-tensor fields,
\\$\Sum_{l\in L_3} a_l C^{l,i_1\dots i_{\gamma+1}}_{g}
(\Omega_1,\dots ,\Omega_p,\phi_1, \dots
,\phi_u)\nabla_{i_1}\phi_{u+1}$ with $u$-simple character
$\vec{\kappa}$ and weak $(u+1)$-character equal to
$Weak(\vec{\kappa}^{+}_{simp})$ so that:

\begin{equation}
\label{nerompogies2}
\begin{split}
&\Sum_{l\in L^\gamma_1\bigcup L^\gamma_2} a_l C^{l,i_1\dots
i_{\gamma}}_{g}(\Omega_1,\dots ,\Omega_p,\phi_1, \dots
,\phi_u)\nabla_{i_1}\phi_{u+1}\nabla_{i_2}\upsilon\dots
\nabla_{i_\gamma}\upsilon
\\& -Xdiv_{i_{\gamma+1}}\Sum_{l\in L^3} a_l
C^{l,i_1\dots i_{\gamma+1}}_{g} (\Omega_1,\dots ,\Omega_p,\phi_1,
\dots ,\phi_u)\nabla_{i_1}\phi_{u+1}\nabla_{i_2}\upsilon\dots
\nabla_{i_\gamma}\upsilon=
\\& \Sum_{j\in J} a_j C^{j,i_1\dots i_\gamma}_{g}
(\Omega_1,\dots ,\Omega_p,\phi_1,\dots ,\phi_u)\nabla_{i_1}
\phi_{u+1}\nabla_{i_2}\upsilon\dots\nabla_{i_\gamma}\upsilon;
\end{split}
\end{equation}
here each $C^{j,i_1\dots i_\gamma}_{g}$ is $u$-subsequent to
$\vec{\kappa}_{simp}$. 
\end{lemma}

We observe that if we can show the above, then our claim
(\ref{adelfia}) follows from the second step of this Lemma.
\newline

{\it Proof of Lemma \ref{addition}:}  We firstly remark that in proving 
Lemma \ref{addition} we will use Lemma \ref{appendix}, which 
is stated and proven in the Appendix of this paper. 
We also easily observe that the third claim above follows
from the first two. So we now prove the first two claims in that Lemma:
\newline

{\it Proof of the second claim of Lemma \ref{addition}:} We now
show that the
 second claim follows from the first one.

\par We will show this by induction. We assume that
$min_{l\in L_2}z_t=\gamma'\ge\gamma$. We denote the index set of
those tensor fields by $L_2^{\gamma'}\subset L_2$. Then, using the
first claim\footnote{Provided that there are no terms indexed in
$L_2^{\gamma'}$ which are forbidden.} and making the
 $\nabla\upsilon$s into $Xdiv$s, we derive that we can write:

\begin{equation}
\label{kuranishi}
\begin{split}
&\Sum_{l\in L_2^{\gamma'}} a_l Xdiv_{i_1}\dots Xdiv_{i_\gamma}
C^{t,i_1,\dots ,i_{z_l}}_{g} (\Omega_1,\dots
,\Omega_p,\phi_1,\dots,\phi_{u+1})=
\\&\Sum_{l\in L_2'} a_l Xdiv_{i_1}\dots Xdiv_{i_{z_t}}
C^{l,i_1,\dots ,i_{z_l}}_{g} (\Omega_1,\dots
,\Omega_p,\phi_1,\dots,\phi_{u+1})+
\\&\Sum_{l\in L_1} a_l Xdiv_{i_1}\dots Xdiv_{i_{z_l}}
C^{l,i_1,\dots ,i_\gamma}_{g} (\Omega_1,\dots
,\Omega_p,\phi_1,\dots,\phi_{u+1}) +\\& \Sum_{j\in J} a_j C^j_{g}
(\Omega_1,\dots ,\Omega_p,\phi_1,\dots,\phi_{u+1}),
\end{split}
\end{equation}
where the tensor fields indexed in $L_2'$ are of the exact same
form as the ones indexed in $L_2$ in (\ref{pskovb}), with the
additional property that $z_l\ge \gamma' +1$. We notice that since we are
 dealing with tensor fields of a given weight $-n$,
iteratively repeating this step we derive our second
 step. ({\it Note:} If at the last step we encounter a
 ``forbidden case'' then clearly  $\gamma'>\gamma$--we then apply Lemma \ref{appendix} 
below with $\Phi=1$). 
\newline

{\it Proof of the first claim of Lemma \ref{addition}:}
 The proof requires only our inductive assumption on
Corollary \ref{corgiade}. We have two cases to consider:
Firstly, when the factor $\nabla\phi_{u+1}$ is contracting
 (in $\vec{\kappa}^{+}_{simp}$) against an internal index (say ${}_i$ with no loss of
 generality) of a factor $\nabla^{(m)}R_{ijkl}$. Secondly,
when the factor $\nabla\phi_{u+1}$ is contracting
 against an index ${}_k$ of a factor $S_{*}\nabla^{(\nu)}R_{ijkl}$.

{\it Proof of first claim of Lemma \ref{addition} in the first
case:} In the first case, we define an operation $Cutsym$ that
acts on the tensor fields indexed in $L_2$ by replacing
 the expression $\nabla^{(m)}_{r_1\dots r_m}R_{ijkl}\nabla^{r_1}
\phi_{t_1}\dots \nabla^{r_a}\phi_{t_a}\nabla^i\phi_{u+1}$ by an
expression $S_{*}\nabla^{(m-a)}_{r_{a+1}\dots r_n}
R_{ijkl}\nabla^i\phi_{u+1}$. We observe that the tensor fields
that arise via this operation have a given simple character which
we will denote by $\vec{\kappa}_{cut}$. For each $l\in L_2$ we
denote by
$$C^{l,i_1\dots i_{z_l}}_{g}(\Omega_1,\dots ,\Omega_p,\phi_1,
\dots \hat{\phi}_{r_{a_1}},\dots ,\hat{\phi}_{r_{a_t}},
\dots ,\phi_u)\nabla_{i_1}\phi_{u+1}$$
the tensor field that we obtain by applying this operation.

\par We also define the operation $Cutsym$ to act
 on the tensor fields indexed in $L_1$ by replacing the
the expression $\nabla^{(m)}_{r_1\dots r_m}R_{ijkl}\nabla^{r_1}
\phi_{t_1}\dots \nabla^{r_a}\phi_{t_a}\nabla^{r_b} \phi_{u+1}$ by
a factor $\nabla^{(m-1)}_{r_{a+1}\dots r_m}R_{ijkl}
\nabla^{r_b}\phi_{u+1}$. Now, by applying the eraser to the
factors $\nabla^{r_1} \phi_{t_1}\dots \nabla^{r_a}\phi_{t_a}$ and
$S_{*}$-symmetrizing,
 we may apply $CutSym$ to (\ref{terramomo}) and derive an
equation:

\begin{equation}
\label{terramomo2}
\begin{split}
& Xdiv_{i_2}\dots Xdiv_{i_{z_l}}\Sum_{l\in L_2} a_l C^{l,i_1\dots
i_{z_l}}_{g}(\Omega_1,\dots ,\Omega_p,\phi_1, \dots
\hat{\phi}_{r_{a_1}},\dots ,\hat{\phi}_{r_{a_t}}, \dots
,\phi_u)\\&\nabla_{i_1}\phi_{u+1}
=\Sum_{j\in J} a_j C^{j,i_1}_{g}
(\Omega_1,\dots ,\Omega_p,\phi_1,
\dots \hat{\phi}_{r_{a_1}},\dots ,\hat{\phi}_{r_{a_t}}, \dots,\phi_u)\nabla_{i_1}\phi_{u+1},
\end{split}
\end{equation}
where each $C^j_{g}$ is $(u-a)$-subsequent to
$\vec{\kappa}_{cut}$. We may then apply Corollary \ref{corgiade} to the 
above.\footnote{Corollary \ref{corgiade} may be applied by virtue of our assumptions 
on various terms  in our Lemma assumption not being 
``forbidden''. This ensures that the terms of minimum rank 
in (\ref{terramomo2}) are not ``forbidden'' in the 
sense of Corollary \ref{corgiade}.} 
This follows since either the weight in (\ref{terramomo2})
 is $-n', n'<n$, (this occurs when we erase factors
 $\nabla\phi_t$ upon performing the operation $CutSym$),
  or the weight is $-n$ and there are $u+1$ factors $\nabla\phi_h$ in
   (\ref{terramomo2}). Thus, our inductive assumption of
    Corollary \ref{corgiade} applies to (\ref{terramomo2}).

\par Thus, by direct application of
 Corollary \ref{corgiade} (which we are now inductively
assuming because either the weight is $>-n$ or there are $(u+1)$
factors $\nabla\phi$) to (\ref{terramomo2}) we derive that there
is
 an acceptable linear combination of $(\gamma+1)$-tensor
fields with a simple character $\vec{\kappa}_{cut}$,
 say $$\Sum_{x\in X} a_x C^{x,i_1\dots i_{\gamma+1}}_{g}
(\Omega_1,\dots ,\Omega_p,\phi_1, \dots \hat{\phi}_{r_{a_1}},\dots
,\hat{\phi}_{r_{a_t}}, \dots ,\phi_u)\nabla_{i_1}\phi_{u+1},$$ so
that:

\begin{equation}
\label{amoremio}
\begin{split}
&\Sum_{l\in L^\gamma_2} a_l C^{l,i_1\dots i_{\gamma}}_{g}
(\Omega_1,\dots ,\Omega_p,\phi_1, \dots \hat{\phi}_{r_{a_1}},\dots
,\hat{\phi}_{r_{a_t}}, \dots
,\phi_u)\nabla_{i_1}\phi_{u+1}\nabla_{i_2}\upsilon\dots
\nabla_{i_\gamma}\upsilon-
\\&\Sum_{x\in X} a_x C^{x,i_1\dots i_{\gamma+1}}_{g}
(\Omega_1,\dots ,\Omega_p,\phi_1,
\dots \hat{\phi}_{r_{a_1}},\dots ,\hat{\phi}_{r_{a_t}},
\dots ,\phi_u)\nabla_{i_1}\phi_{u+1}\nabla_{i_2}\upsilon\dots
\nabla_{i_\gamma}\upsilon
\\&=\Sum_{j\in J} a_j C^{j,i_1\dots i_{\gamma}}_{g}
(\Omega_1,\dots ,\Omega_p,\phi_1,
\dots \hat{\phi}_{r_{a_1}},\dots ,\hat{\phi}_{r_{a_t}},
\dots ,\phi_u)\nabla_{i_1}\phi_{u+1}\nabla_{i_2}\upsilon\dots
\nabla_{i_\gamma}\upsilon,
\end{split}
\end{equation}
where each tensor field $C^j_{g}$ is subsequent to
$\vec{\kappa}_{cut}$.

\par Now, we define an operation $Add$ that acts on the tensor
 fields above by replacing the expression
$S_{*}\nabla^{(m-a)}_{r_{a+1}\dots r_n}R_{ijkl}\nabla^i
\phi_{u+1}$ by an expression
$$\nabla_{r_1\dots r_a}
S_{*}\nabla^{(m-a)}_{r_{a+1}\dots r_m}R_{ijkl}\nabla^{r_1}
\phi_{t_1}\dots \nabla^{r_a}\phi_{t_a}\nabla^i \phi_{u+1}.$$ In
case $\nabla\phi_{u+1}$ is contracting against some derivative
index in some $\nabla^{(m)}R_{ijkl}$, it adds on the factor
$\nabla^{(m)}R_{ijkl}$ against which $\nabla\phi_{u+1}$ is
contracting $a$ derivative indices and contracts them against
factors $\nabla\phi_{a_1},\dots ,\phi_{a_t}$. By applying the
operation $Add$ to (\ref{amoremio}) we derive our desired equation
(\ref{nerompogies}). $\Box$
\newline

\par The second case is treated in a similar fashion. We now define a
formal  operation
 $CutY$ as follows: $CutY$
acts on the tensor fields indexed in $L_2$ by replacing  the
expression $S_{*}\nabla^{(\nu)}_{r_1\dots
r_\nu}R_{ir_{\nu+1}kl}\nabla^{r_1} \phi'_{t_1}\dots
\nabla^{r_a}\phi'_{t_a}\nabla^i \tilde{\phi}_{t_{a+1}}\nabla^k\phi_{u+1}$
by an expression $\nabla^{(\nu-a+2)}_{r_{a}\dots r_{n+1}l}Y\nabla^{r_a}\phi'_{t_a}$, (if
there is at least
 one factor $\nabla\phi'$ contracting against our factor
$S_{*}\nabla^{(\nu)}R_{ijkl}$; if there is no such factor we
replace $S_{*}\nabla^{(\nu)}_{r_1\dots
r_\nu}R_{ir_{\nu+1}kl}\nabla^i \tilde{\phi}_{t_{a+1}}\nabla^k\phi_{u+1}$
by $\nabla^{(\nu+2)}_{r_1\dots r_{\nu+1}l}Y$.
 We will denote the tensor field thus obtained by
\\ $C^{l,i_1\dots
i_{z_l}}_{g}(\Omega_1,\dots ,\Omega_p,Y, \phi_1,\dots
\hat{\phi}_{t_1},\dots ,\hat{\phi}_{t_{a+1}}, \dots ,\phi_u)$. (Observe
that it is acceptable if we treat the function $Y$ as
 a function $\Omega_{p+1}$.
 We observe that all the tensor fields that arise thus
have a given simple character which we will denote by
$\vec{\kappa}_{cut}$). We also define  $CutY$ to act
 on the tensor fields indexed in $L_1$ by replacing them by
 zero. Finally, it follows easily that the operation  $CutY$ either annihilates 
a given complete contraction $C^j_{g}$, or replaces it by a complete
contraction that is subsequent to $\vec{\kappa}_{cut}$.

Now, by virtue of Lemma \ref{technical1} and the ``Eraser'' 
(defined in the Appendix of \cite{alexakis1}), we see that
applying $CutY$ to
 (\ref{terramomo}) produces a true equation which can be written as:

\begin{equation}
\label{terramomo3}
\begin{split}
& Xdiv_{i_2}\dots Xdiv_{i_{z_l}}\Sum_{l\in L_2} a_l C^{l,i_1\dots
i_{z_l}}_{g}(\Omega_1,\dots ,\Omega_p,Y, \phi_1,\dots
\hat{\phi}_{t_1},\dots , \hat{\phi}_{t_{a+1}},\dots ,\phi_u)=
\\&\Sum_{k\in K} a_k C^{k,i_1}_{g}
(\Omega_1,\dots ,\Omega_p,Y, \phi_1, \dots \hat{\phi}_{t_1},\dots
,\hat{\phi}_{t_{a+1}},
 \dots,\phi_u),
\end{split}
\end{equation}
where each $C^k_{g}$ is simply subsequent to
$\vec{\kappa}_{cut}$. Thus, by direct application of Corollary
\ref{corgiade} to (\ref{terramomo2}),\footnote{The observation
in the previous footnote still applies--by virtue of the
assumptions imposed in our Lemma hypothesis, (\ref{terramomo2}) does not fall under a ``forbidden case''.}
 we derive that there is
 an acceptable linear combination of $(\gamma+1)$-tensor
fields with a simple character $\vec{\kappa}_{cut}$,
 \\ $\Sum_{x\in X} a_x C^{x,i_1\dots i_{\gamma+1}}_{g}
(\Omega_1,\dots ,\Omega_p,Y\phi_1, \dots \hat{\phi}_{t_1},\dots
,\hat{\phi}_{t_a}, \dots ,\phi_u)$, so that:

\begin{equation}
\label{amoremio2}
\begin{split}
&\Sum_{l\in L^\gamma_2} a_l C^{l,i_1\dots i_{\gamma}}_{g}
(\Omega_1,\dots ,\Omega_p,Y,\phi_1, \dots \hat{\phi}_{t_1},\dots
,\hat{\phi}_{t_a}, \dots ,\phi_u)\nabla_{i_2}\upsilon\dots
\nabla_{i_\gamma}\upsilon-
\\&\Sum_{x\in X} a_x Xdiv_{i_{\gamma+1}}C^{x,i_1\dots i_{\gamma+1}}_{g}
(\Omega_1,\dots ,\Omega_p,Y,\phi_1, \dots \hat{\phi}_{t_1},\dots
,\hat{\phi}_{t_a}, \dots ,\phi_u)\nabla_{i_2}\upsilon\dots
\nabla_{i_\gamma}\upsilon
\\&=\Sum_{k\in K} a_k C^{k,i_1\dots i_{\gamma}}_{g}
(\Omega_1,\dots ,\Omega_p,Y,\phi_1, \dots \hat{\phi}_{t_1},\dots
,\hat{\phi}_{t_a}, \dots ,\phi_u)\nabla_{i_2}\upsilon\dots
\nabla_{i_\gamma}\upsilon,
\end{split}
\end{equation}
where each tensor field $C^k_{g}$ is subsequent to
$\vec{\kappa}_{cut}$. (Note that the LHS in  (\ref{terramomo3})
has weight $>-n$, hence Corollary \ref{corgiade} applies, thanks
to our inductive assumption).

\par Now, we define a formal  operation $UnY$ that acts on the tensor
 fields above by replacing the expression
$\nabla^{(B)}_{t_1\dots t_B}Y$, $B\ge 2$ by an expression
\\ $\nabla^{(B-2+a)}_{r_1\dots r_at_1\dots t_{B-2+a}}
R_{ijkl}\nabla^{r_1} \phi_{t_1}\dots \nabla^{r_a}\phi_{t_a}
\nabla^i\phi_{t_{a+1}}\nabla^k\phi_{u+1}$. By applying the
operation $UnY$ to (\ref{amoremio2}) (and repeating the
permutations
 by which (\ref{amoremio2}) is made formally zero, modulo introducing
 correction terms by virtue of the Bianchi identities--see (\ref{koichi1}),
 (\ref{koichi2}),(\ref{koichi3}))  we derive our desired
equation (\ref{nerompogies}). $\Box$

\par This completes the proof of step 1
(in the derivation of Proposition \ref{giade}
(in case III) from Lemma \ref{pskovb}. $\Box$
\newline

{\it Proof of step 2 (in the derivation of Proposition \ref{giade}
(in case III)) from Lemma \ref{pskovb}:} We consider
(\ref{esmen'}), (where all the tensor fields are now acceptable,
by definition).
Recall that the $(u+1,\mu-1)$-refined double characters
 that correspond to the index sets $L^z,z\in Z'_{Max}$ in (\ref{esmen'}) (we have denoted them by
$\vec{L}^{z,\sharp}$) are the maximal ones. Now, we can apply our inductive
 assumption of Proposition \ref{giade} to (\ref{esmen'}):\footnote{The inductive
  assumption of Proposition \ref{giade} applies here since
  we have weight $-n$ but have an extra factor $\nabla\phi_{u+1}$.
  Observe that the $(\mu-1)$-tensor fields in that equation
have no special free indices, hence there is no danger of
``forbidden cases''.}

 We derive that for each $z\in Z'_{Max}$,
there is a linear combination of acceptable
 $\mu$-tensor fields (which satisfy the extra restriction if it is applicable),
$$\Sum_{p\in P'} a_p
C^{p,i_1\dots i_\mu}_{g}
(\Omega_1,\dots,\Omega_p,\phi_1,\dots,\phi_u,\phi_{u+1})$$ with a
$(u+1,\mu-1)$-refined double
 character $\vec{L}^{z,\sharp}$, so that for any
$z\in Z'_{Max}$:

\begin{equation}
\label{baldwin}
\begin{split}
&{\alpha\choose{2}}\Sum_{l\in L^z} a_l
\Sum_{r=1}^{k-1} \dot{C}^{l,i_1\dots
\hat{i_{r\alpha+1}}\dots i_\mu,i_{*}}_{g}
(\Omega_1,\dots,\Omega_p,\phi_1,\dots ,\phi_u)\nabla_{i_1}
\phi_{u+1}\nabla_{i_2}\upsilon\dots \nabla_{i_{*}}\upsilon+
\\& \Sum_{p\in P'} a_p Xdiv_{i_\mu}
C^{p,i_1\dots i_\mu}_{g} (\Omega_1,\dots,\Omega_p,\phi_1,
\dots,\phi_u,\phi_{u+1}) \nabla_{i_1}\upsilon\dots
\nabla_{i_{\mu-1}}\upsilon=
\\&\Sum_{k\in K} a_k C^k_{g} (\Omega_1,\dots,\Omega_p,
\phi_1,\dots,\phi_u,\phi_{u+1},\upsilon^{\mu-1}),
\end{split}
\end{equation}
modulo complete contractions of length $\ge\sigma +u+\mu +1$. Here
each $C^k$ is (simply or doubly) subsequent to each
$\vec{L}^{z,\sharp}$,
$z\in Z'_{Max}$.

We then define a formal operation $Op$ that acts
 on the tensor fields in the above by performing two actions: Firstly, we
 pick out a derivative index in the critical
 factor (the {\it unique} factor that is contracting against
 the most factors $\nabla\upsilon$) that contracts
 against a factor $\nabla\upsilon$ and erase it. Secondly, we then add a derivative
 index $\nabla_{i_{+}}$ onto the A-crucial factor
and contract it against the above factor $\nabla\upsilon$.

\par Let us observe that this operation is well-defined,
 and then see that it produces a true equation:
The only thing that could make this operation not well-defined is
if no factor $\nabla\upsilon$ is contracting against a derivative
index in the critical factor (this can only be the case for
factors $S_{*}\nabla^{(\nu)}R_{ijkl}$ with $\nu=0$). But that
cannot happen: Recall the critical factor must start out with at
least two free indices (none of them special), and then we add
another derivative index onto it. Thus in all complete
contractions in (\ref{baldwin}) there are at least three factors
$\nabla\upsilon$ contracting against (non-special) indices in the
critical factor. Thus, our operation $Op$ is well-defined.
 By the same reasoning, observe that our operation $Op$ produces acceptable tensor fields.

\par  We then set
$\phi_{u+1}=\upsilon$. We have observed that $Op$ is well defined,
and we see that after we set $\phi_{u+1}=\upsilon$, we
will have that for each $l\in L^z$, $z\in Z'_{Max}$:

\begin{equation}
\label{baldwin2}
\begin{split}
&Op[\dot{C}^{l,i_1\dots \hat{i_{r\alpha+1}}\dots i_\mu,i_{*}}_{g}
(\Omega_1,\dots,\Omega_p,\phi_1,\dots ,\phi_u)\nabla_{i_1}
\phi_{u+1}\nabla_{i_2}\upsilon\dots \nabla_{i_{*}}\upsilon]=
\\& C^{l,i_1\dots i_\mu}_{g}
(\Omega_1,\dots,\Omega_p,\phi_1,\dots ,\phi_u)
\nabla_{i_1}\upsilon\dots \nabla_{i_\mu}\upsilon.
\end{split}
\end{equation}
Hence, applying $Op$ to (\ref{baldwin}) (which produces a true equation
 since (\ref{baldwin}) holds formally) gives us step 2 and thus we derive
the claim of Proposition \ref{giade} in case $III$ from Lemma \ref{pskovb}. $\Box$

\section{Appendix.}
\subsection{A  weak substitute for Proposition 
\ref{giade} in the ``forbidden cases''.}
\label{weak.sub}

We present a ``substitute'' of sorts of the Proposition \ref{giade} 
in the  ``forbidden cases''. This ``substitute'' (Lemma \ref{appendix} below) 
will rely on a generalized version of the 
Lemma \ref{addition}, Lemma \ref{additiongen}
which is stated below but proven in \cite{alexakis5}. 
The generalized version asserts that the claim of Lemma \ref{addition} remains true, 
for the general case where rather than one ``additional''
 factor $\nabla\phi_{u+1}$ we have $\beta\ge 3$ ``additional'' 
factors $\nabla\phi_{u+1},\dots,\nabla\phi_{u+\beta}$. Moreover, 
in that case there are no ``forbidden cases''.

\begin{lemma}
\label{additiongen}
Let $\sum_{l\in L_1} a_l C^{l,i_1\dots i_\mu,i_{\mu+1}\dots i_{\mu+\beta}}_g
(\Omega_1,\dots,\Omega_p,\phi_1,\dots,\phi_u)$,

$\sum_{l\in L_2} a_l C^{l,i_1\dots i_{b_l},i_{b_l+1}\dots i_{b_l+\beta}}_g
(\Omega_1,\dots,\Omega_p,\phi_1,\dots,\phi_u)$ 
 stand for two linear combinations of acceptable 
tensor fields in the form (\ref{form2}), with a $u$-simple 
character $\vec{\kappa}_{simp}$. We assume that the terms 
indexed in $L_1$ have rank $\mu+\beta$, while 
the ones indexed in $L_2$ have rank greater than $\mu+\beta$. 

Assume an equation:

\begin{equation}
\begin{split}
\label{paragon} & \sum_{l\in L_1} a_l Xdiv_{i_1}\dots
Xdiv_{i_\mu} C^{l,i_1\dots
i_{\mu+\beta}}_g(\Omega_1,\dots,\Omega_p,\phi_1,\dots,\phi_u)
\nabla_{i_1}\phi_{u+1}\dots \nabla_{i_\beta}\phi_{u+\beta}
\\&+\sum_{l\in L_2} a_l Xdiv_{i_1}\dots
Xdiv_{i_{b_l}} C^{l,i_1\dots
i_{b_l+\beta}}_g(\Omega_1,\dots,\Omega_p,\phi_1,\dots,\phi_u)
\nabla_{i_1}\phi_{u+1}\dots \nabla_{i_\beta}\phi_{u+\beta}
\\&+\sum_{j\in J} a_j C^{j}_g(\Omega_1,\dots,\Omega_p,\phi_1,\dots,\phi_{u+\beta})=0,
\end{split}
\end{equation}
modulo terms of length $\ge\sigma+u+\beta+1$.
Furthermore, we assume that the above equation 
falls under the inductive assumption of Proposition 2.1 in \cite{alexakis4} 
(with regard to the parameters weight, $\sigma,\Phi,p$). 
 We are not excluding any  ``forbidden cases''.

We claim that there exists a linear combination of $(\mu+\beta+1)$-tensor fields 
in the form (\ref{form2}) with $u$-simple character $\vec{\kappa}_{simp}$ and 
length $\sigma+u$ (indexed in $H$ below) such that:

\begin{equation}
\begin{split}
\label{paragon.conc} & \sum_{l\in L_1} a_l C^{l,i_1\dots
i_{\mu+\beta}}_g(\Omega_1,\dots,\Omega_p,\phi_1,\dots,\phi_u)
\nabla_{i_1}\phi_{u+1}\dots \nabla_{i_\beta}
\phi_{u+\beta}\nabla_{i_1}\upsilon\dots\nabla_{i_\mu}\upsilon
\\&+\sum_{h\in H} a_h Xdiv_{i_{\mu+\beta+1}} C^{l,i_1\dots
i_{\mu+\beta+1}}_g(\Omega_1,\dots,\Omega_p,\phi_1,\dots,\phi_u)
\nabla_{i_1}\phi_{u+1}\dots \nabla_{i_\beta}
\phi_{u+\beta}\\&\nabla_{i_1}\upsilon\dots\nabla_{i_\mu}\upsilon+
\sum_{j\in J} a_j C^j_g(\Omega_1,\dots,
\Omega_p,\phi_1,\dots,\phi_{u+\beta},\upsilon^\mu)=0,
\end{split}
\end{equation}
modulo terms of length $\ge\sigma+u+\beta+\mu+1$. The terms indexed in 
$J$ here are $u$-simply subsequent to $\vec{\kappa}_{simp}$.
\end{lemma}

A note and a notational convention before we state our Lemma: 
We observe that if some of the $\mu$-tensor fields of maximal refined double character 
in (\ref{hypothese2}) are ``forbidden'', then necessarily {\it all}
 tensor fields in (\ref{hypothese2}) must have rank $\mu$ 
(in other words $L_{>\mu}=\emptyset$).
This follows from weight considerations.
 {\it Also} all the $\mu$-tensor fields must 
have each of the $\mu$ free indices belonging to a 
different factor.  This follows from the definition
 of maximal refined double character.

We introduce the notational convention needed for our Lemma. For each tensor field
 $C^{l,i_1\dots i_\mu}_g$ appearing in (\ref{hypothese2}) we will consider 
its product with an auxilliary function $\Phi$, $C^{l,i_1\dots i_\mu}_g\cdot\Phi$. 
In that contect $Xdiv_{i_a}[C^{l,i_1\dots i_\mu}_g\cdot\Phi]$ 
will stand for the sublinear combination in where $\nabla^{i_a}$ 
is not allowed to hit the factor to which ${}_{i_a}$ belongs 
(but it {\it is} allowed to hit the function $\Phi$).

\begin{lemma}
\label{appendix}
Assume equation (\ref{hypothese2}),  under the additional 
assumption that some of the tensor fields of 
maximal refined double character in $L_\mu$ are ``forbidden'', 
in the sense of definition \ref{forbidden}. Denote by $\vec{\kappa}_{simp}$ 
the $u$-simple character of the tenosr field in (\ref{hypothese2}). 

We then claim that there is then a linear combination
 of acceptable $\mu$-tensor fields 
with a $u$-simple character $\vec{\kappa}_{simp}$ indexed in $H$ below so that:

\begin{equation}
\label{conclhypothese} 
\begin{split}
&\sum_{l\in L_\mu} a_l Xdiv_{i_1}\dots Xdiv_{i_\mu} [C^{l,i_1\dots i_\mu}_g(\Omega_1,\dots,
\Omega_p,\phi_1,\dots,\phi_u)\cdot\Phi]=
\\&\sum_{h\in H} a_h Xdiv_{i_2}\dots Xdiv_{i_\mu}
[ C^{h,i_1\dots i_\mu}_g(\Omega_1,\dots,
\Omega_p,\phi_1,\dots,\phi_u)\nabla_{i_1}\Phi]
\\&+\sum_{j\in J} a_j Xdiv_{i_1}\dots Xdiv_{i_{\mu-1}}
[C^{j,i_1\dots i_{\mu-1}}_g(\Omega_1,\dots,
\Omega_p,\phi_1,\dots,\phi_u)\cdot \Phi],
\end{split}
\end{equation}
\end{lemma}
(modulo longer terms); here the terms indexed in $J$ are 
acceptable $(\mu-1)$-tensor fields in the form (\ref{form1}) 
which are simply subsequent to $\vec{\kappa}_{simp}$.

{\it Proof:} Pick out the
 factor $T_1=S_{*}R_{ijkl}\nabla^i\tilde{\phi}_1$. 
Let $L_\mu^A\subset L_\mu$ stand for the 
index set of terms in which contain a free index in this special factor and 
let $L_\mu^B\subset L_\mu$ stand for the index set of terms with no free index in 
that factor. We assume wlog that for each $l\in L_\mu^A$ the free index that 
belongs to the factor $T_1$ is ${}_{i_1}$. 

We will prove the following statements:

\begin{equation}
\label{thekill1} 
\begin{split}
&\sum_{l\in L_\mu^A} a_l Xdiv_{i_1}C^{l,i_1\dots i_\mu}_g\nabla_{i_2}\upsilon\dots\nabla_{i_\mu}\upsilon=
\sum_{l\in L'^B_\mu} a_l Xdiv_{i_1}C^{l,i_1\dots i_\mu}_g\nabla_{i_2}\upsilon\dots\nabla_{i_\mu}\upsilon
\\&+\sum_{t\in T} a_t C^{t,i_2\dots i_\mu}_g\nabla_{i_2}\upsilon\dots\nabla_{i_\mu}\upsilon+
\sum_{j\in J} a_j C^{j,i_2\dots i_\mu}_g\nabla_{i_2}\upsilon\dots\nabla_{i_\mu}\upsilon,
\end{split}
\end{equation}
where the tensor fields indexed in $L'^B_\mu$ are just like the tensor fields indexed 
in $L_\mu^B$, but the free index ${}_{i_1}$ does not belong to the factor $T_1$. 
The tensor fields indexed in $T$ are acceptable tensor fields of rank $\mu-1$, 
with a simple character $\vec{\kappa}_{simp}$, and moreover they have a 
factor $S_{*}\nabla R_{ijkl}$ (with one derivative)
which does not  contain a free index. 
\newline

We will prove (\ref{thekill1}) momentarily. Let us now check how it implies our claim:
We convert the factors $\nabla\upsilon$'s into $Xdiv$'s
(we are using the last Lemma in the Appendix of \cite{alexakis1} here), 
and replace into our Lemma hypothesis, to derive a new equation:

\begin{equation}
\label{heredia}
\begin{split}
&\sum_{t\in T} a_t Xdiv_{i_2}\dots Xdiv_{i_\mu} C^{t,i_2\dots i_\mu}_g+
\sum_{l\in L^B_\mu\bigcup L'^B_\mu}Xdiv_{i_1}\dots Xdiv_{i_\mu} C^{l,i_1\dots i_\mu}_g
\\&+\sum_{j\in J} a_j C^{j}_g.
\end{split}
\end{equation}
We next pick out the sublinear combination of terms in 
(\ref{heredia}) with a factor $S_{*}R_{ijkl}\nabla^i\tilde{\phi}_1$ (no derivatives);
 this sublinear combination clearly vanishes separately, so we derive:

\begin{equation}
\label{heredia2}
\begin{split}
\sum_{l\in L^B_\mu\bigcup L'^B_\mu}X_{*}div_{i_1}\dots X_{*}div_{i_\mu} C^{l,i_1\dots i_\mu}_g
+\sum_{j\in J} a_j C^{j}_g=0.
\end{split}
\end{equation}
(Here $X_{*}div_i[\dots]$ stands for the sublinear combination in $Xdiv_i[\dots]$ 
where $\nabla^i$ is not allowed to hit the factor $S_{*}R_{ijkl}\nabla^i\tilde{\phi}_1$.
We then define a formal operation $Op[\dots]$ which acts on the terms above by replacing the expression 
$S_{*}R_{ijkl}\nabla^i\tilde{\phi}_1$ by an expression
 $\nabla_j\omega\nabla_k\omega\nabla_l\upsilon-\nabla_j\omega\nabla_l\omega\nabla_k\upsilon$; 
denote the resulting $(u-1)$-simple character 
 (which keeps track of the indices 
$\nabla\tilde{\phi}_2,\dots,\nabla\tilde{\phi}_u$) by $\vec{\kappa}_{simp}'$.
Observe that this produces a new true equation:

\begin{equation}
\label{heredia2'}
\begin{split}
\sum_{l\in L^B_\mu\bigcup L'^B_\mu}Xdiv_{i_1}\dots Xdiv_{i_\mu} Op[C]^{l,i_1\dots i_\mu}_g
+\sum_{j\in J} a_j C^{j}_g=0,
\end{split}
\end{equation}
where the terms $C^j_g$ are simply subsequent to $\vec{\kappa}'_{simp}$. 
We can then apply the ``generalized version'' of Lemma \ref{addition} 
to the above (the above falls 
under the inductive assumption of (the generalized version of)
 Lemma \ref{addition} because the terms above 
have $\sigma_1+\sigma_2+p=\sigma-1$). We derive that:

\begin{equation}
\label{heredia2''}
\begin{split}
\sum_{l\in L^B_\mu\bigcup L'^B_\mu} Op[C]^{l,i_1\dots i_\mu}_g\nabla_{i_1}\upsilon\dots\nabla_{i_\mu}\upsilon=0.
\end{split}
\end{equation}
Now, we formaly replace each  expression $\nabla_a\omega\nabla_b\omega,\nabla_c\upsilon$ 
by an expression \\$S_{*}R_{i(ab)c}\nabla^i\tilde{\phi}_1$ and derive that:

\begin{equation}
\label{heredia2'''}
\begin{split}
\sum_{l\in L^B_\mu\bigcup L'^B_\mu} C^{l,i_1\dots i_\mu}_g\nabla_{i_1}\upsilon\dots\nabla_{i_\mu}\upsilon=0.
\end{split}
\end{equation}

Thus, we may return to (\ref{heredia}) and erase the sublinear combination
in $L^B_\mu\bigcup L'^B_\mu$. We then pick out the sublinear combination in that equation with 
a factor $S_{*}\nabla_aR_{ijkl}\nabla^i\tilde{\phi}_1$. We then derive  a new true equation:

\begin{equation}
\label{peteg}
\begin{split}
&\sum_{t\in T} a_t X_{*}div_{i_2}\dots X_{*}div_{i_\mu} C^{t,i_2\dots i_\mu}_g+
\sum_{j\in J} a_j C^{j}_g;
\end{split}
\end{equation}
($X_{*}div_i\dots $ now means that $\nabla^{i}$ is not allowed to hit the factor 
$S_{*}\nabla_aR_{ijkl}$). 
We then define a formal operation $Op'[\dots]$ which acts on 
the terms above by replacing the expression 
$S_{*}\nabla_aR_{ijkl}\nabla^{i_1}\tilde{\phi}_1$ by an expression
 $\nabla_a\omega\nabla_j\omega\nabla_k\omega\nabla_l\upsilon-
\nabla_a\omega\nabla_j\omega\nabla_l\omega\nabla_k\upsilon$;
  denote the resulting $(u-1)$-simple character 
 (which keeps track of the indices $\nabla\tilde{\phi}_2,\dots,\nabla\tilde{\phi}_u$)
  by $\vec{\kappa}_{simp}'$. Then by the same argument as above, 
  we derive that:

\begin{equation}
\label{peteg'}
\begin{split}
&\sum_{t\in T} a_t  Op'[C]^{t,i_2\dots i_\mu}_g\nabla_{i_2}
\upsilon\dots\nabla_{i_\mu}\upsilon=0,
\end{split}
\end{equation}
  and therefore:

\begin{equation}
\label{peteg''}
\begin{split}
&\sum_{t\in T} a_t  C^{t,i_2\dots i_\mu}_g\nabla_{i_2}
\upsilon\dots\nabla_{i_\mu}\upsilon=0.
\end{split}
\end{equation}

Thus, we derive our claim by just multiplying 
the equations (\ref{thekill1}), (\ref{heredia2'''}), 
(\ref{peteg''}) by $\Phi$, converting
 the $\nabla\upsilon$'s into $Xdiv$'s (we are here applying 
the relevant Lemma from the Appendix of 
\cite{alexakis1}),\footnote{Recall that in this setting
 the derivative $\nabla^i$ in each $Xdiv_i$ {\it is} 
allowed to hit the factor $\Phi$.} and then adding the resulting equations. 
\newline

Thus, matters are reduced to proving (\ref{thekill1}). We do this as follows: Refer to 
our Lemma assumption and pick out the sublinear combination of terms with a factor 
$S_{*}R_{ijkl}\nabla^i\tilde{\phi}_1$ (with no derivatives). This sublinear 
combination vanishes separaely, thus we derive a new true equation:

\begin{equation}
\label{nixon}
\begin{split}&\sum_{l\in L^A_1} a_lXdiv_{i_1}X_{*}div_{i_2}\dots 
X_{*}div_{i_\mu} C^{l,i_1\dots i_\mu}_g
+\sum_{l\in L^B_1} a_l X_{*}div_{i_1}\dots X_{*}div_{i_\mu} 
C^{l,i_1\dots i_\mu}_g\\&+\sum_{j\in J} a_j C^j_g=0.
\end{split}
\end{equation} 
Again, applying the operation $Op$ (defined above) to (\ref{nixon}) we derive a new true equation:

\begin{equation}
\label{nixon2}
\begin{split}
&\sum_{l\in L^A_1} a_lX_{*}div_{i_2}\dots X_{*}div_{i_\mu} 
\{Xdiv_{i_1}Op[C]^{l,i_1\dots i_\mu}_g\}
\\&+\sum_{l\in L^B_1} a_l X_{*}div_{i_1}\dots X_{*}div_{i_\mu} 
Op[C]^{l,i_1\dots i_\mu}_g+\sum_{j\in J} a_j C^j_g=0.
\end{split}
\end{equation} 
Here $Xdiv_{i_1}Op[C]^{l,i_1\dots i_\mu}_g$ stands for the sublinear combination where the 
derivative $\nabla^i$ is not allowed to hit any of the 
factors $\nabla\phi_h$ nor any of the factors $\nabla\omega,\nabla\upsilon$. 
In fact, if we treat the $Xdiv_iOp[C]^{l,i_1\dots i_\mu}_g$ as a sum 
of $(\mu-1)$-tensor fields (so we forget its $Xdiv$-structure). We then 
apply the inductive assumption of Lemma \ref{addition} to derive 
that there exists a linear combination of $\mu$-tensor fields
 with a $(u-1)$-simple character $\vec{\kappa}_{simp}'$, 
such that:

\begin{equation}
\label{nixon3}
\begin{split}&\sum_{l\in L^A_1} a_l \{Xdiv_{i_1}Op[C]^{l,i_1\dots i_\mu}_g\}
\nabla_{i_2}\upsilon\dots\nabla\upsilon_{i_\mu}
+\sum_{h\in H} a_l X_{*}div_{i_\mu}C^{h,i_1\dots i_\mu}_g\nabla_{i_2}
\upsilon\dots\nabla\upsilon_{i_\mu}\\&+\sum_{j\in J} a_j C^j_g=0.
\end{split}
\end{equation} 
Now, we act on the above by another formal operation 
$Op^{-1}[\dots]$ which replaces each expression
 $\nabla_a\omega\nabla_b\omega\nabla_c\upsilon$ by $S_{*}R_{i(ab)c}\nabla^i\tilde{\phi}_1$.
The result precisely is our desired (\ref{thekill1}).

\subsection{Mini-Appendix: Proof that  the ``delicate assumption'' (in case I)
can be made with no loss of generality:}
\label{ma1}
 We let $M$ stand for the number of free indices in the critical factor, 
for the terms of maximal refined double character in (\ref{assumpcion}). 
We denote by
$L_{\mu,*}\subset \bigcup_{z\in Z'_{Max}}L^z$ the index
set of $\mu$-tensor fields in (\ref{assumpcion}) with $M$ free
indices in the critical factor and with only one index (the index
${}_l$) in the critical factor contracting against another factor,
{\it in particular against a special index in some (simple) factor
$S_{*}\nabla^{(\rho)}R_{abcd}$}.\footnote{Denote this other factor by $T'$ (while the
critical factor will be denoted by $T_{*}$).}  We
will show that:

\begin{equation}
\label{trexeis?}
\begin{split}
&\sum_{l\in L_{\mu,*}} a_l C^{l,i_1\dots i_\mu}_g
 (\Omega_1,\dots,\Omega_p,\phi_1,\dots,\phi_u)\nabla_{i_1}\upsilon\dots\nabla_{i_\mu}\upsilon-
\\&\sum_{h\in H} a_h  Xdiv_{i_{\mu+1}}C^{h,i_1\dots i_{\mu+1}}_g
 (\Omega_1,\dots,\Omega_p,\phi_1,\dots,\phi_u)
\nabla_{i_1}\upsilon\dots\nabla_{i_\mu}\upsilon=
\\& \sum_{t\in T} a_t C^{t,i_1\dots i_\mu}_g
 (\Omega_1,\dots,\Omega_p,\phi_1,\dots,\phi_u)
\nabla_{i_1}\upsilon\dots\nabla_{i_\mu}\upsilon,
\end{split}
\end{equation}
where the tensor fields in the RHS have all the features of the
tensor fields in the first line but in addition the 
index ${}_l$ in the critical factor is contracting against a non-special index.
 If we can show (\ref{trexeis?}), it then
follows that the ``delicate assumption'' can be made with no loss
of generality.

{\it Proof of (\ref{trexeis?}):}  We divide the index set
$L_{\mu,*}$ according to {\it which} factor
$S_{*}\nabla^{(\rho)}R_{abcd}$ the index ${}_l$ in the critical
factor is contracting against: We say that $l\in L_{\mu,*,k}, k\in
K$ if and only if the index ${}_l$ is contracting against a
special index in the factor
$S_{*}\nabla^{(\rho)}R_{ijcd}\nabla^i\tilde{\phi}_k$ (denote this
factor by $T'_k$)--say the index ${}_l$ in $T'_k$.

We prove (\ref{trexeis?}) for the terms in $L_{\mu,*,k}$. Clearly,
if we can prove this then the whole of (\ref{trexeis?}) will
follow. We denote by $\tilde{C}^{l,i_1\dots
i_\mu}_g(\Omega_1,\dots,\Omega_p,Y_1,Y_2,\phi_1,\dots,\phi_u)$ the
tensor field that arises from $C^{l,i_1\dots
i_\mu}_g(\Omega_1,\dots,\Omega_p,\phi_2,\dots,\hat{\phi}_k,\dots,\phi_u)$ by
replacing the expression $S_{*}\nabla^{(\nu)}_{r_1\dots
r_\nu}R_{ijkl}S_{*}\nabla^{(\rho)}_{y_1\dots
y_\rho}{R_{i'j'k'}}^l\nabla^i\tilde{\phi}_1\nabla^{i'}\tilde{\phi}_k$
by \\$\nabla^{(\nu+2)}_{r_1\dots r_\nu
jk}Y_1\nabla^{(\rho+2)}_{y_1\dots y_\rho j'k'}Y_2$; denote the
resulting simple character by $Cut[\vec{\kappa}_{simp}]$.
Considering the second conformal variation of (\ref{assumpcion})
 and pick out the terms of length $\sigma+u$ with the 
factors $\nabla\tilde{\phi}_1,\nabla\tilde{\phi}_k$ contracting against each other,
we derive a new true equation:

\begin{equation}
\label{stateofterror}
\begin{split}
&[\sum_{l\in L_{\mu,*,k}} a_l Xdiv_{i_1}\dots
Xdiv_{i_\mu}\tilde{C}^{l,i_1\dots
i_\mu}_g(\Omega_1,\dots,\Omega_p,Y_1,Y_2,\phi_2,\dots,\hat{\phi}_k,\dots,\phi_u)+
\\&\sum_{h\in H} a_h Xdiv_{i_1}\dots Xdiv_{i_a} C^{h,i_1\dots
i_a}_g(\Omega_1,\dots,\Omega_p,Y_1,Y_2,\phi_2,\dots,\hat{\phi}_k,\dots,\phi_u)+
\\&\sum_{j\in J} a_j
C^j_g(\Omega_1,\dots,\Omega_p,Y_1,Y_2,\phi_2,\dots,\hat{\phi}_k,\dots,\phi_u)]
\nabla^s\tilde{\phi}_1\nabla_s\tilde{\phi}_k=0.
\end{split}
\end{equation}
The terms indexed in $H$ have length $>\mu$ and are acceptable
with a simple character $Cut[\vec{\kappa}_{simp}]$. The complete
contractions indexed in $J$ are simply subsequent to
$Cut[\vec{\kappa}_{simp}]$.

\par Now, we apply our inductive assumption of Corollary
\ref{corgiade} to the above,\footnote{Notice that
(\ref{stateofterror}) does not fall under any of the  ``forbidden
cases'', since the for the factor $\nabla^{(A)}Y$ we have  $\Phi_1+M\ge 2$, 
where $M$ is the number of free indices that belong 
to that factor and $\Phi_1$ is the number of factors 
$\nabla\phi_h$ that contract against it.} and we pick out the sublinear
combination of maximal refined double character--denote the index
set of those terms by $\tilde{L}_{\mu,*,k}$ (notice that the
sublinear combination $\sum_{l\in L_{\mu,*,k}}$ will be included
in those terms). We then derive that there exists a linear
combination of acceptable $(\mu+1)$-tensor fields with a simple
character $Cut[\vec{\kappa}_{simp}]$ so that:
\begin{equation}
\label{tore}
\begin{split}
&\sum_{l\in \tilde{L}_{\mu,*,k}} a_l \tilde{C}^{l,i_1\dots
i_\mu}_g(\Omega_1,\dots,\Omega_p,Y_1,Y_2,\phi_2,\dots,\hat{\phi}_k,\dots,\phi_u)\nabla_{i_1}
\upsilon\dots\nabla_{i_\mu}\upsilon-
\\&\sum_{h\in H} a_h Xdiv_{i_{\mu+1}} C^{h,i_1\dots
i_{\mu+1}}_g(\Omega_1,\dots,\Omega_p,Y_1,Y_2,\phi_2,\dots,\hat{\phi}_k,\dots,\phi_u)\nabla_{i_1}
\upsilon\dots\nabla_{i_\mu}\upsilon
\\&+\sum_{j\in J} a_j C^{j,i_1\dots
i_\mu}_g(\Omega_1,\dots,\Omega_p,Y_1,Y_2,\phi_2,\dots,\hat{\phi}_k,\dots,\phi_u)\nabla_{i_1}
\upsilon\dots\nabla_{i_\mu}\upsilon=0.
\end{split}
\end{equation}

\par Now, formally replacing the expression
$$\nabla^{(A)}_{y_1\dots y_A}Y_1\otimes\nabla^{(B)}_{r_1\dots r_B}Y_2
\otimes\nabla^{y_q}\phi_{z}\otimes\dots\otimes
\nabla^{y_{x-1}}\phi_{\chi}\otimes\nabla^{y_x}\upsilon
\dots\nabla^{y_A}\otimes\upsilon$$ by an expression
$$S_{*}\nabla^{(A-2)}_{y_1\dots y_{A-2}}R_{iy_{A-1}y_Al}\otimes
S_{*}\nabla^{(B-2)}_{r_1\dots
r_{B-2}}{R_{i'r_{B-1}r_B}}^l\otimes\nabla^i\phi_1\otimes\nabla^{i'}\phi_k$$
and repeating the permutations that make the above hold formally,\footnote{See the 
argument or the proof of (\ref{sarnak}).} we derive our claim. $\Box$

\end{document}